\documentclass{tran-l}

\usepackage{verbatim}
\usepackage{appendix}
\usepackage{amssymb}
\usepackage{xcolor}
\usepackage{xparse}

\newtheorem{theorem}{Theorem}[section]
\newtheorem{lemma}[theorem]{Lemma}
\newtheorem{proposition}[theorem]{Proposition}
\newtheorem{corollary}[theorem]{Corollary}

\theoremstyle{definition}
\newtheorem{definition}[theorem]{Definition}

\theoremstyle{remark}
\newtheorem{remark}[theorem]{Remark}

\numberwithin{equation}{section}

\newif\ifShowText
\ShowTexttrue   
\newcommand{\maybe}[1]{%
  \ifShowText%
    #1%
  \fi%
}

\newif\ifShowRed
\ShowRedfalse

\newcommand{\red}[1]{%
  \ifShowRed
    \textcolor{red}{#1}%
  \else
    #1%
  \fi
}

\newenvironment{redblock}
  {\par\ifShowRed\color{red}\fi}
  {\par}

\newif\ifpubinfo      \pubinfofalse
\makeatletter
\ifpubinfo
\else
  \let\ISSN\@empty
  \def\@serieslogo{}%
  \def\@setcopyright{}%
  \let\copyrightinfo\@gobbletwo
  \gdef\@volume{}%
  \gdef\@number{}%
  \gdef\@pages{}%
  \gdef\@date{}%
  \def\ps@firstpage{\ps@empty}%
\fi
\makeatother

\allowdisplaybreaks[1]

\begin{document}

\title{Entropy Rigidity for Maximal Representations}


\author{\maybe{Zhufeng Yao}}
\thanks{\maybe{The author was supported by the NUS-MOE grants  A-8000458-00-00 and A-8001950-00-00.}}
\curraddr{}
\email{}


\date{}

\dedicatory{}

\begin{abstract}
Let $\Gamma\subset \mathrm{PSL}(2,\mathbb{R})$ be a lattice and $\rho:\Gamma\to \mathrm{Sp}(2n,\mathbb{R})$ be a maximal representation with continuous limit map. We show that $\rho$ satisfies a measurable $(1,1,2)$-hypertransversality condition when $n\ge3$. With this we define a measurable Gromov product and the Bowen--Margulis--Sullivan measure associated to the unstable Jacobian introduced by Pozzetti, Sambarino and Wienhard. As a main application, we prove a strong entropy rigidity result for $\rho$.
\end{abstract}

\maketitle

\setcounter{tocdepth}{1}
\tableofcontents

\section{Introduction}
Let $S$ be a finite-type hyperbolic surface, $\mathrm{T}^1 S$ its unit tangent bundle, and $\phi_t$ the geodesic flow on $\mathrm{T}^1S$. Let $\mathcal{CG}$ denote the set of closed geodesics on $S$, and let $l(\gamma)$ for $\gamma \in \mathcal{CG}$ denote the length of $\gamma$. The \emph{topological entropy} $h^{\phi_t}$ of the flow can be computed via the \emph{length spectrum entropy} $h^l$ of the hyperbolic structure, that is,
\[
h^{\phi_t} = h^{l} := \lim_{T\to \infty} \frac{\log \# \{\gamma \mid \gamma \in \mathcal{CG},\ l(\gamma)\le T\} }{T},
\]
and Sullivan~\cite{sullivanoriginal} showed that this equals the Hausdorff dimension of the limit set of the Fuchsian model of $S$ (see also Patterson~\cite{patterson1976limit} and Schapira~\cite[Theorem 3.10]{schapira:hal-01882452}; see also Margulis~\cite{margulisthesis}). In particular, when $S$ has finite area, the limit set is a circle, so the identity $h^{\phi_t} = 1$ holds rigidly under deformation of the hyperbolic structures in the Teichmüller space.

When $S$ is a closed surface, Potrie--Sambarino~\cite{potrie2017eigenvaluesentropyhitchinrepresentation} studied the entropy rigidity phenomenon in higher Teichmüller spaces known as Hitchin components. These were first introduced by Hitchin~\cite{HITCHIN1992449} using the theory of Higgs bundles, and later developed by Labourie~\cite{labourie2005anosovflowssurfacegroups} from a dynamical perspective via Anosov flows, and by Fock and Goncharov~\cite{fock2006modulispaceslocalsystems} from an algebraic perspective using positivity. Roughly speaking, Hitchin components are the connected components of the representation space $\mathrm{Hom}(\pi_1(S), \mathrm{PSL}(d,\mathbb{R}))/\mathrm{PSL}(d,\mathbb{R})$ that contain the irreducibly Fuchsian representations, and conjugacy classes in these components are called Hitchin representations. In this setting, Potrie and Sambarino used thermodynamic formalism to show that when $d > 1$, the renormalized symmetric space orbit entropy—a suitable generalization of the topological entropy of the geodesic flow—is at most $1$, and that equality is attained if and only if the Hitchin representation is irreducibly Fuchsian. Thus, in the higher rank setting, entropy is no longer invariant under deformation, but it is rigid in the sense that only Fuchsian representations achieve maximal entropy. For the case when $S$ is not closed, see Canary--Zhang--Zimmer~\cite{CaZhZi}.

On the other hand, building on earlier work of Goldman~\cite{Goldman1988}, Toledo~\cite{toledo89}, Hern\'andez~\cite{Hernandez} and others, Burger--Iozzi--Wienhard~\cite{maximalearlier,burger2008surfacegrouprepresentationsmaximal} and Burger--Iozzi--Labourie--Wienhard~\cite{maximal} systematically developed the theory of \emph{maximal representations} of $\pi_1(S)$ into Hermitian Lie groups, most notably those of tube type such as $\mathrm{Sp}(2n,\mathbb{R})$ (see also Bradlow--Garc\'ia-Prada--Gothen~\cite{bgg1,bgg2} and others via Higgs bundles). These are representations whose Toledo invariant attains its maximal absolute value. They share many similar properties with Hitchin representations, such as the positivity of the associated limit map. Our main goal in this paper is to establish an entropy rigidity result analogous to that of Potrie--Sambarino~\cite{potrie2017eigenvaluesentropyhitchinrepresentation}, but for maximal representations into $\mathrm{Sp}(2n,\mathbb{R})$.

We state our results as follows:

\subsection{Entropy rigidity}
Let $\mathbb{D}$ denote the 2-dimensional hyperbolic disk, $\Gamma \subset \mathrm{PSL}(2,\mathbb{R}) = \mathrm{Isom}^+(\mathbb{D})$ a non-elementary discrete subgroup, and $\Lambda(\Gamma) \subset \partial \mathbb{D} =S^1$ its limit set. We assume that $\Gamma$ is a lattice in most main results of this paper, i.e., $\Gamma$ is finitely generated and $\Lambda(\Gamma) = \partial \mathbb{D}$. 

Let $\mathbb{R}^{2n}$ be the Euclidean space equipped with the standard inner product and an orthonormal basis $\{e_1,\dots,e_{2n}\}$. The matrix
\(
\Omega = \begin{pmatrix} 0 & \mathrm{Id}_n \\ -\mathrm{Id}_n & 0 \end{pmatrix}
\)
represents a symplectic form. The corresponding symplectic group is $\mathrm{Sp}(2n,\mathbb{R}) = \{ g \in \mathrm{SL}(2n,\mathbb{R}) \mid g^T \Omega g = \Omega \}$. For any $g\in \mathrm{Sp}(2n,\mathbb{R})$, we define the \emph{Cartan projection} $\kappa(g)$ as
\[
\kappa(g) = \mathrm{diag}(\log \sigma_1(g), \dots, \log \sigma_n(g), -\log \sigma_1(g), \dots, -\log \sigma_n(g)),
\]
where $\sigma_1(g) \ge \cdots \ge \sigma_n(g) \ge 1 \ge \sigma_n(g)^{-1} \ge \cdots \ge \sigma_1(g)^{-1}$ are the singular values of $g$. We also define $\alpha(\kappa(g)) := 2\log \sigma_n(g)$ and $\hat{\omega}_{\alpha}(\kappa(g)) := \frac{2}{n}(\log \sigma_1(g)+\dots+\log \sigma_n(g))$.

Let $\mathcal{F}_\alpha$ be the flag manifold consisting of $n$-dimensional Lagrangian in $(\mathbb{R}^{2n},\Omega)$. Since  $\partial{\mathbb{D}} = S^1$, we can define positive $N$-tuples on the boundary via cyclic order. On the other hand, $\mathcal{F}_\alpha$ can be identified with the Shilov boundary of the Siegel upper half-plane, and positive $N$-tuples can be defined there as well (see Section~\ref{posandmax}). 

A representation $\rho: \Gamma \to \mathrm{Sp}(2n,\mathbb{R})$ is called \emph{maximal} in this paper if there exists a continuous $\rho$-equivariant map $\xi^\alpha: \Lambda(\Gamma) \to \mathcal{F}_\alpha$ that sends positive tuples to positive tuples. (In some contexts, maximal representations are defined without requiring a continuous boundary map. For example, for a finite-area surface with cusps, a representation with maximal relative Toledo invariant need not admit such a continuous map, see \cite[discussions before Lemma~8.8]{burger2008surfacegrouprepresentationsmaximal}.) When $\Gamma$ is geometrically finite, maximal representations are relatively $\alpha$-Anosov with respect to the cusp subgroups of $\Gamma$ (see Burger--Iozzi--Labourie--Wienhard~\cite{maximal}, Guichard--Labourie--Wienhard~\cite{GLW}, and Yao~\cite[Proposition~1.1]{yzf}), and $\xi^{\alpha}$ serves as the limit map, which is strongly dynamics-preserving (see Canary--Zhang--Zimmer~\cite[Theorem~1.2 and Appendix~B]{CaZhZi2}). For any such $\rho$ and any fixed base point $b_0\in\mathbb{D}$, there exist constants $a>1$ and $A\ge0$ such that
\begin{equation}\label{anosovgrowth}
   \frac{1}{a} d_{\mathbb{D}}(b_0, \gamma(b_0)) - A \le \alpha(\kappa(\rho(\gamma))) \le a d_{\mathbb{D}}(b_0, \gamma(b_0)) + A. 
\end{equation}
If we define the $\alpha$- and $\hat{\omega}_{\alpha}$-Poincaré series of $\rho$ by
\[
Q_{\rho}^{\alpha}(s) = \sum_{\gamma \in \Gamma} e^{-s \alpha(\kappa(\rho(\gamma)))},\quad Q_{\rho}^{\hat{\omega}_{\alpha}}(s) = \sum_{\gamma \in \Gamma} e^{-s \hat{\omega}_{\alpha}(\kappa(\rho(\gamma)))}
\]
and define the \emph{critical exponents} by
\[
\delta_{\rho}(\alpha) = \inf\left\{s \in \mathbb{R}^{\ge 0}  \,\middle|\, Q_{\rho}^{\alpha}(s) < \infty \right\},\quad \delta_{\rho}(\hat{\omega}_{\alpha}) = \inf\left\{s \in \mathbb{R}^{\ge 0}  \,\middle|\, Q_{\rho}^{\hat{\omega}_{\alpha}}(s) < \infty \right\},
\]
then Inequality~(\ref{anosovgrowth}) implies that $\delta_{\rho}(\alpha)$ is finite.

It can be further shown that if $\Gamma$ is a lattice, then $\delta_{\rho}(\alpha) = 1$ (see Pozzetti--Sambarino--Wienhard~\cite[Corollary 9.6]{Liplim} and Yao~\cite[Theorem 1.2]{yzf}). Since $\hat{\omega}_\alpha(\kappa(g)) \ge \alpha(\kappa(g))$ holds for any $g\in \mathrm{Sp}(2n,\mathbb{R})$, it follows from the definition of critical exponents that $\delta_{\rho}(\hat{\omega}_\alpha) \le 1$.
\begin{redblock}

The equality can be attained. Define the diagonal embedding $\rho_D: \mathrm{SL}(2,\mathbb{R}) \to \mathrm{Sp}(2n,\mathbb{R})$ by
\[
\begin{pmatrix} a & b \\ c & d \end{pmatrix}
\mapsto
\begin{pmatrix} a \mathrm{Id}_n & b \mathrm{Id}_n \\ c \mathrm{Id}_n & d \mathrm{Id}_n \end{pmatrix}.
\]
Observe that the $\mathrm{Sp}(2n,\mathbb R)$-normalizer $N_{D}$ of the image subgroup is generated by
\[
\mathrm{diag}(Q,Q), Q\in \mathrm{O}(n) \quad \text{and} \quad \rho_D(\mathrm{SL}(2,\mathbb{R})).
\]
Then $\delta_{\rho}(\hat{\omega}_\alpha) = 1$ holds whenever $\rho(\Gamma)\subset N_D$. Indeed, it is direct to check that for any $g\in N_D$,
\[
\alpha(\kappa(g)) = \hat{\omega}_{\alpha}(\kappa(g)).
\]
So $\delta_{\rho}(\hat{\omega}_\alpha) = \delta_{\rho}(\alpha)=1$.

Our main result states that this is the only case in which equality holds.

\begin{theorem}\label{main}
If $\rho: \Gamma \to \mathrm{Sp}(2n,\mathbb{R})$ is a maximal representation from a lattice $\Gamma$, then $\delta_{\rho}(\hat{\omega}_\alpha) \le 1$. Equality holds if and only if $\rho(\Gamma)$ can be conjugated into $N_{D}$.
\end{theorem}   
\end{redblock}

In the case where $\Gamma$ is torsion-free and $S = \mathbb{D}/\Gamma$ is a closed surface, let $[\Gamma]_{\mathrm{hyp}}$ denote the set of conjugacy classes of non-identity elements in $\Gamma$. Note that each such class consists purely of hyperbolic elements. Recall that $\mathrm{T}^1 S$ denotes the unit tangent bundle and $\phi_t$ denotes the geodesic flow on it. There is a one-to-one correspondence between elements of $[\Gamma]_{\mathrm{hyp}}$ and closed geodesics on $S$. Using thermodynamic formalism, one can reparameterize $\phi_t$ so that each closed geodesic corresponding to $[\gamma] \in [\Gamma]_{\mathrm{hyp}}$ has length $\hat{\omega}_\alpha(\lambda(\rho(\gamma)))$, where $\hat{\omega}_\alpha(\lambda(\rho(\gamma))) := \frac{2}{n}\sum_{i= 1}^{n} \log \lambda_i(\rho(\gamma))$ is the renormalized sum of the logarithms of the first $n$ eigenvalue moduli of $\rho(\gamma)$. If we define the length entropy $h_{\rho}(\hat{\omega}_\alpha)$ by
\[
h_{\rho}(\hat{\omega}_\alpha) := \limsup_{T \to \infty} \frac{\log \# R_T^{\hat{\omega}_\alpha}(\Gamma)}{T} \text{, where } R_T^{\hat{\omega}_\alpha}(\Gamma) := \{[\gamma] \in [\Gamma]_{\mathrm{hyp}} \mid \hat{\omega}_\alpha(\lambda(\rho(\gamma))) \le T\},
\]
then it turns out that $h_{\rho}(\hat{\omega}_\alpha)$ coincides with the topological entropy of this reparameterized flow (see Potrie--Sambarino~\cite[Section 2]{potrie2017eigenvaluesentropyhitchinrepresentation}). On the other hand, the relative $\alpha$-Anosov property of maximal representations implies the equality $h_{\rho}(\hat{\omega}_\alpha) = \delta_{\rho}(\hat{\omega}_\alpha)$ (see Canary--Zhang--Zimmer~\cite[Proposition~9.1]{CaZhZi}). Hence, rigidity of the critical exponent for $\hat{\omega}_\alpha$ is equivalent to rigidity of the topological entropy for this reparameterized flow.

\subsection{The symmetric space exponent}

\begin{redblock}
We can also interpret Theorem \ref{main} more geometrically as an orbit growth rigidity in the associated symmetric space.    
\end{redblock}

The maximal compact subgroup of $\mathrm{Sp}(2n,\mathbb{R})$ can be chosen as $\mathrm{U}(n) = \mathrm{O}(2n)\cap \mathrm{Sp}(2n,\mathbb{R})$ and the associated symmetric space $X = \mathrm{Sp}(2n,\mathbb{R})/\mathrm{U}(n)$ is a Hermitian symmetric space which can be identified as the \emph{Siegel upper half plane}. We choose a left invariant metric $d_{X}(-,-)$ compatible with the Killing form $B$ on the Lie algebra $\mathfrak{sp}(2n,\mathbb{R})$ where
\begin{redblock}
$$
B(\kappa(g),\kappa(g)) = 4(n+1)\sum_{i = 1}^{n} \log^2(\sigma_i(g))
$$
and
$$
d_{X}(\mathrm{U}(n),g\mathrm{U}(n)):= \frac{1}{\sqrt{n(n+1)}}\sqrt{B(\kappa(g),\kappa(g))}.
$$
Note that $d_{X}$ is normalized so that $\rho_{D}$ induces an isometric embedding from $\mathbb{D}$ to $X$. We denote this embedded disk by $\mathbb D_{D}$, and define the $\mathrm{Sp}(2n,\mathbb R)$ translates of $\mathbb D_D$ as the \emph{diagonal disks}.

Define the \emph{symmetric space exponent} by
$$
\delta_{X}(\rho) = \inf \{s\in \mathbb{R}^{\ge0}\mid \sum_{\gamma\in \Gamma} e^{-sd_{X}(\mathrm{U}(n),\rho(\gamma)\mathrm{U}(n))} <\infty\}.
$$
Since $B(\kappa(g),\kappa(g))\ge 4(n+1) \log^2(\sigma_n(g)) = (n+1) (\alpha(\kappa(g)))^2$ for any $g\in \mathrm{Sp}(2n,\mathbb{R})$, inequality~\eqref{anosovgrowth} implies that $\delta_{X}(\rho)$ is finite.

\begin{theorem}\label{symrig}
If $\rho:\Gamma\to \mathrm{Sp}(2n,\mathbb{R})$ is a maximal representation from a lattice, then $\delta_{X}(\rho)\le 1$. The equality is attained if and only if $\rho(\Gamma)$ preserves a diagonal disk.
\end{theorem}

\begin{proof}
Note that the $\mathrm{Sp}(2n,\mathbb{R})$ stabilizer of $\mathbb{D}_{D}$ coincides with $N_D$.

We first prove the "only if" direction. By the Cauchy--Schwarz inequality, we have
\begin{align*}
d_{X}(\mathrm{U}(n),\rho(\gamma) \mathrm{U}(n))
&= \frac{1}{\sqrt{n(n+1)}}\sqrt{4(n+1)\sum_{i =1}^{n}(\log \sigma_i(\rho(\gamma)))^2} \\
&\ge \frac{2}{n}\left(\sum_{i = 1}^{n} \log \sigma_i(\rho(\gamma))\right) = \hat{\omega}_\alpha(\kappa(\rho(\gamma))).
\end{align*}
Thus, $\delta_X(\rho)\le \delta_{\rho}(\hat{\omega}_{\alpha})$. If $\delta_{X}(\rho) = 1$, then $\delta_{\rho}(\hat{\omega}_{\alpha}) \ge 1$. Theorem~\ref{main} then implies that $\delta_{\rho}(\hat{\omega}_{\alpha}) = 1$ and that $\rho(\Gamma)$ can be conjugated into $N_D$. So $\rho(\Gamma)$ preserves a diagonal disk.

Conversely, if $\rho(\Gamma)$ preserves a diagonal disk, then up to conjugation, we may assume without loss of generality that it preserves $\mathbb{D}_{D}$, thus $\rho(\Gamma)\subset N_D$. In this case, it can be verified that
\[
\delta_{X}(\rho) = \delta_{\rho}(\hat{\omega}_{\alpha}) = \delta_{\rho}(\alpha) = 1.
\]
\end{proof}

\subsection{The Manhattan curve theorem}

A small application of Theorem \ref{main} is to prove a weak version of the Manhattan curve theorem (see Burger~\cite{burgermanhattan}; see also Bishop--Steger~\cite{bishopsteger}).

Fix a torsion-free lattice $\Gamma\subset \mathrm{PSL}(2,\mathbb{R})$ and suppose that there are two type-preserving Fuchsian representations $\rho_i:\Gamma\to \mathrm{PSL}(2,\mathbb{R}),\, i = 1,2$. Define $l_i(\gamma) := 2\log \sigma(\rho_i(\gamma))$, where $\sigma(g) \ge 1$ denotes the largest singular value of any lift $g \in \mathrm{SL}(2,\mathbb{R})$. Let $\delta_i$ be the critical exponent of $l_i$, and let $\delta'$ be the critical exponent of the averaged function $\frac{l_1 + l_2}{2}$. It is known that $\delta_1 = \delta_2 = 1$ (see Sullivan~\cite{sullivanoriginal}), and since the map $s \mapsto e^{-s}$ is convex, the definition of the critical exponent implies that $\delta' \le 1$.

\begin{corollary}
    We have $\delta' = 1$ if and only if $\rho_1$ and $\rho_2$ are conjugate in $\mathrm{PGL}(2,\mathbb R)$.
\end{corollary}
\begin{proof}
We only need to prove the ``only if" direction. 

By applying a $\begin{pmatrix} 1&0\\ 0&-1\end{pmatrix}$ conjugation if necessary, we assume that the boundary maps of both $\rho_1$ and $\rho_2$ have the same orientation.

Since $\Gamma$ is torsion-free, we lift $\rho_i$ to $\tilde{\rho}_i:\Gamma\to \mathrm{SL}(2,\mathbb R)$, $i = 1,2$. Write 
\[
\tilde{\rho}_i(\gamma) = \begin{pmatrix}
    a_i(\gamma) & b_i(\gamma) \\
    c_i(\gamma) & d_i(\gamma)
\end{pmatrix},
\]
and define the diagonal representation
\[
\rho_{1,2} : \Gamma \to \mathrm{Sp}(4,\mathbb{R})
\]
by
\[
\gamma \mapsto \begin{pmatrix}
    a_1(\gamma) & 0 & b_1(\gamma) & 0 \\
    0 & a_2(\gamma) & 0 & b_2(\gamma) \\
    c_1(\gamma) & 0 & d_1(\gamma) & 0 \\
    0 & c_2(\gamma) & 0 & d_2(\gamma)
\end{pmatrix}
\]
which is maximal. Notice that
\[
\hat\omega_{\alpha}(\kappa(\rho_{1,2}(\gamma))) = \frac{l_1(\gamma) + l_2(\gamma)}{2}.
\]
By Theorem \ref{main}, $\delta' = 1$ if and only if $\rho_{1,2}(\Gamma)$ can be conjugated into $N_{D}$, which implies the two largest eigenvalue moduli coincide. So $\mathbb{D}/\rho_1(\Gamma)$ and $\mathbb{D}/\rho_2(\Gamma)$ have the same marked length spectrum, which implies that $\rho_1$ and $\rho_2$ are conjugate (see, for example, Imayoshi--Taniguchi~\cite[Theorem~3.12]{imayoshi}).
\end{proof}

\end{redblock}

\subsection{Historical comments}
Wienhard~\cite[Task~22]{wienhardtalk} raised the problem of finding upper bounds for the topological entropy of maximal representations and geometrically characterizing the representations that attain these bounds. For maximal representations into $\mathrm{SO}^{\circ}(2,n)$ (which has real rank~2), Collier--Tholozan--Toulisse~\cite[Corollary~5]{rk2maximal} have already shown a strong entropy rigidity result. Indeed, their result is much stronger (see \cite[Theorem~4]{rk2maximal}): representations outside the Fuchsian locus have length spectra dominating that of some Fuchsian representation. For $\Theta$-positive representations into $\mathrm{SO}(p,q)$ with $q>p>1$, Pozzetti--Sambarino--Wienhard~\cite[Theorem~10.7]{Liplim} showed a strong rigidity result. Although $\mathrm{SO}(p,q)$ has arbitrarily high real rank as $p$ and $q$ vary, its $\Theta$-positive structure still arises from rank~$2$ Hermitian tube-type Lie groups (see Guichard–Wienhard~\cite[Proposition~3.8 and Section~5.1]{GW2}). Indeed, for simple Lie groups, the crucial difference between a $\Theta$-positive structure (not of split real type) and a Hitchin structure (i.e., a $\Theta$-positive structure of split real type) arises from the unique special root in the Dynkin diagram. This unique special root always arises from an associated Hermitian tube-type Lie group; for $\mathrm{SO}(p,q)$, this associated group is $\mathrm{SO}(2,q-p+2)$. However, in our paper, $\mathrm{Sp}(2n,\mathbb{R})$ is a Hermitian tube-type Lie group of arbitrarily high rank as $n$ varies.

It remains an interesting problem how to prove entropy rigidity results for all $\Theta$-positive representations (see \cite{GW2}).

\textbf{Note.}
About three months after the first version of this paper was made public,
Kim and Zimmer posted a preprint~\cite{kim2026measurableboundarymapspattersonsullivan} in which they use measurable boundary
maps to study Patterson--Sullivan measures on the Furstenberg boundary. With some additional justification, one can apply their Theorem~10.1 to show that
if $\delta_{\rho}(\hat{\omega}_{\alpha})=1$, then
$\alpha(\lambda(\rho(\gamma)))=\hat{\omega}_{\alpha}(\lambda(\rho(\gamma)))$
for every $\gamma\in\Gamma$, where $\lambda$ denotes the Jordan projection.
Combining this with the final argument in the proof of Theorem~\ref{main}
gives an alternative proof of our main result.

However, this paper also establishes a new geometric property of maximal
representations, namely the measurable $(1,1,2)$-hypertransversality
introduced in Section~\ref{sec3}. It may offer new perspectives on
the study of other properties of maximal representations.

\subsection{Outline of the paper}
In Section \ref{sec2} we discuss some preliminaries on $\mathrm{Sp}(2n,\mathbb{R})$, the definition of positivity, and maximal representations.

In Section \ref{sec3} we show that the limit map $\xi^{\alpha}$ satisfies a measurable version of Pozzetti--Sambarino--Wienhard's $(1,1,2)$-hypertransversality (see~\cite{pozzetti2020conformalityrobustclassnonconformal}).

In Section \ref{sec4} we recall the definition of $\alpha$-Patterson--Sullivan measure $\mu_{\alpha}$, which was introduced by Pozzetti--Sambarino--Wienhard~\cite[Section 6]{Liplim} using unstable Jacobians. We then use the measurable $(1,1,2)$-hypertransversality in Section~\ref{sec3} to construct a measurable Gromov product, and produce the corresponding Bowen--Margulis--Sullivan measure $m_{\alpha}$. Following Canary--Zhang--Zimmer~\cite{CaZhZi3,CaZhZi-relatively}, we also introduce the $\hat{\omega}_\alpha$-Patterson--Sullivan measure $\mu_{\hat{\omega}_\alpha}$ and the Bowen--Margulis--Sullivan measure $m_{\hat{\omega}_\alpha}$. Finally, we present the shadow lemmas for $\mu_{\alpha}$ and $\mu_{\hat{\omega}_\alpha}$, which are crucial in the proof of Theorem~\ref{main}.

In Section \ref{mainproof} we prove Theorem \ref{main}; the method follows Canary--Zhang--Zimmer~\cite[Proof of Theorem~13.1]{CaZhZi3} (see also Sullivan~\cite[Section~8]{sullivanacta} and Burger~\cite[Proof of Theorem~1(a)]{burgermanhattan}). The assumption $\delta_{\rho}(\hat{\omega}_{\alpha}) = 1$, together with the shadow lemma, implies that $m_{\hat{\omega}_\alpha}$ is absolutely continuous with respect to $m_{\alpha}$. Moreover, this assumption also ensures that the ergodicity of $m_{\hat{\omega}_\alpha}$ transfers to $m_{\alpha}$, so the two measures are equivalent up to a constant multiple. A further calculation gives a quantitative estimate on the Radon--Nikodym derivative of the corresponding Patterson--Sullivan measures. We then apply the shadow lemma again to show that the Jordan projection is restricted to a one-dimensional subspace, so the Benoist limit cone of $\rho(\Gamma)$ is one-dimensional. It follows that the semisimple part of the Zariski closure has real rank one, which concludes the proof.

\vspace{3mm}

\textbf{Acknowledgement.} The author sincerely thanks his advisor\maybe{, Prof.~Tengren Zhang,} for consistent guidance, many helpful discussions, and for reading drafts and providing comments that improved the author's writing skills.

The author also sincerely thanks Prof.~Beatrice Pozzetti for her discussions, kindness, and encouragement, and Tianqi Wang for his useful comments on this paper.

\section{Notation and preliminaries}\label{sec2}

\subsection{Lattice}

Let $\mathbb{D}$ denote the hyperbolic disk; its orientation–preserving isometry group is $\mathrm{PSL}(2,\mathbb{R})$. Let $\Gamma \subset \mathrm{PSL}(2,\mathbb{R})$ be a discrete subgroup. Fix a base point $b_0 \in \mathbb{D}$ and define the limit set
\[
\Lambda(\Gamma) \subset \partial\mathbb{D}
\]
to be the set of accumulation points of the orbit $\Gamma\cdot b_0$. This set is independent of the choice of $b_0$. If $\#\Lambda(\Gamma)>2$, then $\Lambda(\Gamma)$ is a compact perfect set; in this case we say that $\Gamma$ is \emph{non-elementary}.

We call $\Gamma$ a \emph{lattice} if the hyperbolic surface $\mathbb{D}/\Gamma$ has finite hyperbolic area. It is well known that $\Gamma$ is a lattice if and only if it is finitely generated and its limit set equals the ideal boundary of $\mathbb{D}$, i.e.\ $\Lambda(\Gamma)=\partial\mathbb{D}$.

Since $\partial \mathbb{D} = S^1$, a tuple $(x_1, x_2, \dots, x_N)$ of points on the boundary is said to be \emph{positive} if the points are cyclically ordered on the circle.

\subsection{Basics of $\mathrm{Sp}(2n,\mathbb{R})$}\label{section:notationsforsp}

Throughout this paper, we denote $\mathsf{c}_g$ the $g$-conjugate action on the Lie group and denote $\mathrm{Ad}_g$ the $g$-adjoint action on the Lie algebra. For any matrix $M$, let $M^{T}$ denote its transpose. We choose an orthonormal basis $\{e_1,e_2,\dots,e_{2n}\}$ in $\mathbb{R}^{2n}$ and define the symplectic form $\omega(v,w) = v^{T} \Omega w$ where \(
\Omega = \begin{pmatrix} 0 & \mathrm{Id}_n \\ -\mathrm{Id}_n & 0 \end{pmatrix}
\). The Lie algebra of $\mathrm{Sp}(2n,\mathbb{R})$ can be written as
\[
\mathfrak{sp}(2n,\mathbb{R}) := \left\{\begin{pmatrix}A&B\\C&D \end{pmatrix} : A+D^{T} = 0,\ B = B^{T},\ C = C^{T} \right\}.
\]

The maximal compact subgroup $K \subset \mathrm{Sp}(2n,\mathbb{R})$ can be taken as
\[
\mathrm{U}(n) = \mathrm{O}(2n,\mathbb{R}) \cap \mathrm{Sp}(2n,\mathbb{R}).
\]

Choose the Cartan subalgebra
\[
\mathfrak{a}
:=\Bigl\{\mathrm{diag}(\lambda_1,\dots,\lambda_n,-\lambda_1,\dots,-\lambda_n)\ :\ \lambda_i\in\mathbb{R}\Bigr\}
\subset \mathfrak{sp}(2n,\mathbb{R}).
\]
For $i=1,\dots,n-1$, define linear forms on $\mathfrak{a}$ by
\[
\beta_i\!\left(\mathrm{diag}(\lambda_1,\dots,\lambda_n,-\lambda_1,\dots,-\lambda_n)\right)=\lambda_i-\lambda_{i+1},
\qquad
\alpha\!\left(\mathrm{diag}(\lambda_1,\dots,\lambda_n,-\lambda_1,\dots,-\lambda_n)\right)=2\lambda_n.
\]
Thus $\Delta=\{\alpha,\beta_1,\dots,\beta_{n-1}\}$ is the set of simple positive roots.
The corresponding positive Weyl chamber is
\[
\mathfrak{a}^{+}
=\Bigl\{\mathrm{diag}(\lambda_1,\dots,\lambda_n,-\lambda_1,\dots,-\lambda_n)\in\mathfrak{a}\ :\
\lambda_1>\lambda_2>\cdots>\lambda_n>0\Bigr\},
\]
whose closure is denoted by $\overline{\mathfrak{a}}^{+}$.

We will call $\alpha$ the \emph{long root}. Its fundamental weight $\omega_\alpha \in \mathfrak{a}^*$ is
\[
\omega_\alpha(\mathrm{diag}(\lambda_1, ..., \lambda_n, -\lambda_1, ..., -\lambda_n)) = \lambda_1 + \cdots + \lambda_n,
\]
and we denote by $\hat{\omega}_\alpha := \frac{2}{n} \omega_\alpha$ the \emph{renormalized fundamental weight}. Note that $\hat{\omega}_\alpha \ge \alpha$ on $\mathfrak{a}^+$.

Define the \emph{Cartan projection} $\kappa : \mathrm{Sp}(2n,\mathbb{R}) \to \overline{\mathfrak{a}^+}$ by
\[
\kappa(g) = \mathrm{diag}(\log \sigma_1(g), ..., \log \sigma_n(g), -\log \sigma_1(g), ..., -\log \sigma_n(g)),
\]
where $\sigma_1(g) \ge \cdots \ge \sigma_n(g) \ge 1 \ge \sigma_n(g)^{-1} \ge \cdots \ge \sigma_1(g)^{-1}$ are the singular values of $g$. Define the \emph{Jordan projection} $\lambda : \mathrm{Sp}(2n,\mathbb{R}) \to \overline{\mathfrak{a}^+}$ by
\[
\lambda(g) = \mathrm{diag}(\log \lambda_1(g), ..., \log \lambda_n(g), -\log \lambda_1(g), ..., -\log \lambda_n(g)),
\]
where $\lambda_1(g) \ge \cdots \ge \lambda_n(g) \ge 1 \ge \lambda_n(g)^{-1} \ge \cdots \ge \lambda_1(g)^{-1}$ are the moduli of the eigenvalues of $g$.

Let $\mathrm{Sym}(n)$ denote the set of $n\times n$ real symmetric matrices, and let $\mathrm{Pos}(n)$ denote the set of $n\times n$ positive-definite symmetric matrices. Note that $\mathrm{Pos}(n)$ is an open, acute, convex cone in $\mathrm{Sym}(n)$ (here “acute” means that the closure of $\mathrm{Pos}(n)$ contains no nontrivial linear subspace).

Define the parabolic subgroup $\mathsf{P}_{\alpha}$ and its opposite $\mathsf{P}_{\alpha}^{\mathrm{opp}}$ associated to $\alpha$ by
\[
\mathsf{P}_{\alpha} := \left\{ \begin{pmatrix}A&B\\0&A^{-T} \end{pmatrix} : A \in \mathrm{GL}(n,\mathbb R),\ B A^{T} = A B^{T} \right\},
\]
\[
\mathsf{P}_{\alpha}^{\mathrm{opp}} := \left\{ \begin{pmatrix}A&0\\C&A^{-T} \end{pmatrix} : A \in \mathrm{GL}(n,\mathbb R),\ A^{T} C = C^{T} A \right\}.
\]

Their Lie algebras are
\[
\mathfrak{p}_{\alpha} := \left\{ \begin{pmatrix}A&M\\0&-A^{T} \end{pmatrix} : A \in \mathfrak{gl}(n,\mathbb R),\ M \in \mathrm{Sym}(n) \right\},
\]
\[
\mathfrak{p}_{\alpha}^{\mathrm{opp}} := \left\{ \begin{pmatrix}A&0\\M&-A^{T} \end{pmatrix} : A \in \mathfrak{gl}(n,\mathbb R),\ M \in \mathrm{Sym}(n) \right\}.
\]

The Levi subgroup $\mathsf{L}_{\alpha}$ associated to $\alpha$ is defined as
\[
\mathsf{L}_{\alpha} := \left\{ \begin{pmatrix}A&0\\0&A^{-T} \end{pmatrix} : A \in \mathrm{GL}(n,\mathbb R) \right\},
\]
with Lie algebra
\[
\mathfrak{l}_{\alpha} := \left\{ \begin{pmatrix}A&0\\0&-A^{T} \end{pmatrix} : A \in \mathfrak{gl}(n,\mathbb R) \right\}.
\]

The unipotent subgroup $\mathsf{U}_{\alpha}$ and its opposite $\mathsf{U}_{\alpha}^{\mathrm{opp}}$ are given by
\[
\mathsf{U}_{\alpha} := \left\{ \begin{pmatrix}\mathrm{Id}_n&M\\0&\mathrm{Id}_n \end{pmatrix} : M \in \mathrm{Sym}(n) \right\}, \quad
\mathsf{U}_{\alpha}^{\mathrm{opp}} := \left\{ \begin{pmatrix}\mathrm{Id}_n&0\\M&\mathrm{Id}_n \end{pmatrix} : M \in \mathrm{Sym}(n) \right\},
\]
with Lie algebras
\[
\mathfrak{u}_{\alpha} := \left\{ \begin{pmatrix}0&M\\0&0 \end{pmatrix} : M \in \mathrm{Sym}(n) \right\}, \quad
\mathfrak{u}_{\alpha}^{\mathrm{opp}} := \left\{ \begin{pmatrix}0&0\\M&0 \end{pmatrix} : M \in \mathrm{Sym}(n) \right\}.
\]

Finally, define
\[
\mathsf{U}_{\alpha}^{>0} := \left\{ \begin{pmatrix}\mathrm{Id}_n&M\\0&\mathrm{Id}_n \end{pmatrix} : M \in \mathrm{Pos}(n) \right\}, \quad
\mathsf{U}_{\alpha}^{\mathrm{opp},>0} := \left\{ \begin{pmatrix}\mathrm{Id}_n&0\\M&\mathrm{Id}_n \end{pmatrix} : M \in \mathrm{Pos}(n) \right\}.
\]

In the remainder of this paper, for any $M \in \mathrm{Sym}(n)$, we will use the following notation:
\[
U^M = \begin{pmatrix} \mathrm{Id}_n & M \\ 0 & \mathrm{Id}_n \end{pmatrix} \in \mathsf{U}_{\alpha}, \quad U_M = \begin{pmatrix} \mathrm{Id}_n & 0 \\ M & \mathrm{Id}_n \end{pmatrix} \in \mathsf{U}_{\alpha}^{\mathrm{opp}},
\]
\[
u^M = \begin{pmatrix} 0 & M \\ 0 & 0 \end{pmatrix} \in \mathfrak{u}_{\alpha}, \quad u_M = \begin{pmatrix} 0 & 0 \\ M & 0 \end{pmatrix} \in \mathfrak{u}_{\alpha}^{\mathrm{opp}}.
\]

Note the following properties:
\begin{enumerate}
\item \(
\mathsf{P}_{\alpha}
= \{ g \in \mathrm{Sp}(2n,\mathbb{R}) : \mathrm{Ad}_g(\mathfrak{p}_{\alpha}) = \mathfrak{p}_{\alpha} \}
= \{ g \in \mathrm{Sp}(2n,\mathbb{R}) : \mathrm{Ad}_g(\mathfrak{u}_{\alpha}) = \mathfrak{u}_{\alpha} \}.
\)
Similarly,
\(
\mathsf{P}_{\alpha}^{\mathrm{opp}}
= \{ g \in \mathrm{Sp}(2n,\mathbb{R}) : \mathrm{Ad}_g(\mathfrak{p}_{\alpha}^{\mathrm{opp}}) = \mathfrak{p}_{\alpha}^{\mathrm{opp}} \}
= \{ g \in \mathrm{Sp}(2n,\mathbb{R}) : \mathrm{Ad}_g(\mathfrak{u}_{\alpha}^{\mathrm{opp}}) = \mathfrak{u}_{\alpha}^{\mathrm{opp}} \}.
\)

    \item $\mathsf{L}_{\alpha} = \mathsf{P}_{\alpha} \cap \mathsf{P}_{\alpha}^{\mathrm{opp}}$, and $\mathsf{c}_{\mathsf{L}_{\alpha}}$ preserves both $\mathsf{U}_{\alpha}^{>0}$ and $\mathsf{U}_{\alpha}^{\mathrm{opp},>0}$.
    \item There are semidirect product decompositions: $\mathsf{P}_{\alpha} = \mathsf{L}_{\alpha} \mathsf{U}_{\alpha}$, \quad $\mathsf{P}_{\alpha}^{\mathrm{opp}} = \mathsf{L}_{\alpha} \mathsf{U}_{\alpha}^{\mathrm{opp}}$.
\end{enumerate}

\subsection{Flag manifold}\label{flag}
Recall 
\[
\Omega=\begin{pmatrix}
0&\mathrm{Id}_n\\ -\mathrm{Id}_n&0
\end{pmatrix}\in \mathrm{Sp}(2n,\mathbb{R}).
\]
For fixed $n\in\mathbb{N}$, define the $\alpha$-flag manifold by the homogeneous space
\[
\mathcal{F}_{\alpha}=\mathrm{Sp}(2n,\mathbb{R})/\mathsf{P}_{\alpha}.
\]
\subsubsection{Various identifications of $\mathcal{F}_{\alpha}$}\label{identificationsofflags}

First note that $\mathcal{F}_{\alpha}$ can be identified with the manifold of Lagrangian subspaces, i.e., $n$-dimensional subspaces of $\mathbb{R}^{2n}$ on which the symplectic form $\omega$ vanishes identically. Let $\mathcal{L}=\mathrm{Span}(e_1,\ldots,e_n)$ and $\mathcal{L}^{\mathrm{opp}}=\mathrm{Span}(e_{n+1},\ldots,e_{2n})$ be the standard pair of  Lagrangian subspaces. Then for each $g\in \mathrm{Sp}(2n,\mathbb{R})$ we naturally identify $g\mathsf{P}_{\alpha}$ with $g\mathcal{L}$. Note that $\mathcal{L}^{\mathrm{opp}} = \Omega\mathcal{L}$.

We also note that, since
\[
\mathsf{P}_{\alpha}=\{\,g\in \mathrm{Sp}(2n,\mathbb{R}) \mid \mathrm{Ad}_g(\mathfrak{p}_{\alpha})=\mathfrak{p}_{\alpha}\,\},
\]
the space $\mathcal{F}_{\alpha}$ can be identified with the $\mathrm{Sp}(2n,\mathbb{R})$-conjugacy class of $\mathfrak{p}_{\alpha}$ in $\mathfrak{sp}(2n,\mathbb{R})$. Thus for each $g\in \mathrm{Sp}(2n,\mathbb{R})$ we also identify $g\mathsf{P}_{\alpha}$ with $\mathrm{Ad}_g(\mathfrak{p}_{\alpha})$. Also note that $\mathfrak{p}_{\alpha}^{\mathrm{opp}} = \mathrm{Ad}_\Omega (\mathfrak{p}_{\alpha})$.

Let
\[
\Lambda^n:\ \mathrm{SL}(2n,\mathbb{R}) \to \mathrm{SL}\!\left(\bigwedge^{n}\mathbb{R}^{2n}\right)
\]
be the $n$-th wedge representation, and define the bilinear form $\omega'$ on $\bigwedge^{n}\mathbb{R}^{2n}$ by
\[
\omega'(v,w)=\frac{v\wedge w}{e_1\wedge e_2\wedge \cdots \wedge e_{2n}}.
\]
It has to be noted that $v\wedge w$ is an element in $\bigwedge^{2n} \mathbb{R}^{2n}$ but not in $\bigwedge^2(\bigwedge^{n}\mathbb{R}^{2n})$. And this form is nondegenerate and symmetric or skew-symmetric depending on whether $n$ is even or odd, respectively.

Let $V_0=\mathrm{Span}(e_1\wedge \cdots \wedge e_n)$, and let $[V_0]$ denote the corresponding projective point in $\mathbb{P}(\bigwedge^n \mathbb{R}^{2n})$. We then observe that
\[
\mathsf{P}_{\alpha}
=\{\, g \in \mathrm{Sp}(2n,\mathbb{R}) \mid \Lambda^n(g)[V_0]=[V_0] \,\},
\]
so the flag manifold $\mathcal{F}_{\alpha}$ can be identified with the orbit
\[
\mathcal{F}_{\alpha}\cong \Lambda^n\bigl(\mathrm{Sp}(2n,\mathbb{R})\bigr)\!\cdot [V_0]
\subset \mathbb{P}\!\left(\bigwedge^n \mathbb{R}^{2n}\right),
\]
which is a closed submanifold. This identification defines the \emph{Plücker embedding}
\[
\mathcal{P}_{\alpha}:\ \mathcal{F}_{\alpha}\to \mathbb{P}\!\left(\bigwedge^{n}\mathbb{R}^{2n}\right),
\qquad
\mathcal{P}_{\alpha}(g\mathsf{P}_{\alpha})=\Lambda^n(g)[V_0].
\]

As $\omega'$ is nondegenerate, we have a similar dual embedding:

\[
\mathcal{P}_{\alpha}^{\mathrm{opp}}:\ \mathcal{F}_{\alpha}\to
\mathbb{P}\!\left((\bigwedge^n \mathbb{R}^{2n})^*\right),\qquad
\mathcal{P}_{\alpha}^{\mathrm{opp}}(g\mathsf{P}_{\alpha})
=\Lambda^n(g)\,\mathrm{Ann}([V_0]).
\]
Here $\mathrm{Ann}([V_0])$ denotes the annihilator hyperplane of $[V_0]$ with respect to $\omega'$, i.e., $\mathrm{Ann}([V_0]): = \{w\in \bigwedge^n \mathbb{R}^{2n}\mid \omega'(V_0,w) = 0\}$ (note that $\omega'$ is symmetric or skew-symmetric so we do not distinguish left or right annihilator), and $\mathbb{P}\!\left((\bigwedge^n \mathbb{R}^{2n})^*\right)$ is the dual projective space, which can be identified as the space of projective hyperplanes in $\bigwedge^n \mathbb{R}^{2n}$. Note that $V_0\subset \mathrm{Ann}([V_0])$, so for any $x\in\mathcal{F}_{\alpha}$, we have $\mathcal{P}_{\alpha}(x)\subset \mathcal{P}_{\alpha}^{\mathrm{opp}}(x)$ when viewed as linear subspaces.

\subsubsection{The automorphism group}

To talk about positivity, we need a Lie group acting more transitively on the space of triples in $\mathcal{F}_{\alpha}$, so we introduce $\mathrm{Aut}(\mathfrak{sp}(2n,\mathbb{R}))$, which denotes the automorphism group of the Lie algebra $\mathfrak{sp}(2n,\mathbb{R})$. 

$\mathrm{Aut}(\mathfrak{sp}(2n,\mathbb{R}))$ acts on $\mathcal{F}_{\alpha}$ via the identification with the $\mathrm{Sp}(2n,\mathbb{R})$-conjugacy class of $\mathfrak{p}_{\alpha}$. Since the Dynkin diagram of $\mathfrak{sp}(2n,\mathbb{R})$ has no nontrivial automorphisms, for any $g\in \mathrm{Aut}(\mathfrak{sp}(2n,\mathbb{R}))$ the subalgebra $g(\mathfrak{p}_{\alpha})$ is $\mathrm{Sp}(2n,\mathbb{R})$ conjugated to $\mathfrak{p}_{\alpha}$ (see Guichard--Wienhard~\cite[Section~5.3]{GW2}).

By transporting this action to the identification of $\mathcal{F}_{\alpha}$ with the space of Lagrangian subspaces, we also obtain an action of $\mathrm{Aut}(\mathfrak{sp}(2n,\mathbb{R}))$ on $\mathcal{F}_{\alpha}$ in the Lagrangian model. For any $g\in \mathrm{Aut}(\mathfrak{sp}(2n,\mathbb{R}))$, we use the following shorthand to identify these actions:
\[
g(\mathfrak{p}_{\alpha})=g\mathfrak{p}_{\alpha}=g\mathcal{L},\qquad
g(\mathfrak{p}_{\alpha}^{\mathrm{opp}})=g\mathfrak{p}_{\alpha}^{\mathrm{opp}}=g\mathcal{L}^{\mathrm{opp}}.
\]
And for $g\in \mathrm{Sp}(2n,\mathbb{R})$, we identify the following notation:
\[
g\mathsf{P}_{\alpha}
=\mathrm{Ad}_g(\mathfrak{p}_{\alpha})
=\mathrm{Ad}_g\mathfrak{p}_{\alpha}
= g\mathfrak{p}_{\alpha}
= g\mathcal{L}
= \Lambda^n(g)[V_0]
\]
and
\[
\mathrm{Ad}_g(\mathfrak{p}_{\alpha}^{\mathrm{opp}})
=\mathrm{Ad}_g\mathfrak{p}_{\alpha}^{\mathrm{opp}}
= g\mathfrak{p}_{\alpha}^{\mathrm{opp}} = g\mathcal{L}^{\mathrm{opp}}.
\]

We have some useful calculations for actions on the flag manifold. Let
\[
\mathcal{N} =
\begin{pmatrix}
0 & \mathrm{Id}_n \\
\mathrm{Id}_n & 0
\end{pmatrix}
\in \mathrm{GL}(2n,\mathbb{R}).
\]
It is straightforward to verify that, when we regard $\mathfrak{sp}(2n,\mathbb{R}) \subset \mathfrak{gl}(2n,\mathbb{R})$, the automorphism $\mathrm{Ad}_{\mathcal{N}}$ preserves $\mathfrak{sp}(2n,\mathbb{R})$, and hence $\mathrm{Ad}_{\mathcal{N}}|_{\mathfrak{sp}(2n,\mathbb{R})} \in \mathrm{Aut}(\mathfrak{sp}(2n,\mathbb{R}))$.  
And recall that we let
\[
\Omega =
\begin{pmatrix}
0 & \mathrm{Id}_n \\
- \mathrm{Id}_n & 0
\end{pmatrix}
\in \mathrm{Sp}(2n,\mathbb{R}).
\]
The following lemma can be verified directly, and its computations will be used frequently throughout the paper.

\begin{lemma}\label{lemma:computations}
Following the notation in Section~\ref{section:notationsforsp}, we have:
\begin{enumerate}
    \item 
    For any $M\in \mathrm{Sym}(n)$,
    \[
    \mathsf{c}_{\mathcal{N}}(U^{M}) = U_M, 
    \quad 
    \mathsf{c}_{\mathcal{N}}(U_{M}) = U^M,
    \qquad
    \mathsf{c}_{\Omega}(U^{M}) = U_{-M}, 
    \quad 
    \mathsf{c}_{\Omega}(U_{M}) = U^{-M}.
    \]

    \item 
    \[
    \mathsf{c}_{\mathcal{N}}(\mathsf{U}_{\alpha}^{>0}) 
    = \mathsf{U}_{\alpha}^{\mathrm{opp},>0},
    \quad
    \mathsf{c}_{\mathcal{N}}(\mathsf{U}_{\alpha}^{\mathrm{opp},>0}) 
    = \mathsf{U}_{\alpha}^{>0},
    \]
    \[
    \mathsf{c}_{\Omega}(\mathsf{U}_{\alpha}^{>0}) 
    = (\mathsf{U}_{\alpha}^{\mathrm{opp},>0})^{-1},
    \quad
    \mathsf{c}_{\Omega}(\mathsf{U}_{\alpha}^{\mathrm{opp},>0}) 
    = (\mathsf{U}_{\alpha}^{>0})^{-1}.
    \]

    \item 
    \[
    \mathrm{Ad}_{\mathcal{N}}\mathcal{L} 
    = \mathrm{Ad}_{\Omega}\mathcal{L} 
    = \mathcal{L}^{\mathrm{opp}},
    \qquad
    \mathrm{Ad}_{\mathcal{N}}\mathcal{L}^{\mathrm{opp}} 
    = \mathrm{Ad}_{\Omega}\mathcal{L}^{\mathrm{opp}} 
    = \mathcal{L}.
    \]

    \item 
\[
\forall\, M \in \mathrm{Sym}(n) \text{ invertible}, \qquad 
 U_{M}\mathcal{L}=U^{M^{-1}}\mathcal{L}^{\mathrm{opp}}.
\]

\end{enumerate}
\end{lemma}

\subsubsection{General flag varieties}

We now discuss flag manifolds in general linear spaces. Let $V$ be a finite-dimensional real vector space and fix a sequence of integers $1 \leq i_1 < i_2 < \cdots < i_k < \dim V$. Let $I = \{i_1, i_2, \ldots, i_k\}$. The associated flag manifold is
\[
\mathcal{F}_I(V)
:= \bigl\{(V_{i_1},\dots,V_{i_k}) \mid V_{i_1}\subset V_{i_2}\subset \cdots \subset V_{i_k}\subset V,\ \dim V_{i_j}=i_j \text{ for each } j\bigr\}.
\]
For each flag $F\in \mathcal{F}_I(V)$, let $F_{i_s}$ denote the subspace of dimension $i_s$.

 Let $I^{\mathrm{opp}} = \{\dim V-i_k,\dots,\dim V-i_1\}$.
Given $(F_1,F_2)\in \mathcal{F}_I(V)\times\mathcal{F}_{I^{\mathrm{opp}}}(V)$, we say that $(F_1,F_2)$ is \emph{transverse} if, for every $i_s\in I$, the subspaces $F_{1,i_s}$ and $F_{2,\dim V - i_s}$ are transverse in $V$, i.e., \ $F_{1,i_s}\oplus F_{2,\dim V - i_s}=V$. And we say that the index set $I$ is \emph{symmetric} if $I = I^{\mathrm{opp}}$ (and in this setting we can talk about transverse pairs in $\mathcal{F}_I(V)\times \mathcal{F}_I(V)$).

Now suppose that $(x,y)\in \mathcal{F}_{\alpha}\times \mathcal{F}_{\alpha}$, we say they are \emph{transverse} if they satisfy any of the following equivalent conditions under the various identifications of $\mathcal{F}_{\alpha}$:
\begin{enumerate}
    \item $x$ and $y$ are transverse as $n$-dimensional subspaces of $\mathbb{R}^{2n}$.
    \item $(\mathcal{P}_\alpha(x), \mathcal{P}_{\alpha}^{\mathrm{opp}}(x))$ and $(\mathcal{P}_\alpha(y), \mathcal{P}_{\alpha}^{\mathrm{opp}}(y))$ are transverse in $\mathcal{F}_{1,d-1}(\bigwedge^n \mathbb{R}^{2n})$, where $d=\dim \bigwedge^n \mathbb{R}^{2n}$.
    \item There exists $g\in \mathrm{Sp}(2n,\mathbb{R})$ such that $(x,y)=g(\mathfrak{p}_{\alpha},\mathfrak{p}_{\alpha}^{\mathrm{opp}})$ (equivalently, under the other identification, $(x,y)=g(\mathcal{L},\mathcal{L}^{\mathrm{opp}})$).
\end{enumerate}

\subsection{Positivity in $\mathcal{F}_{\alpha}$}\label{posandmax}

We now define positive tuples in $\mathcal{F}_{\alpha}$.

\begin{definition}[{Guichard--Wienhard~\cite[Definition~13.23 and Lemma~13.25]{GW2}}]\label{postuple}

An $N$-tuple $(x_1, x_2, \ldots, x_N) \subset \mathcal{F}_{\alpha}$, $N\ge 3$, is \emph{positive} if there exist $g \in \mathrm{Aut}(\mathfrak{sp}(2n, \mathbb{R}))$ and elements $u_1, u_2, \ldots, u_{N-2} \in \mathsf{U}_{\alpha}^{>0}$ such that
\[
(x_1, \ldots, x_N) = g\big(\mathcal{L}, u_{1} u_{2} \cdots u_{N-2} \mathcal{L}^{\mathrm{opp}}, \ldots,  u_1u_2 \mathcal{L}^{\mathrm{opp}}, u_1 \mathcal{L}^{\mathrm{opp}}, \mathcal{L}^{\mathrm{opp}}\big).
\]
\end{definition}

Here are properties of positive $N$-tuples.

\begin{proposition}[{Guichard--Wienhard~\cite[Lemma~13.24]{GW2}}]\label{proposition:propertyofpositivetuple}
    \begin{enumerate}
        \item The set of positive $N$-tuples is invariant under the dihedral group $D_N$.
        \item Any $N'$-subtuple of a positive $N$–tuple ($3\le N'\le N$) is still a positive $N'$–tuple.
    \end{enumerate}
\end{proposition}

We state some equivalent characterizations of positive $N$–tuples.

\begin{lemma}\label{lemma:equivalentdefinitionffpositivity}
    Let $(x_1, \ldots, x_N) \subset \mathcal{F}_{\alpha}$ be an $N$–tuple ($N \ge 3$). The following are equivalent:
    \begin{enumerate}
        \item There exist $g \in \mathrm{Aut}(\mathfrak{sp}(2n, \mathbb{R}))$ and elements $u_1, u_2, \ldots, u_{N-2} \in \mathsf{U}_{\alpha}^{>0}$ such that
        \[
        (x_1, \ldots, x_N)
        = g\big(\mathcal{L}, u_{1}u_{2}\cdots u_{N-2}\mathcal{L}^{\mathrm{opp}}, \ldots, u_1u_2\mathcal{L}^{\mathrm{opp}}, u_1\mathcal{L}^{\mathrm{opp}}, \mathcal{L}^{\mathrm{opp}}\big).
        \]

        \item There exist $g \in \mathrm{Aut}(\mathfrak{sp}(2n, \mathbb{R}))$ and elements $u_1, u_2, \ldots, u_{N-2} \in \mathsf{U}_{\alpha}^{\mathrm{opp},>0}$ such that
        \[
        (x_1, \ldots, x_N)
        = g\big(\mathcal{L}^{\mathrm{opp}}, u_1u_2\cdots u_{N-2}\mathcal{L}, \ldots, u_1u_2\mathcal{L}, u_1\mathcal{L}, \mathcal{L}\big).
        \]

        \item There exist $g \in \mathrm{Aut}(\mathfrak{sp}(2n, \mathbb{R}))$ and elements $u_1, u_2, \ldots, u_{N-2} \in \mathsf{U}_{\alpha}^{>0}$ such that
        \[
        (x_1, \ldots, x_N)
        = g\big(\mathcal{L}^{\mathrm{opp}}, u_1\mathcal{L}^{\mathrm{opp}}, u_1u_2\mathcal{L}^{\mathrm{opp}}, \ldots, u_1u_2\cdots u_{N-2}\mathcal{L}^{\mathrm{opp}}, \mathcal{L}\big).
        \]

        \item There exist $g \in \mathrm{Aut}(\mathfrak{sp}(2n, \mathbb{R}))$ and elements $u_1, u_2, \ldots, u_{N-2} \in \mathsf{U}_{\alpha}^{\mathrm{opp},>0}$ such that
        \[
        (x_1, \ldots, x_N)
        = g\big(\mathcal{L}, u_1\mathcal{L}, u_1u_2\mathcal{L}, \ldots, u_1u_2\cdots u_{N-2}\mathcal{L}, \mathcal{L}^{\mathrm{opp}}\big).
        \]
    \end{enumerate}

Moreover, in all the above statements, we may replace $\mathsf{U}_{\alpha}^{>0}$ (respectively, $\mathsf{U}_{\alpha}^{\mathrm{opp},>0}$) by $(\mathsf{U}_{\alpha}^{>0})^{-1}$ (respectively, $(\mathsf{U}_{\alpha}^{\mathrm{opp},>0})^{-1}$).
\end{lemma}

\begin{proof}
\noindent
$(1)\Rightarrow(2)$:  
If $g \in \mathrm{Aut}(\mathfrak{sp}(2n,\mathbb{R}))$ and $u_i \in \mathsf{U}_{\alpha}^{>0}$ for $i = 1, \ldots, N-2$, then Lemma~\ref{lemma:computations}~(1) shows that $u_i' := \mathcal{N}u_i\mathcal{N} \in \mathsf{U}_{\alpha}^{\mathrm{opp},>0}$.  
Thus
\[
(x_1, \ldots, x_N)
= g\big(\mathcal{L}, u_1u_2\cdots u_{N-2}\mathcal{L}^{\mathrm{opp}}, \ldots, u_1u_2\mathcal{L}^{\mathrm{opp}}, u_1\mathcal{L}^{\mathrm{opp}}, \mathcal{L}^{\mathrm{opp}}\big)
\]
\[
= g\big(\mathrm{Ad}_{\mathcal{N}}\mathcal{L}^{\mathrm{opp}}, \mathrm{Ad}_{\mathcal{N}}\cdot u'_1u'_2\cdots u'_{N-2}\cdot \mathrm{Ad}_{\mathcal{N}}\mathcal{L}^{\mathrm{opp}}, \ldots, \mathrm{Ad}_{\mathcal{N}}\cdot u'_1\cdot \mathrm{Ad}_{\mathcal{N}}\mathcal{L}^{\mathrm{opp}}, \mathrm{Ad}_{\mathcal{N}}\mathcal{L}\big)
\]
\[
= (g\mathrm{Ad}_{\mathcal{N}})\big(\mathcal{L}^{\mathrm{opp}}, u'_1u'_2\cdots u'_{N-2}\mathcal{L}, \ldots, u'_1u'_2\mathcal{L}, u'_1\mathcal{L}, \mathcal{L}\big).
\]

\noindent
$(2)\Rightarrow(3)$:  
If $g \in \mathrm{Aut}(\mathfrak{sp}(2n, \mathbb{R}))$ and $u_1, \ldots, u_{N-2} \in \mathsf{U}_{\alpha}^{\mathrm{opp},>0}$, write $u_i = U_{M_i}$ with $M_i \in \mathrm{Pos}(n)$, and set $P_i = (M_1 + \cdots + M_i)^{-1} \in \mathrm{Pos}(n)$.  
Then Lemma~\ref{lemma:computations}~(4) gives
\[
g\big(\mathcal{L}^{\mathrm{opp}}, u_1u_2\cdots u_{N-2}\mathcal{L}, \ldots, u_1u_2\mathcal{L}, u_1\mathcal{L}, \mathcal{L}\big)
= g\big(\mathcal{L}^{\mathrm{opp}}, U^{P_{N-2}}\mathcal{L}^{\mathrm{opp}}, \ldots, U^{P_2}\mathcal{L}^{\mathrm{opp}}, U^{P_1}\mathcal{L}^{\mathrm{opp}}, \mathcal{L}\big).
\]
Since $P_i - P_{i+1} \in \mathrm{Pos}(n)$ for $1 \le i \le N-3$, letting $u'_1 = U^{P_{N-2}} \in \mathsf{U}_{\alpha}^{>0}$ and $u'_i = U^{P_{N-i-1} - P_{N-i}} \in \mathsf{U}_{\alpha}^{>0}$ for $2 \le i \le N-2$, the equality continues to get
\[
g\big(\mathcal{L}^{\mathrm{opp}}, u'_1\mathcal{L}^{\mathrm{opp}}, u'_1u'_2\mathcal{L}^{\mathrm{opp}}, \ldots, u'_1u'_2\cdots u'_{N-2}\mathcal{L}^{\mathrm{opp}}, \mathcal{L}\big).
\]

\noindent
$(3)\Rightarrow(4)$:  
Similar to $(1)\Rightarrow(2)$.

\noindent
$(4)\Rightarrow(1)$:  
Similar to $(2)\Rightarrow(3)$.

\smallskip
The “Moreover” part follows directly by conjugating $\Omega$ to all the statements in $(1)$–$(4)$.  
Using $(1)$ as an example: if $g \in \mathrm{Aut}(\mathfrak{sp}(2n,\mathbb{R}))$ and $u_i \in \mathsf{U}_{\alpha}^{>0}$ for $i = 1, \ldots, N-2$, define $u_i' = \Omega u_i \Omega^{-1} \in (\mathsf{U}_{\alpha}^{\mathrm{opp},>0})^{-1}$. Then
\[
g\big(\mathcal{L}, u_1u_2\cdots u_{N-2}\mathcal{L}^{\mathrm{opp}}, \ldots, u_1u_2\mathcal{L}^{\mathrm{opp}}, u_1\mathcal{L}^{\mathrm{opp}}, \mathcal{L}^{\mathrm{opp}}\big)
\]
\[
= g\big(\mathrm{Ad}_{\Omega^{-1}}\mathcal{L}^{\mathrm{opp}}, \mathrm{Ad}_{\Omega^{-1}}\cdot u'_1u'_2\cdots u'_{N-2}\cdot \mathrm{Ad}_{\Omega}\mathcal{L}^{\mathrm{opp}}, \ldots, \mathrm{Ad}_{\Omega^{-1}}\cdot u'_1\cdot \mathrm{Ad}_{\Omega}\mathcal{L}^{\mathrm{opp}}, \mathrm{Ad}_{\Omega^{-1}}\mathcal{L}\big)
\]
\[
= g\mathrm{Ad}_{\Omega^{-1}}\big(\mathcal{L}^{\mathrm{opp}}, u'_1u'_2\cdots u'_{N-2}\mathcal{L}, \ldots, u'_1\mathcal{L}, \mathcal{L}\big),
\]
which yields the variant of (2) in which $\mathsf{U}_{\alpha}^{\mathrm{opp},>0}$ is replaced by $(\mathsf{U}_{\alpha}^{\mathrm{opp},>0})^{-1}$.
\end{proof}

\begin{lemma}[{Guichard--Wienhard~\cite[Proposition~13.1]{GW2}}]\label{lemma:changeofcoordinates}
Let $\mathcal{TF} \subset \mathcal{F}_{\alpha}$ be the subset of flags simultaneously transverse to $\mathcal{L}$ and $\mathcal{L}^{\mathrm{opp}}$.  
The map $\mathsf{U}_{\alpha}^{>0} \to \mathcal{F}_{\alpha}$, $u \mapsto u\mathcal{L}^{\mathrm{opp}}$, is a diffeomorphism onto a connected component of $\mathcal{TF}$.  
Similarly, $\mathsf{U}_{\alpha}^{\mathrm{opp},>0} \to \mathcal{F}_{\alpha}$, $u \mapsto u\mathcal{L}$, is also a diffeomorphism onto a connected component of $\mathcal{TF}$, and Lemma~\ref{lemma:computations}~(4) gives a diffeomorphism
\[
\mathsf{U}_{\alpha}^{>0}\mathcal{L}^{\mathrm{opp}} \longrightarrow \mathsf{U}_{\alpha}^{\mathrm{opp},>0}\mathcal{L}.
\]
\end{lemma}

Let $\Lambda\subset S^1$ be a closed subset such that $\#\Lambda\ge 3$. We say that a map $\xi^{\alpha} : \Lambda\to \mathcal{F}_{\alpha}$ is \emph{positive} if it sends every positive $N$–tuple in $\Lambda$ to a positive $N$–tuple in $\mathcal{F}_{\alpha}$. According to Guichard--Wienhard~\cite[Lemma 13.30]{GW2}, under the assumption that $\#\Lambda > 3$, this is equivalent to requiring that all positive quadruples are mapped to positive quadruples.

Here is a parameterization lemma for positive maps.
\begin{lemma}[{Burger--Iozzi--Labourie--Wienhard~\cite[Section~8, Proof of Corollary~6.3]{maximal}}]\label{lemma:parameterizationofmaps_combined}
Assume $\xi:\Lambda \to \mathcal{F}_{\alpha}$ is a continuous positive map. Let $x \neq y \in \Lambda$ be distinct points, and let $I_1, I_2 \subset S^1$ be two closed, disjoint arcs containing $x$ and $y$, respectively. Let $E_1 \subset I_1 \cap \Lambda$ and $E_2 \subset I_2 \cap \Lambda$ be closed subsets such that $E_1$ contains $x$ and $E_2$ contains $y$.

Then there exist $g \in \mathrm{Aut}(\mathfrak{sp}(2n,\mathbb{R}))$ and continuous maps 
\[
f_1: E_1 \to \mathrm{Pos}(n) \cup \{0\} \cup -\mathrm{Pos}(n), 
\qquad 
f_2: E_2 \to \mathrm{Pos}(n) \cup \{0\} \cup -\mathrm{Pos}(n),
\]
such that the following hold:
\begin{enumerate}
    \item $f_1(x) = 0$. For any $z \in E_1$,
    \[
    \xi(x,z,y) = g\big(\mathcal{L}^{\mathrm{opp}},\, U^{f_1(z)}\mathcal{L}^{\mathrm{opp}},\, \mathcal{L}\big),
    \]
    and whenever $y, z_1, z_2$ are clockwise ordered (where $z_1, z_2 \in E_1$), we have 
    \[
    f_1(z_2) - f_1(z_1) \in \mathrm{Pos}(n).
    \]
    
    \item $f_2(y)=0$. For any $z \in E_2$,
    \[
    \xi(x,z,y) = g\big(\mathcal{L}^{\mathrm{opp}},\, U_{f_2(z)}\mathcal{L},\, \mathcal{L}\big),
    \]
    and whenever $x, z_1, z_2$ are clockwise ordered (where $z_1, z_2 \in E_2$), we have 
    \[
    f_2(z_1) - f_2(z_2) \in \mathrm{Pos}(n).
    \]
\end{enumerate}
\end{lemma}

\begin{proof}
\textbf{Part (1):}
Pick $z_0\in E_1$ distinct from $x$. Since $\xi$ is a positive map, we can choose $g \in \mathrm{Aut}(\mathfrak{sp}(2n,\mathbb{R}))$ such that $\xi(x) = g\mathcal{L}^{\mathrm{opp}}$, $\xi(y) = g\mathcal{L}$, and $\xi(z_0) = gU^{M_0}\mathcal{L}^{\mathrm{opp}}$, where $M_0\in \mathrm{Pos}(n)$ if $(y,x,z_0)$ is clockwise, and $M_0\in -\mathrm{Pos}(n)$ if $(y,z_0,x)$ is clockwise.

Basic Lie theory shows that the map $\mathsf{U}_{\alpha} \to \mathcal{F}_{\alpha}$, $u \mapsto u\mathcal{L}^{\mathrm{opp}}$, is a diffeomorphism onto the set of flags transverse to $\mathcal{L}$. Since $\xi(E_1)$ lies in the set of flags transverse to $\xi(y) = g\mathcal{L}$, there exists a map $f_1: E_1 \to \mathrm{Sym}(n)$ such that for any $z \in E_1$,
\[
\xi(x, z, y) = g\big(\mathcal{L}^{\mathrm{opp}},\, U^{f_1(z)}\mathcal{L}^{\mathrm{opp}},\, \mathcal{L}\big).
\]
Clearly, $f_1(x) = 0$ and $f_1(z_0) = M_0$.

Now assume $y, z_1, z_2$ are clockwise ordered with $z_1, z_2 \in E_1$. To verify the positivity of the difference $f_1(z_2) - f_1(z_1)$, we distinguish cases based on the cyclic order of the quadruple $(x, y, z_1, z_2)$.
We specifically prove the case where the quadruple $(z_2, x, z_1, y)$ is counter-clockwise ordered; other cases follow a similar logic.

By the definition of positive maps, there exist $g_1 \in \mathrm{Aut}(\mathfrak{sp}(2n,\mathbb{R}))$ and $u_1, u_2 \in \mathsf{U}_{\alpha}^{>0}$ such that
\[
\xi(z_2, x, z_1, y) = g_1\big(\mathcal{L}^{\mathrm{opp}},\, u_1\mathcal{L}^{\mathrm{opp}},\, u_1u_2\mathcal{L}^{\mathrm{opp}},\, \mathcal{L}\big).
\]
On the other hand, substituting our coordinate expressions:
\begin{align*}
\xi(z_2, x, z_1, y) &= g\big(U^{f_1(z_2)}\mathcal{L}^{\mathrm{opp}},\, \mathcal{L}^{\mathrm{opp}},\, U^{f_1(z_1)}\mathcal{L}^{\mathrm{opp}},\, \mathcal{L}\big),
\end{align*}
which can be rewritten as:
\[
\xi(z_2, x, z_1, y) = g\,U^{f_1(z_2)}\big(\mathcal{L}^{\mathrm{opp}},\, U^{-f_1(z_2)}\mathcal{L}^{\mathrm{opp}},\, U^{f_1(z_1)-f_1(z_2)}\mathcal{L}^{\mathrm{opp}},\, \mathcal{L}\big).
\]
Let $h = g_1^{-1}(g\,U^{f_1(z_2)})$. Comparing the above expressions for $\xi(z_2, x, z_1, y)$ yields
\[
\big(\mathcal{L}^{\mathrm{opp}},\, u_1\mathcal{L}^{\mathrm{opp}},\, u_1u_2\mathcal{L}^{\mathrm{opp}},\, \mathcal{L}\big) = h\big(\mathcal{L}^{\mathrm{opp}},\, U^{-f_1(z_2)}\mathcal{L}^{\mathrm{opp}},\, U^{f_1(z_1)-f_1(z_2)}\mathcal{L}^{\mathrm{opp}},\, \mathcal{L}\big),
\]
so it follows that $h \in \mathrm{Stab}_{\mathrm{Aut}(\mathfrak{sp}(2n,\mathbb{R}))}(\mathcal{L}^{\mathrm{opp}}, \mathcal{L})$, and
\[
u_1 = \exp\big( h(u^{-f_1(z_2)}) \big) 
\quad \text{and} \quad 
u_1u_2 = \exp\big( h(u^{f_1(z_1)-f_1(z_2)}) \big) ,
\]
which is equivalent to
\[
u_1 = \exp\big( h(u^{-f_1(z_2)}) \big)
\quad \text{and} \quad 
u_2 = \exp\big( h(u^{f_1(z_1)}) \big).
\]

Recall that any element in $\mathrm{Stab}(\mathcal{L}^{\mathrm{opp}}, \mathcal{L})$ either preserves the positive cone $u^{\mathrm{Pos}(n)}$ or maps it to $u^{-\mathrm{Pos}(n)}$ (see Guichard--Wienhard~\cite[Corollary~5.3]{GW2}). Since $u_1, u_2 \in \mathsf{U}_{\alpha}^{>0}$, the elements $f_1(z_2)$ and $-f_1(z_1)$ must lie in $\mathrm{Pos}(n)$ or $-\mathrm{Pos}(n)$ simultaneously. Note that this discussion applies to all counter-clockwise configurations $(z_2,x,z_1,y)$. Since we restricted $f_1(z_0) = M_0$ to lie in $\mathrm{Pos}(n)$ or $-\mathrm{Pos}(n)$ (depending on the position of $z_0$) at the beginning of the proof, we can substitute $z_1$ or $z_2$ for $z_0$ in the argument above. We conclude that for any such $z_1, z_2$, we have $f_1(z_2) \in \mathrm{Pos}(n)$ and $-f_1(z_1) \in \mathrm{Pos}(n)$, so $f_1(z_2) - f_1(z_1) \in \mathrm{Pos}(n)$.

\textbf{Part (2):}
The proof for $f_2$ is analogous. We choose the chart $\mathsf{U}_{\alpha}^{\mathrm{opp}} \to \mathcal{F}_{\alpha}$, $u \mapsto u\mathcal{L}$, and using the $g$ chosen in Part (1), we define $f_2: E_2 \to \mathrm{Sym}(n)$ such that $\xi(z) = g U_{f_2(z)}\mathcal{L}$, with $f_2(y)=0$.
Using the automorphism $\mathrm{Ad}_{\mathcal{N}}$ (where $\mathcal{N} = (\begin{matrix} 0 & \mathrm{Id}_n \\ \mathrm{Id}_n & 0 \end{matrix})$) which swaps $\mathcal{L}$ and $\mathcal{L}^{\mathrm{opp}}$ and maps $\mathsf{U}_{\alpha}^{>0}$ to $\mathsf{U}_{\alpha}^{\mathrm{opp}, >0}$, the statement reduces to the case proved in Part (1). Thus, for $x, z_1, z_2$ clockwise, we obtain $f_2(z_1) - f_2(z_2) \in \mathrm{Pos}(n)$.
\end{proof}

\begin{remark}\label{remark:variationofparameterization}
We will sometimes apply the above lemma by taking $E_1$ or $E_2$ to be the singleton $\{x\}$ or $\{y\}$, and sometimes we will let $E_1$ or $E_2$ be an arc in $S^1$ ending at $x$ or $y$.

We may also replace the condition $f_1(z_2) - f_1(z_1) \in \mathrm{Pos}(n)$ by $f_1(z_2) - f_1(z_1) \in -\mathrm{Pos}(n)$ without affecting the validity of the results.
Moreover, we note the following equivalent variations of the parameterization:
\begin{align*}
g\big(\mathcal{L}^{\mathrm{opp}},\, U^{f(z)}\mathcal{L}^{\mathrm{opp}},\, \mathcal{L}\big)
&= g\,\mathrm{Ad}_{\mathcal{N}}\big(\mathcal{L},\, U_{f(z)}\mathcal{L},\, \mathcal{L}^{\mathrm{opp}}\big) \\
&= g\,\mathrm{Ad}_{\Omega}\big(\mathcal{L},\, U_{-f(z)}\mathcal{L},\, \mathcal{L}^{\mathrm{opp}}\big) \\
&= g\big(\mathcal{L}^{\mathrm{opp}},\, U_{f(z)^{-1}}\mathcal{L},\, \mathcal{L}\big),
\end{align*}
where the last equality holds under the assumption that $f(z)$ is invertible.
\end{remark}

\subsection{Maximal Representations}\label{subsecction:maximalrepresentations}
\begin{definition}\label{definitionofmaximal}
Assume $\Gamma\subset \mathrm{PSL}(2,\mathbb{R})$ is a non-elementary discrete subgroup. A representation $\rho : \Gamma \to \mathrm{Sp}(2n, \mathbb{R})$ is called \emph{maximal} if there exists a continuous, positive, and $\rho$-equivariant map $\xi^{\alpha} : \Lambda(\Gamma) \to \mathcal{F}_{\alpha}$.
\end{definition}

Inspired by Goldman~\cite{Goldman1988}, the terminology “maximal representations” originated with Burger--Iozzi--Wienhard~\cite{maximalearlier, burger2008surfacegrouprepresentationsmaximal} to describe representations in $\mathrm{Hom}(\pi_1(S),G)$ that maximize the absolute value of the Toledo invariant, where $S$ is a closed surface and $G$ is a Hermitian Lie group of tube type. It turns out that these representations form connected components of $\mathrm{Hom}(\pi_1(S),G)/G$ (Burger--Iozzi--Wienhard~\cite[Section~1, Corollary~14]{burger2008surfacegrouprepresentationsmaximal}), which motivated higher Teichmüller theory. It can also be shown that for closed surface groups, the definition via the Toledo invariant is equivalent to Definition~\ref{definitionofmaximal} (for example see \cite[Section~1, Theorem~8]{burger2008surfacegrouprepresentationsmaximal} ), and the latter definition applies to general non-elementary discrete subgroups of $\mathrm{PSL}(2,\mathbb{R})$.

On the other hand, in the development of Hitchin representations (another prominent family in higher Teichmüller theory), Labourie~\cite{labourie2005anosovflowssurfacegroups} showed that each Hitchin representation is associated with an Anosov system and introduced the concept of Anosov representations. Moreover, Burger--Iozzi--Labourie--Wienhard~\cite[Theorem~6.1]{maximal} showed that when $G=\mathrm{Sp}(2n,\mathbb{R})$, maximal representations are $\alpha$-Anosov, and moreover, the $\rho$-equivariant positive limit map $\xi^{\alpha}$ coincides with the limit map of the associated Anosov representation.

There are numerous characterizations of Anosov representations (and their generalizations—relatively Anosov representations and transverse representations). Here we adopt the one that originated from Kapovich–Leeb–Porti~\cite{kapovich2017anosov}, which uses the dynamical properties of the limit map.

\begin{definition}[{Canary--Zhang--Zimmer~\cite[Theorem~1.2]{CaZhZi2} and ~\cite[Section~6.1]{CaZhZi3}}]\label{definitionofanosov}

Let $\Gamma\subset \mathrm{PSL}(2,\mathbb{R})$ be a non-elementary discrete subgroup and let $\rho':\Gamma\to \mathrm{SL}(d,\mathbb{R})$ be a representation. Let $I$ be a symmetric subset of $\{1,\dots,d-1\}$. We call $\rho'$ \emph{$I$-transverse} if there exists a $\rho'$-equivariant, continuous map $\xi^{I}:\Lambda(\Gamma)\to \mathcal{F}_I(\mathbb{R}^d)$ such that:
\begin{enumerate}
    \item $\xi^{I}$ is transverse; that is, for any $x\neq y\in \Lambda(\Gamma)$, $\xi^I(x)$ is transverse to $\xi^I(y)$.
    \item $\xi^I$ is strongly dynamics-preserving; that is, for any fixed $b_0\in \mathbb{D}$, if $(\gamma_n)\subset \Gamma$ satisfies $\gamma_n(b_0)\to x\in \Lambda(\Gamma)$ and $\gamma_n^{-1}(b_0)\to y\in \Lambda(\Gamma)$, then $\rho'(\gamma_n)F\to \xi^{I}(x)$ for any $F$ transverse to $\xi^{I}(y)$.
\end{enumerate}
And $\xi^{I}$ is called the $I$-limit map of $\rho'$ which can be directly checked to be unique.

Furthermore, if $\Gamma$ is geometrically finite, then $I$-transverse representations are called relatively $I$-Anosov representations. And if additionally $\Gamma$ has a compact fundamental domain, then  $I$-transverse representations are called $I$-Anosov representations.
\end{definition}

Under the geometrically finite assumption (so that $\rho'$ is relatively $I$-Anosov), we have
\begin{theorem}[{Canary--Zhang--Zimmer~\cite[Theorem~1.1]{CaZhZi2}}]
For each fixed $b_0\in \mathbb{D}$ and $k\in I$, there exist constants $a\ge1$ and $A\ge0$ such that for any $\gamma\in \Gamma$,
\[
\frac{1}{a}\,d_{\mathbb{D}}(b_0,\gamma(b_0)) - A
\;\le\;
\log \frac{\sigma_k(\rho'(\gamma))}{\sigma_{k+1}(\rho'(\gamma))}
\;\le\;
a\,d_{\mathbb{D}}(b_0,\gamma(b_0)) + A,
\]
where $\sigma_k(g)$ denotes the $k$-th singular value of $g\in \mathrm{SL}(d,\mathbb{R})$. 
\end{theorem}

And we have that maximal representations are transverse:
\begin{proposition}[{Yao~\cite[Proposition~1.1]{yzf} 
}]\label{anosov}
If $\Gamma\subset \mathrm{PSL}(2,\mathbb{R})$ is a non-elementary subgroup and $\rho:\Gamma\to \mathrm{Sp}(2n,\mathbb{R})$ is a maximal representation, then the positive map $\xi^{\alpha}$ is strongly dynamics-preserving. Hence $\rho$ is $n$-transverse when we regard $\mathrm{Sp}(2n,\mathbb{R})$ as a subgroup of $\mathrm{SL}(2n,\mathbb{R})$ and $\mathcal{F}_{\alpha}$ as a subset of $\mathcal{F}_{n}(\mathbb{R}^{2n})$. 
\end{proposition}
Recall $\Lambda^n$ denotes the $n$-th wedge representation. Under the same assumptions as the above proposition, it is direct to check:
\begin{corollary}\label{anosov2}
     Define $\rho':\Gamma\to \mathrm{SL}(\bigwedge^n \mathbb{R}^{2n})$ by $\rho'=\Lambda^n\circ \rho$, then $\rho'$ is $I_1$-transverse, where $I_1=\{1,\,\dim(\bigwedge^n \mathbb{R}^{2n})-1\}$. The limit map is given by
    \[
    \xi^{I_1}:\ \Lambda(\Gamma)\to \mathcal{F}_{I_1},\qquad
    x\longmapsto \bigl(\mathcal{P}_{\alpha}(\xi^{\alpha}(x)),\,\mathcal{P}_{\alpha}^{\mathrm{opp}}(\xi^{\alpha}(x))\bigr).
    \]
\end{corollary}

\section{Measurable $(1,1,2)$-Hypertransversality}\label{sec3}

In this section, we assume that $n>1$. Davalo~\cite[Theorem~5.2]{Colin} proved that if a maximal representation from a closed surface group is further \(\{n-1,n+1\}\)-Anosov (when view it as a representation to $\mathrm{SL}(2n,\mathbb{R})$), then its limit map satisfies a hypertransverse property, called property \(H_n\) introduced by Beyrer--Pozzetti~\cite{beyrerpozzettisopq} (which implies that the limit curve is $C^{1,\alpha}$ for some $\alpha>0$) following Labourie~\cite{labourie2005anosovflowssurfacegroups}. Our goal here is to generalize the hypertransversality to a measurable setting, without the $\{n-1,n+1\}$ Anosov hypothesis. 

As before, let $V = \bigwedge^{n} \mathbb{R}^{2n}$, and let $d$ denote its dimension. We assume that $\Gamma\subset \mathrm{PSL}(2,\mathbb{R})$ is a non-elementary discrete subgroup with $\Lambda(\Gamma)=\partial\mathbb{D}$, and that $\rho:\Gamma\to \mathrm{Sp}(2n,\mathbb{R})$ is a maximal representation with limit map $\xi^{\alpha}$. Define $\xi^{1}:\Lambda(\Gamma)\to \mathbb{P}(\bigwedge^n \mathbb{R}^{2n})$ by $\xi^{1}=\mathcal{P}_{\alpha}\circ \xi^{\alpha}$. For any choice of Riemannian metric on the ambient space, the images $\xi^{\alpha}(\Lambda(\Gamma))\subset \mathcal{F}_{\alpha}$ and $\xi^{1}(\Lambda(\Gamma))\subset \mathbb{P}(\bigwedge^{n}\mathbb{R}^{2n})$ are rectifiable curves (see Burger--Iozzi--Labourie--Wienhard~\cite[Corollary~6.3]{maximal}). Hence there is a well-defined class of Lebesgue measures on them—namely, the finite, nonvanishing $1$-dimensional Hausdorff measure (see Stein–Shakarchi~\cite[Chapter~3, Section~4; Chapter~7, Theorem~2.4]{steinreal}). And for almost every point $x$ of $\xi^{1}(\Lambda(\Gamma))$, a tangent line exists. Since $\xi^{1}$ is a homeomorphism between $\Lambda(\Gamma)$ and its image, any Lebesgue measure $m$ on $\xi^{1}(\Lambda(\Gamma))$ induces a pullback measure $m^{*}$ on $\Lambda(\Gamma)$, which is a regular Borel measure.

The $m^*$-almost everywhere defined tangent line determines a measurable map:
\[
\xi^{1,2}: \Lambda(\Gamma) \to \mathcal{F}_{1,2}(\bigwedge^{n} \mathbb{R}^{2n}),
\]
\[
x \mapsto (\xi^1(x), \ell(x)),
\]
where $\ell(x)$ is the projective tangent line to $\xi^1(\Lambda(\Gamma))$ at the point $\xi^1(x)$ (when it exists). Each tangent line can be identified with a $2$-dimensional subspace of $\bigwedge^{n} \mathbb{R}^{2n}$.

\begin{proposition}\label{tangent}
If  $\ell(x)$ is defined, then there exist $g \in \mathrm{Aut}(\mathfrak{sp}(2n,\mathbb{R}))$ and a nonzero matrix $M \in \overline{\mathrm{Pos}(n)}$ such that
\[
\ell(x) = g \cdot \mathrm{Span}\left(e_1 \wedge \cdots \wedge e_n,\ \sum_{i=1}^{n} e_1 \wedge \cdots \wedge e_{i-1} \wedge u_M e_i \wedge e_{i+1} \wedge \cdots \wedge e_n \right),
\]
where $u_M = \begin{pmatrix} 0 & 0 \\ M & 0 \end{pmatrix}$.
\end{proposition}

\begin{proof}
Pick $y \neq x \in \Lambda(\Gamma)$ and let $[x,y]$ denote the clockwise arc segment on $\Lambda(\Gamma)$ connecting $x$ to $y$. Choose a parameterization $p:[0,1] \to [x,y]$ with $p(0) = x$ and $p(1) = y$.

From Lemma~\ref{lemma:parameterizationofmaps_combined} and Remark~\ref{remark:variationofparameterization}, there exist $g \in \mathrm{Aut}(\mathfrak{sp}(2n,\mathbb{R}))$ and a continuous map $f:[0,1) \to \{0\}\cup\mathrm{Pos}(n)$ such that:
\begin{itemize}
    \item $\xi^1(x) = g \cdot \mathrm{Span}(e_1 \wedge \cdots \wedge e_n)$,
    \item $\xi^1(y) = g \cdot \mathrm{Span}(e_{n+1} \wedge \cdots \wedge e_{2n})$,
    \item $f(0) = 0$,
    \item $f(t_2) - f(t_1) \in \mathrm{Pos}(n)$ for $t_2 > t_1$,
    \item for all $t \in [0,1)$,
    \[
    \xi^1(p(t)) = g \cdot \mathrm{Span}\left(U_{f(t)} e_1 \wedge \cdots \wedge U_{f(t)} e_n\right),
    \]
    where $U_{f(t)} = \begin{pmatrix} \mathrm{Id}_n & 0 \\ f(t) & \mathrm{Id}_n \end{pmatrix}$.
\end{itemize}

If the tangent line $\ell(x)$ exists, we may choose the parametrization $p$ so that $f$ is differentiable at $t=0$ with nonzero derivative. The positivity condition on $f$ implies that $f'(0) = M$ for some $M \in \overline{\mathrm{Pos}(n)}-\{0\}$.

Differentiating the expression for $\xi^1(p(t))$ at $t=0$ gives the claimed formula:
\[
\ell(x) = g \cdot \mathrm{Span}\left(e_1 \wedge \cdots \wedge e_n,\ \sum_{i=1}^{n} e_1 \wedge \cdots \wedge e_{i-1} \wedge u_M e_i \wedge e_{i+1} \wedge \cdots \wedge e_n \right).
\]
\end{proof}

Recall from Subsection~\ref{identificationsofflags} that we defined a bilinear form $\omega'$ on $\bigwedge^{n} \mathbb{R}^{2n}$. Let $d=\dim \bigwedge^{n}\mathbb{R}^{2n}$, and recall that for any subspace $V\subset \bigwedge^{n}\mathbb{R}^{2n}$ we define its bi-annihilator (as $\omega'$ is symmetric or skew-symmetric) by
\[
\mathrm{Ann}(V):=\left\{\,w\in \bigwedge^{n}\mathbb{R}^{2n}\ \middle|\ \forall v\in V,\ \omega'(w,v)=0\,\right\}.
\]
We define the dual map:
\[
\xi^{d-2,d-1}:\ \Lambda(\Gamma)\to \mathcal{F}_{d-2,d-1}\!\left(\bigwedge^{n}\mathbb{R}^{2n}\right)
\]
by setting
\[
x\longmapsto \bigl(\mathrm{Ann}(\ell(x)),\,\mathrm{Ann}(\xi^{1}(x))\bigr).
\]
\begin{remark}\label{ngreaterthan2}
    Using the notation and formula in Proposition~\ref{tangent}, one can verify the following:
    \begin{enumerate}
        \item $\omega'|_{\xi^1(x)\times \xi^1(x)} = 0$.

        \item If $n>2$ and $\xi^{1}(\Lambda(\Gamma))$ is differentiable at $\xi^1(x)$, then $\omega'|_{\ell(x)\times \ell(x)}=0$.

        \item If $n=2$, then in the above notation, 
    \[
        \omega'\!\big( (u_M e_1) \wedge e_2 + e_1\wedge u_M e_2 ,(u_M e_1) \wedge e_2 + e_1\wedge u_M e_2\big)
        = -2\det M.
    \]
    \end{enumerate}
    So for $n>2$, we can indeed define a symmetric map
    \[
    \xi^{1,2,d-2,d-1}:\ \Lambda(\Gamma)\to \mathcal{F}_{1,2,d-2,d-1}\!\left(\bigwedge^{n}\mathbb{R}^{2n}\right)
    \]
    by 
    \[
    x\longmapsto \bigl(\xi^1(x), \ell(x), \mathrm{Ann}(\ell(x)),\,\mathrm{Ann}(\xi^{1}(x))\bigr).
    \]
\end{remark}

The limit map has a measurable $(1,1,2)$-hypertransversality in the following sense:

\begin{proposition}\label{measurablehypertrans}
For three distinct points $x, y, z \in \Lambda(\Gamma)$, if $\ell(z)$ exists, then
\[
\xi^1(x) \oplus \xi^1(y) \oplus \mathrm{Ann}(\ell(z))
\]
is a direct sum.
\end{proposition}

As noted above, when $\Gamma$ is a closed surface group and $\rho$ is $\{n-1, n+1\}$-Anosov, one can directly verify that the wedged representation $\rho' = \Lambda^n \circ \rho$ is $\{1, 2, d-2, d-1\}$-Anosov and $\xi^1(\Lambda(\Gamma))$ is a $C^1$ curve (see Davalo~\cite{Colin}). In this setting, it can be checked that $\ell(x)\subset \mathrm{Ann}(\ell(x))$ always holds, and $\xi^{1,2,d-2,d-1}$ introduced above serves as the limit map. In this setting, the hypertransversality has been proved \cite{Colin}.

\begin{proof}
From Lemma~\ref{lemma:parameterizationofmaps_combined}, Remark~\ref{remark:variationofparameterization}, and Proposition~\ref{tangent}, we can find $g \in \mathrm{Aut}(\mathfrak{sp}(2n,\mathbb{R}))$, $M \in \mathrm{Pos}(n)$, and $0\not =N  \in \overline{\mathrm{Pos}(n)}$ such that:
\begin{align*}
\xi^1(x) &= g \, \mathrm{Span}(e_{n+1} \wedge \cdots \wedge e_{2n}), \\
\xi^1(y) &= g \, U_M \, \mathrm{Span}(e_1 \wedge \cdots \wedge e_n), \\
\xi^1(z) &= g \, \mathrm{Span}(e_1 \wedge \cdots \wedge e_n), \\
\ell(z) &= g \, \mathrm{Span}\left( e_1 \wedge \cdots \wedge e_n, \sum_{i=1}^{n} e_1 \wedge \cdots \wedge e_{i-1} \wedge u_N e_i \wedge e_{i+1} \cdots \wedge e_n \right),
\end{align*}
where $U_M = \begin{pmatrix} \mathrm{Id}_n & 0 \\ M & \mathrm{Id}_n \end{pmatrix}$ and $u_N = \begin{pmatrix} 0 & 0 \\ N & 0 \end{pmatrix}$.

Denote
\[
v_x := e_{n+1} \wedge \cdots \wedge e_{2n}, \quad v_y := U_M(e_1 \wedge \cdots \wedge e_n).
\]
Suppose the sum $\xi^1(x) + \xi^1(y) + \mathrm{Ann}(\ell(z))$ is not direct. Then there exists a nonzero pair $(c_1, c_2) \in \mathbb{R}^2$ such that for all $w \in g^{-1}\ell(z)$,
\[
\omega'(w, c_1 v_x + c_2 v_y) = 0.
\]
Take $w_1 := e_1 \wedge \cdots \wedge e_n$. Then:
\[
\omega'(w_1, v_x) = \det \begin{pmatrix} \mathrm{Id}_n & 0 \\ 0 & \mathrm{Id}_n \end{pmatrix} = 1, \quad
\omega'(w_1, v_y) = \det \begin{pmatrix} \mathrm{Id}_n & \mathrm{Id}_n \\ 0 & M \end{pmatrix} = \det M.
\]
So we have $c_1 + \det M \cdot c_2 = 0$. Now denote
\[
w_2 := \sum_{i=1}^{n} e_1 \wedge \cdots \wedge e_{i-1} \wedge u_N e_i \wedge e_{i+1} \cdots \wedge e_n,
\]
and notice that
\[
w_2 = \lim_{t \to 0} \frac{U_{tN}(e_1 \wedge \cdots \wedge e_n) - e_1 \wedge \cdots \wedge e_n}{t}.
\]
Hence,
\[
\omega'(w_2, v_x) = \lim_{t \to 0} \frac{\det \begin{pmatrix} \mathrm{Id}_n & 0 \\ tN & \mathrm{Id}_n \end{pmatrix} - \det \begin{pmatrix} \mathrm{Id}_n & 0 \\ 0 & \mathrm{Id}_n \end{pmatrix}}{t} = 0
\]
and
\begin{align*}
\omega'(w_2, v_y) &= \lim_{t \to 0} \frac{\det \begin{pmatrix} \mathrm{Id}_n & \mathrm{Id}_n \\ tN & M \end{pmatrix} - \det \begin{pmatrix} \mathrm{Id}_n & \mathrm{Id}_n \\ 0 & M \end{pmatrix}}{t} \\
&= \left. \frac{d}{dt} \right|_{t=0} \det \begin{pmatrix} \mathrm{Id}_n & \mathrm{Id}_n \\ tN & M \end{pmatrix} \\
&=\left. \frac{d}{dt} \right|_{t=0} \det \begin{pmatrix} \mathrm{Id}_n & \mathrm{Id}_n \\ 0 & M-tN \end{pmatrix}
\\
&=\left. \frac{d}{dt} \right|_{t=0} \det (M(\mathrm{Id}_n-tM^{-1}N))
\\
&= -\det M \cdot \mathrm{tr}(M^{-1}N).
\end{align*}

We now prove the following lemma:
\begin{lemma}\label{lemmaprodtrace}
If $A \in \mathrm{Pos}(n)$ and $B \in \overline{\mathrm{Pos}}(n) \setminus \{0\}$, then $\mathrm{tr}(AB) > 0$.
\end{lemma}
\begin{proof}
As $A\in \mathrm{Pos}(n)$, we can pick $A'\in \mathrm{Pos}(n)$ such that $A = (A')^2$, so $\mathrm{tr}(AB) = \mathrm{tr}(A'A'B) = \mathrm{tr}(A'BA')$. As $B \in \overline{\mathrm{Pos}}(n) \setminus \{0\}$, $\mathrm{tr}(B)>0$, so the congruence-transformed trace $\mathrm{tr}(A'BA')>0$. 
\end{proof}

Returning to the proof, since $M$ and $M^{-1}$ are positive definite, and $N \in \overline{\mathrm{Pos}}(n) \setminus \{0\}$, we conclude from Lemma~\ref{lemmaprodtrace} that $\mathrm{tr}(M^{-1}N) > 0$. Therefore:
\[
\omega'(w_2, c_1 v_x + c_2 v_y) = c_1 \cdot 0 - c_2 \cdot \det M \cdot \mathrm{tr}(M^{-1}N) = -c_2 \cdot \det M \cdot \mathrm{tr}(M^{-1}N).
\]
Hence, the vanishing of $\omega'$ implies $c_2 = 0$, and from $c_1 + \det M \cdot c_2 = 0$, it follows that $c_1 = 0$ as well, contradicting the assumption that $(c_1, c_2) \neq (0, 0)$. Thus,
\[
\xi^1(x) \oplus \xi^1(y) \oplus \mathrm{Ann}(\ell(z))
\]
is indeed a direct sum.
\end{proof}

Let $\mathcal{ET}$ be the set of distinct pairs $(x,y) \in \Lambda(\Gamma) \times \Lambda(\Gamma)$ for which both $\ell(x)$ and $\ell(y)$ are defined and $\ell(y)$ is transverse to $\mathrm{Ann}(\ell(x))$. Notice that $(x,y) \in \mathcal{ET}$ is equivalent to the nondegeneracy of the pairing $\omega'(\ell(x), \ell(y))$. Consequently, this is equivalent to $(y,x) \in \mathcal{ET}$.

\begin{corollary}\label{transverse}
We have $m^* \times m^*(\mathcal{ET}) > 0$.
\end{corollary}

\begin{proof}
Fix a point $x \in \Lambda(\Gamma)$ such that $\ell(x)$ exists. Consider the projection
\[
\pi_x' : \bigwedge\nolimits^n \mathbb{R}^{2n} \to \big(\bigwedge\nolimits^n \mathbb{R}^{2n}\big) / \mathrm{Ann}(\ell(x)),
\]
and let $\pi_x$ denote the induced map on projective space:
\[
\pi_x : \mathbb{P}\left( \bigwedge\nolimits^n \mathbb{R}^{2n} -\mathrm{Ann}(\ell(x))\right) \to \mathbb{P} \left( \bigwedge\nolimits^n \mathbb{R}^{2n} / \mathrm{Ann}(\ell(x)) \right).
\]

By Proposition~\ref{measurablehypertrans}, the composition $\pi_x \circ \xi^1$ is injective on $\Lambda(\Gamma) \setminus \{x\}$. $\pi_x$ is smooth, so clearly the image curve $\pi_x \circ \xi^1(\Lambda(\Gamma) \setminus \{x\})$ is also rectifiable. Fix a Riemannian metric on the quotient projective space. 

Let $E_x \subset \Lambda(\Gamma)$ be the set of points $y$ for which $\ell(y)$ exists but is \textbf{not} transverse to $\mathrm{Ann}(\ell(x))$. If we parameterize $\xi^1(\Lambda(\Gamma))$ in unit speed, then $\pi_x \circ \xi^1(\Lambda(\Gamma) \setminus \{x\})$ is parameterized in a local bounded speed, and for any $p\in E_x$, the speed at $p$ is zero.

If $E_x$ had full $m^*$-measure, then the speed of this parameterization is $0$ almost everywhere, so the curve $\pi_x \circ \xi^1$ would be constant almost everywhere, contradicting the injectivity of $\pi_x \circ \xi^1$ on $\Lambda(\Gamma) \setminus \{x\}$. Thus, the complement of $E_x$—the set of $y$ such that $\ell(y)$ exists and is transverse to $\mathrm{Ann}(\ell(x))$—has positive $m^*$-measure.

Since the set of $x$ for which $\ell(x)$ exists has full $m^*$-measure, applying Fubini's theorem we finish the proof.
\end{proof}

\section{Patterson--Sullivan Theory}\label{sec4}
As a preliminary step to proving the main theorem, we review the construction and dynamical properties of Patterson--Sullivan measures.

We will use $\tilde{\rho}$ to denote a general representation to $\mathrm{SL}(d,\mathbb{R})$, and use $\rho$ to denote a maximal representation to $\mathrm{Sp}(2n,\mathbb{R})$.

\subsection{Partial Iwasawa Cocycle}

Choose an orthonormal basis $\{e_1, e_2, \dots, e_d\}$ of $\mathbb{R}^d$. Let $\mathrm{N}(d)$ be the group of upper triangular matrices with all diagonal entries equal to $1$, let $\mathrm{A}(d)$ be the group of positive diagonal matrices in $\mathrm{SL}(d, \mathbb{R})$, and let $\mathrm{K}(d) = \mathrm{SO}(d)$. The \emph{Iwasawa decomposition}
\[
\mathrm{SL}(d, \mathbb{R}) = \mathrm{K}(d)\, \mathrm{A}(d)\, \mathrm{N}(d)
\]
asserts that for each $g \in \mathrm{SL}(d, \mathbb{R})$, there is a unique triple $(k(g), a(g), u(g)) \in \mathrm{K}(d) \times \mathrm{A}(d) \times \mathrm{N}(d)$ such that
\[
g = k(g)\, a(g)\, u(g).
\]

Choose the Cartan subalgebra $\mathfrak{a}_d \subset \mathfrak{sl}(d, \mathbb{R})$ by
\[
\mathfrak{a}_d := \left\{ \mathrm{diag}(\lambda_1, \lambda_2, \dots, \lambda_d) \mid \sum_{i=1}^{d} \lambda_i = 0 \right\}.
\]

For each $1 \leq i < d$, define the $i$-th \emph{simple root} $\alpha_i \in \mathfrak{a}_d^*$ by
\[
\alpha_i\left(\mathrm{diag}(\lambda_1, \dots, \lambda_d)\right) = \lambda_i - \lambda_{i+1},
\]
and the $i$-th \emph{fundamental weight} $\omega_i \in \mathfrak{a}_d^*$ by
\[
\omega_i\left(\mathrm{diag}(\lambda_1, \dots, \lambda_d)\right) = \lambda_1 + \dots + \lambda_i.
\]

The \emph{positive Weyl chamber} $\mathfrak{a}_d^+ \subset \mathfrak{a}_d$ is defined as
\[
\mathfrak{a}_d^+ := \left\{ \mathrm{diag}(\lambda_1, \dots, \lambda_d) \mid \sum_{i=1}^{d} \lambda_i = 0,\, \lambda_1 > \lambda_2 > \dots > \lambda_d \right\},
\]
and $(\mathfrak{a}_d^+)^* \subset \mathfrak{a}_d^*$ is the cone of linear functionals that are positive on $\mathfrak{a}_d^+$, which is exactly the set of non-negative coefficient combinations of $\alpha_i$.

For a multi-index $I = \{i_1, i_2, \dots, i_k\}$ with $1 \leq i_1 < i_2 < \dots < i_k < d$, define the subspace
\[
\mathfrak{a}_I := \bigcap_{j \notin I} \ker \alpha_j \subset \mathfrak{a}_d,
\]
and define the projection $\pi_I : \mathfrak{a}_d \to \mathfrak{a}_I$ by the condition
\[
\forall i_s \in I,\; \forall a \in \mathfrak{a}_d,\quad \omega_{i_s}(\pi_I(a)) = \omega_{i_s}(a).
\]

Recall that $\mathcal{F}_I$ denotes the $I$-flag manifold. For any $F \in \mathcal{F}_I$ and $i_s \in I$, denote by $F_{i_s}$ the corresponding $i_s$-dimensional subspace.

\begin{definition}[Partial Iwasawa Cocycle]\label{iwasawa}
For each $I$, the \emph{partial Iwasawa cocycle} is the map
\[
B_I : \mathrm{SL}(d, \mathbb{R}) \times \mathcal{F}_I \to \mathfrak{a}_I
\]
defined by the condition
\[
\forall\, g\in \mathrm{SL}(d,\mathbb{R}),\ F\in \mathcal{F}_I,\ i_s\in I,\quad
\omega_{i_s}\!\left(B_I(g,F)\right)
= \log \frac{\left\|\,g\, v_{i_s}\right\|}{\left\| v_{i_s}\right\|},
\]
where $v_{i_s}$ is a non-zero vector in $\wedge^{i_s}F_{i_s}$ (unique up to scalar), and $\|\cdot\|$ denotes the norm on $\wedge^{i_s}\mathbb{R}^d$ induced by the Euclidean norm on $\mathbb{R}^d$.
\end{definition}

There is an equivalent, more geometric definition: Let $F' \in \mathcal{F}_I$ be the standard flag
\[
\mathrm{Span}(e_1, \dots, e_{i_1}) \subset \mathrm{Span}(e_1, \dots, e_{i_2}) \subset \dots \subset \mathrm{Span}(e_1, \dots, e_{i_k}).
\]
Then for any $F =k' F'$ with $k' \in \mathrm{K}(d)$, we have
\[
B_I(g, F) = \pi_I(\log a(gk')),
\]
where $a(gk')$ denotes the $\mathrm{A}(d)$ component in the Iwasawa decomposition of $gk'$. This expression is independent of the choice of $k'$ satisfying $F = k' F'$.

Finally, $B_I$ satisfies the cocycle identity:
\[
B_I(gh, F) = B_I(g, hF) + B_I(h, F).
\]
For further details, see for example Quint~\cite{quintoriginal} or Pozzetti--Sambarino--Wienhard~\cite[Subsection 4.6]{Liplim}.

\subsection{Definition of Patterson--Sullivan Measure}

We make the following assumptions:
$\Gamma$ is a non-elementary discrete subgroup of $\mathrm{PSL}(2,\mathbb{R})$,
$\tilde{\rho} : \Gamma \to \mathrm{SL}(d, \mathbb{R})$ is a representation,
$I \subset \{1, \dots, d-1\}$,
and there exists a Borel measurable, $\tilde{\rho}$-equivariant map
\[
\xi^I : \Lambda(\Gamma) \to \mathcal{F}_I.
\]

For each $\phi \in (\mathfrak{a}_d^+)^*$, define the \emph{critical exponent} $\delta_{\tilde{\rho}}(\phi)$ of the associated Poincaré series by
\[
\delta_{\tilde{\rho}}(\phi) := \inf \left\{ s \in [0, \infty) \;\middle|\; \sum_{\gamma \in \Gamma} e^{-s\, \phi(\kappa_d(\tilde{\rho}(\gamma)))} < \infty \right\},
\]
where $\kappa_d : \mathrm{SL}(d, \mathbb{R}) \to \overline{\mathfrak{a}_d^+}$ is the \emph{Cartan projection}, given by
\[
\kappa_d(g) = \mathrm{diag}\left(\log \sigma_1(g), \log \sigma_2(g), \dots, \log \sigma_d(g)\right).
\]
Here $\sigma_i(g)$ denotes the $i$-th singular value of $g$.

We identify the dual space $\mathfrak{a}_I^* \subset \mathfrak{a}_d^*$ with the span of the fundamental weights $\omega_{i_s}$ for $i_s \in I$:
\[
\mathfrak{a}_I^* := \mathrm{Span}_{i_s \in I}(\omega_{i_s}) \subset \mathfrak{a}_d^*.
\]

\begin{definition}[Patterson--Sullivan Measure]
Let $\phi \in \mathfrak{a}_I^*$ be such that $\delta_{\tilde{\rho}}(\phi) < \infty$. A \emph{$\phi$-Patterson--Sullivan measure} (for $\tilde{\rho}$) is a regular Borel probability measure $\mu$ supported on the image $\xi^I(\Lambda(\Gamma))$ satisfying the transformation rule
\[
\frac{d\, \gamma_* \mu}{d\mu}(-) = e^{-\delta_{\tilde{\rho}}(\phi) \cdot \phi\left(B_I(\tilde{\rho}(\gamma)^{-1}, -)\right)},
\]
where $\gamma_* \mu$ denotes the pushforward of $\mu$ under the action of $\gamma$, and the equality holds in the sense of Radon--Nikodym derivatives.
\end{definition}

Note that in many contexts, $\xi^{I}$ is assumed to be a \textbf{continuous} limit map of an Anosov representation. However, in this work, we only assume it to be measurable, owing to the construction of the Unstable Jacobian we introduce below.

\subsection{The Unstable Jacobian} 

We introduce two Patterson--Sullivan measures associated to a maximal representation \(\rho:\Gamma \to \mathrm{Sp}(2n,\mathbb{R})\) where $\Gamma$ is a lattice. Let \(V := \bigwedge^n(\mathbb{R}^{2n})\), \(\rho' := \Lambda^n \circ \rho\), and let \(d := \dim V = \binom{2n}{n}\). 

Let \(\xi^{\alpha}:\Lambda(\Gamma)\to \mathcal{F}_{\alpha}\) be the limit map associated with the maximal representation. It induces a limit map \(\xi^1:\Lambda(\Gamma)\to \mathcal{F}_1(V)\) via the \(n\)-fold wedge product. Fix a Riemannian metric on $\mathcal{F}_1(V)$. We saw in Section~\ref{sec3} that \(\xi^{1}(\Lambda(\Gamma))\) is a rectifiable curve. Let \(m\) be the 1-dimensional Hausdorff measure on \(\xi^{1}(\Lambda(\Gamma))\), and let \(m^*\) be the pullback of \(m\) to \(\Lambda(\Gamma)\).
We extend \(\xi^1\) to an equivariant and \(m^*\)-Borel measurable map
\[
\xi^{1,2}(x):=(\xi^1(x),\,\ell(x)).
\]

Let \(\mathcal{J}_2 := 2\omega_1 - \omega_2 \in \mathfrak{a}_{1,2}^*\). Using properties of the wedge product, one verifies that for all \(g \in \mathrm{Sp}(2n,\mathbb{R})\),
\[
\alpha(\kappa(g)) = \mathcal{J}_2(\kappa_d(\Lambda^n g)).
\]

We have:
\begin{proposition}[{Pozzetti--Sambarino--Wienhard~\cite[Theorem~1.4]{Liplim} and Yao~\cite[Theorem~1.2]{yzf}}]\label{propcritical1} For a maximal representation \(\rho\) from a lattice,  \(\delta_{\rho}(\alpha) = 1\). As a result, \(\delta_{\rho'}(\mathcal{J}_2) = 1\).
\end{proposition}

Pozzetti, Sambarino, and Wienhard constructed a \(\mathcal{J}_2\)-Patterson–Sullivan measure supported on \(\xi^{1,2}(\Lambda(\Gamma))\) in \cite[Section~6]{Liplim} as follows. For \(m^*\)-almost every \(x\in \Lambda(\Gamma)\), the tangent line to the curve \(\xi^{1}(\Lambda(\Gamma))\) at \(\xi^1(x)\) is identified with \(\mathrm{Hom}(\xi^1(x),\xi^1(x)^\perp)\), where \(\xi^1(x)^\perp\) is the one-dimensional subspace of \(\ell(x)\) perpendicular to \(\xi^1(x)\). Here “perpendicular’’ is with respect to the Euclidean inner product on \(\bigwedge^n \mathbb{R}^{2n}\) induced from that on \(\mathbb{R}^{2n}\). Define the \emph{Unstable Jacobian} as the volume form \(\Omega_{\xi^{1}(x)}\) on \(\mathrm{Hom}(\xi^1(x),\xi^1(x)^\perp)\) by
\[
\Omega_{\xi^{1}(x)}(\phi_1):=\frac{v\wedge \phi_1(v)}{\|v\|^2},
\]
for any \(\phi_1\in \mathrm{Hom}(\xi^1(x),\xi^1(x)^\perp)\) and any nonzero \(v\in \xi^1(x)\). Here \(\|\cdot\|\) denotes the Euclidean norm, and \(v\wedge \phi_1(v)\) is identified with the scalar obtained by evaluating against the unit volume form induced by the Euclidean structure on \(\ell(x)\). This definition is independent of the choice of \(v\).

The following proposition is about the existence of $\mathcal{J}_2$-Patterson--Sullivan measures associated to $\rho'$.
\begin{proposition}[{Pozzetti--Sambarino--Wienhard~\cite[Lemma~6.3 and proof of Proposition~6.4]{Liplim}}]\label{defofjac}
After scalar renormalization, the integration of the form \(\Omega_{\xi^{1}(x)}\) induces a measure on \(\xi^{1}(\Lambda(\Gamma))\) which is mutually absolutely continuous with the 1-dimensional Hausdorff measure, which pushes forward to a \(\mathcal{J}_2\)-Patterson–Sullivan measure on \(\xi^{1,2}(\Lambda(\Gamma))\) via
\[
\xi^{1,2}\circ (\xi^{1})^{-1}:\ \xi^{1}(\Lambda(\Gamma))\to \xi^{1,2}(\Lambda(\Gamma)),
\]
denoted by \(\mu_{\mathcal{J}_2}\).
\end{proposition}

\subsection{The \(\hat{\omega}_{\alpha}\)-Patterson--Sullivan Measure}

We still assume that $\rho: \Gamma\to \mathrm{Sp}(2n,\mathbb{R})$ is a maximal representation from a lattice and let $\rho' = \Lambda^n\circ \rho$. Let $d$ denote the dimension of $V =\bigwedge^n\mathbb{R}^{2n} $. Recall that \(\mathfrak{a} \subset \mathfrak{sp}(2n, \mathbb{R})\) denotes the Cartan subalgebra, which consists of diagonal elements. The fundamental weight \(\omega_{\alpha} \in \mathfrak{a}^*\) is given by
\[
\omega_{\alpha}(\mathrm{diag}(\lambda_1, \dots, \lambda_n, -\lambda_1, \dots, -\lambda_n)) = \sum_{i=1}^n \lambda_i.
\]
It is direct to check that for all \(g \in \mathrm{Sp}(2n, \mathbb{R})\),
\[
\omega_{\alpha}(\kappa(g)) = \omega_1(\kappa_d(\Lambda^n g)).
\]

Define the renormalized weight \(\hat{\omega}_{\alpha} := \frac{2}{n} \omega_{\alpha}\in \mathfrak{a}^{*}\), and let \(\hat{\omega}_1 := \frac{2}{n} \omega_1\in \mathfrak{a}_d^*\). Since \(\hat{\omega}_{\alpha} \ge \alpha\) on \(\overline{\mathfrak{a}^+}\), we have
\[
\delta_{\rho}(\hat{\omega}_{\alpha}) \le \delta_{\rho}(\alpha) = 1.
\]
So for \(\rho' = \Lambda^n \circ \rho\), 
\[
\delta_{\rho'}(\hat{\omega}_1) = \delta_{\rho}(\hat{\omega}_{\alpha}) \le 1.
\]

Corollary~\ref{anosov2} shows that $\rho'$ is relatively $\{1,d-1\}$-Anosov. From Canary--Zhang--Zimmer~\cite[Theorem~1.1]{CaZhZi-relatively} and \cite[Theorem~1.4]{CaZhZi3} (when $\Gamma$ is a closed surface group, see also Sambarino~\cite{sambarino2014orbitalcountingproblemhyperconvex} and Dey--Kapovich~\cite{Dey_2022}), we have

\begin{proposition}\label{defofps}
There exists a \(\hat{\omega}_1\)-Patterson--Sullivan measure \(\mu_{\hat{\omega}_1}\) on \(\xi^{1}(\Lambda(\Gamma))\).
\end{proposition}

\subsection{Gromov product}\label{gromovprod}

The Patterson--Sullivan measure is not invariant under the action of \(\Gamma\); it is only quasi-invariant. Although we can define ergodicity for quasi-invariant measures, classical ergodic theory about invariant measures does not directly apply. However, one can use the Patterson--Sullivan measures to construct invariant measures, and the key ingredient is the Gromov product to be introduced below.

Let \(V\) be a \(d\)-dimensional real vector space with a Euclidean inner product. Fix a symmetric index set \(I\subset \{1,...,d-1\}\), and let \(\mathcal{TF}_{I}(V) \subset \mathcal{F}_{I}(V) \times \mathcal{F}_{I}(V)\) denote the open subset consisting of transverse pairs. There exists a continuous map
\[
[-,-]_I : \mathcal{TF}_{I}(V) \to \mathfrak{a}_I,
\]
called the \emph{Gromov product}, defined by the condition
\[
\omega_{s}([F_1,F_2]_I)
= \log \|\frac{\bigwedge_{i=1}^{s} e_{1,i} }{\|\bigwedge_{i=1}^{s} e_{1,i}\|}\wedge \frac{\bigwedge_{j=1}^{d-s} e_{2,j}}
{ \|\bigwedge_{j=1}^{d-s} e_{2,j}\| }\|,
\qquad s \in I,
\]
where \(\{e_{1,i}\}_{i=1}^{s}\) is a basis of \((F_{1})_{s}\), \(\{e_{2,j}\}_{j=1}^{d-s}\) is a basis of \((F_{2})_{d-s}\), and \(\|\cdot\|\) denotes the norms on the exterior powers induced by the Euclidean product on \(V\). The definition of \([-, -]_I\) does not depend on the choices of \(e_{1,i}\) and \(e_{2,j}\). And clearly, \(\omega_{s}([F_1,F_2]_I) \le 0\) always holds.

It is straightforward to check that an equivalent definition of the Gromov product is as follows. For each $s\in I$, choose a basis $\{e_{1,i}\}_{i=1}^s$ of $(F_1)_s$ and a basis $\{e'_{2,j}\}_{j=1}^s$ of $((F_{2})_{d-s})^{\perp}$ (orthogonal with respect to the fixed Euclidean inner product $\langle-,-\rangle$ on $V$). Then
\[
\omega_{s}([F_1,F_2]_I)
= \log \frac{\bigl|\det\bigl(\langle e_{1,i},e'_{2,j}\rangle\bigr)_{1\le i,j\le s}\bigr|}
{\| \bigwedge_{i=1}^{s} e_{1,i} \|  \| \bigwedge_{j=1}^{s} e'_{2,j} \|}.
\]

Moreover, it is straightforward to check that the Gromov product satisfies the cocycle identity
\[
[g(F_1), g(F_2)]_I - [F_1, F_2]_I \;=\; -\, B_I(g, F_1) \;-i|_{\mathfrak{a}_I}\circ \; B_I(g, F_2),
\]
for all \(g \in \mathrm{SL}(V)\) and all \((F_1, F_2) \in \mathcal{TF}_I(V)\). Here, \(i\) denotes the involution on \(\mathfrak{a}_d\) given by
\[
i\left(\mathrm{diag}(\lambda_1, \lambda_2, \dots, \lambda_d)\right) = \mathrm{diag}(-\lambda_d, -\lambda_{d-1}, \dots, -\lambda_1).
\]
Since \(I\) is symmetric, the involution \(i\) restricts naturally to \(\mathfrak{a}_I\). 

We now prove a computational lemma for later use. We let $V = \bigwedge^n \mathbb{R}^{2n}$ whose dimension is denoted $d$. For any \(M \in \overline{\mathrm{Pos}}(n)\), recall that
\[
u_M = \begin{pmatrix} 0 & 0 \\ M & 0 \end{pmatrix}, \quad u^M = \begin{pmatrix} 0 & M \\ 0 & 0 \end{pmatrix}.
\]
Let \(N, M \in \overline{\mathrm{Pos}}(n)-\{0\}\), and define the following vectors:
\[
v_1 = e_1 \wedge \cdots \wedge e_n, \quad w_1 = \sum_{i=1}^{n} e_1 \wedge \cdots \wedge u_N e_i \wedge \cdots \wedge e_n,
\]
\[
v_2 = e_{n+1} \wedge \cdots \wedge e_{2n}, \quad w_2 = \sum_{i=1}^{n} e_{n+1} \wedge \cdots \wedge u^M e_{n+i} \wedge \cdots \wedge e_{2n}.
\]
Let $I_2 = \{1,2,d-1,d-2\}$. If $n>2$, then from Remark~\ref{ngreaterthan2} we can define $F_1,F_2\in \mathcal{F}_{I_2}(V)$ by
\[
F_i = \mathrm{Span}(v_i) \subset \mathrm{Span}(v_i, w_i) \subset \mathrm{Ann}(\mathrm{Span}(v_i, w_i)) \subset \mathrm{Ann}(\mathrm{Span}(v_i)),\quad i = 1,2.
\]

\begin{lemma}\label{Gromovproduct}
The subspaces \(F_1\) and \(F_2\) are transverse if and only if \(\mathrm{tr}(NM) > 0\). In this case, we have
\[
e^{\omega_2[F_1,F_2]_{I_2}} = \frac{\mathrm{tr}(NM)}{\|M\|_2 \cdot \|N\|_2},
\]
where \(\|\cdot\|_2\) denotes the Frobenius norm compatible with the Euclidean product on $\mathbb{R}^{2n}$. Explicitly, $\|M\|_2 = \sqrt{\mathrm{tr}(M^TM)}$.
\end{lemma}
\begin{proof}
Note that the Euclidean inner product on $\bigwedge^n \mathbb{R}^{2n}$ can be written as
\begin{equation}\label{equationofnorm}
\langle\,\cdot,\cdot\,\rangle=\omega'(\,\cdot\,,P(\,\cdot\,)),    
\end{equation}
where $P$ is the linear map (indeed it is the Hodge star operator) sending a basis vector $e_{i_1}\wedge\cdots\wedge e_{i_n}$ to $\pm\,e_{j_1}\wedge\cdots\wedge e_{j_n}$, with
\(
\{i_1,\dots,i_n,j_1,\dots,j_n\}=\{1,\dots,2n\}
\)
(the sign $\pm$ depends on the indices $i_1,\dots,i_n,j_1,...,j_n$ such that $(e_{i_1}\wedge\cdots\wedge e_{i_n})\wedge (\pm\,e_{j_1}\wedge\cdots\wedge e_{j_n}) = e_1\wedge \cdots\wedge e_{2n}$).

Moreover, one checks directly that for any $j>n$ and $i\le n$,
\[
P\!\bigl(e_{1}\wedge\cdots\wedge e_{i-1}\wedge e_j\wedge e_{i+1}\wedge\cdots\wedge e_{n}\bigr)
= -\,e_{n+1}\wedge\cdots\wedge e_{j-1}\wedge e_{i}\wedge e_{j+1}\wedge\cdots\wedge e_{2n},
\]
and
\[
P\!\bigl(e_{1}\wedge\cdots\wedge e_{n}\bigr)=e_{n+1}\wedge\cdots\wedge e_{2n}.
\]
Note that $\omega'(v,w)=0$ is equivalent to $\langle v,P^{-1}(w)\rangle=0$. Set
\[
v_2' := P^{-1}(v_2)=e_1\wedge\cdots\wedge e_n,\qquad
w_2' := P^{-1}(w_2)=-\sum_{i=1}^{n} e_1\wedge\cdots\wedge (u_M e_{i})\wedge\cdots\wedge e_n,
\]
where the equality for $w_2'$ comes from the fact that
\[
P^{-1}(e_{n+1}\wedge \cdots\wedge u^{E_{ji}}e_{n+i}\wedge\cdots e_{2n}) = -e_1\wedge\cdots\wedge (u_{E_{ij}} e_{j})\wedge\cdots\wedge e_n
\]
holds as $E_{ij}$ running through all $n\times n$ matrix units. Clearly, $v_2'$ and $w_2'$ span the orthogonal complement of $\mathrm{Ann}(\mathrm{Span}(v_2,w_2))$.

The 2-dimensional component of $F_1$ and $(d-2)$-dimensional component of $F_2$ are transverse if and only if
\[
\det \begin{pmatrix}
\langle v_1,v'_2\rangle & \langle w_1,v'_2\rangle \\[2pt]
\langle v_1,w'_2\rangle & \langle w_1,w'_2\rangle
\end{pmatrix}= \det \begin{pmatrix}
\omega'(v_1,v_2) & \omega'(w_1,v_2) \\[2pt]
\omega'(v_1,w_2) & \omega'(w_1,w_2)
\end{pmatrix}\neq 0.
\]

Using the fact that 
\[
w_1 = \frac{d}{dt}|_{t= 0} U_{tN} (e_{1}\wedge \cdots\wedge e_n),\quad w_2 = \frac{d}{dt}|_{t = 0}  U^{tM}(e_{n+1}\wedge \cdots \wedge e_{2n}),
\]
one can repeat the computation in Proposition~\ref{measurablehypertrans} to get $\omega'(v_1,v_2) = 1$, $\omega'(w_1,v_2)=0$, $\omega'(v_1,w_2) = 0$ and $\omega'(w_1,w_2) = -\mathrm{tr}(NM)$. As a result we have
\[
|\det \begin{pmatrix}
\omega'(v_1,v_2) & \omega'(w_1,v_2) \\[2pt]
\omega'(v_1,w_2) & \omega'(w_1,w_2)
\end{pmatrix}|
= \mathrm{tr}(NM).
\]
So the $2$-dimensional component of $F_1$ is transverse to the $(d-2)$-dimensional component of $F_2$ if and only if $\mathrm{tr}(NM)>0$. Similarly, the $2$-dimensional component of $F_2$ is transverse to the $(d-2)$-dimensional component of $F_1$ if and only if $\mathrm{tr}(MN)>0$, which is equivalent to $\mathrm{tr}(NM)>0$.

Under this assumption, by the equivalent determinant definition of the Gromov product, we obtain
\begin{align*}
e^{\omega_2\,[F_1,F_2]_{I_2}}
&= 
\frac{\left|
\det \begin{pmatrix}
\langle v_1,v'_2\rangle & \langle w_1,v'_2\rangle \\[2pt]
\langle v_1,w'_2\rangle & \langle w_1,w'_2\rangle
\end{pmatrix}
\right|}{\|M\|_2\,\|N\|_2} \\[6pt]
&=
\frac{\left|
\det \begin{pmatrix}
\omega'(v_1,v_2) & \omega'(w_1,v_2) \\[2pt]
\omega'(v_1,w_2) & \omega'(w_1,w_2)
\end{pmatrix}
\right|}{\|M\|_2\,\|N\|_2} \\[6pt]
&=
\frac{\mathrm{tr}(NM)}{\|M\|_2\,\|N\|_2}.
\end{align*}

This proves the claim.
\end{proof}

\subsection{Bowen--Margulis--Sullivan Measure}\label{bmsdef}

We now construct two invariant measures associated with a maximal representation \(\rho \colon \Gamma \to \mathrm{Sp}(2n,\mathbb{R})\), where \(\Gamma\) is a lattice and \textbf{\(n\ge 3\)}. The case \(n=2\) is excluded in view of Remark~\ref{ngreaterthan2} (we want to construct symmetric maps), however, this will not cause trouble as we can reduce maximal representations to $\mathrm{Sp}(4,\mathbb{R})$ to general cases in later discussions. 

We retain the notation $V = \bigwedge^{n} \mathbb{R}^{2n}$, set $d = \dim(V)$, and let $\rho' = \Lambda^n \circ \rho$. Let $\mathfrak{a}_{d}$ be the Cartan subspace of $\mathfrak{sl}(d, \mathbb{R})$. Recalling that $\mathcal{J}_2 := 2\omega_1-\omega_2$ and $\hat{\omega}_1 := \frac{2}{n}\omega_1$ in $\mathfrak{a}_{d}^*$, we observe that for any $\gamma \in \Gamma$, $\mathcal{J}_2(\kappa_d(\rho'(\gamma))) = \alpha(\kappa(\rho(\gamma)))$ and $\hat{\omega}_{1}(\kappa_d(\rho'(\gamma))) = \hat{\omega}_{\alpha}(\kappa(\rho(\gamma)))$.

Recall from Proposition~\ref{propcritical1} that $\delta_{\rho'}(\mathcal{J}_2)=1$. Also recall that \(\mu_{\mathcal{J}_2}\) is the \(\mathcal{J}_2\)-Patterson--Sullivan measure on \(\xi^{1,2}(\Lambda(\Gamma))\) introduced in Proposition~\ref{defofjac}, and \(\mu_{\hat{\omega}_{1}}\) is the \(\hat{\omega}_{1}\)-Patterson--Sullivan measure on \(\xi^{1}(\Lambda(\Gamma))\) introduced in Proposition~\ref{defofps}. Pulling these measures back via \(\xi^{1,2}, \xi^{1}\), we obtain regular Borel probability measures on \(\Lambda(\Gamma)\), denoted \(\mu_{\alpha}\) and \(\mu_{\hat{\omega}_{\alpha}}\), respectively. The Radon--Nikodym cocycle relations for these measures are
\[
\frac{d(\gamma_* \mu_{\alpha})}{d\mu_{\alpha}}(x)
= e^{- \mathcal{J}_2\big(B_{1,2}\!\big(\rho'(\gamma^{-1}),\,\xi^{1,2}(x)\big)\big)},
\qquad
\frac{d(\gamma_* \mu_{\hat{\omega}_{\alpha}})}{d\mu_{\hat{\omega}_{\alpha}}}(x)
= e^{- \delta_{\rho'}(\hat{\omega}_1) \hat\omega_1\big( B_{1}\!\big(\rho'(\gamma^{-1}),\,\xi^{1}(x)\big)\big)}.
\]
We set \(I_1=\{1,d-1\}\) and \(I_2=\{1,2,d-2,d-1\}\). Recall that we defined
\[
\xi^{I_1}:\Lambda(\Gamma)\to \mathcal{F}_{I_1}(V),\quad
x\mapsto \big(\xi^1(x),\,\mathrm{Ann}(\xi^1(x))\big),
\]
and, when \(n>2\),
\[
\xi^{I_2}:\Lambda(\Gamma)\to \mathcal{F}_{I_2}(V),\quad
x\mapsto \big(\xi^1(x),\,\ell(x),\,\mathrm{Ann}(\ell(x)),\,\mathrm{Ann}(\xi^1(x))\big).
\]
With this notation, the cocycle identities can equivalently be written as
\[
\frac{d(\gamma_* \mu_{\alpha})}{d\mu_{\alpha}}(x)
= e^{- \mathcal{J}_2\big(B_{I_2}\!\big(\rho'(\gamma^{-1}),\,\xi^{I_2}(x)\big)\big)},
\qquad
\frac{d(\gamma_* \mu_{\hat{\omega}_{\alpha}})}{d\mu_{\hat{\omega}_{\alpha}}}(x)
= e^{- \delta_{\rho'}(\hat{\omega}_1) \hat\omega_1\big( B_{I_1}\!\big(\rho'(\gamma^{-1}),\,\xi^{I_1}(x)\big)\big)},
\]

Fix a basepoint \(b_0 \in \mathbb{D}\), and let \(S := \mathbb{D}/\Gamma\). Denote by \(\mathrm{T}^1 S\) and \(\mathrm{T}^1 \mathbb{D}\) the unit tangent bundles of \(S\) and \(\mathbb{D}\), respectively. Let \(\tilde{\mathcal{G}} \subset \partial \mathbb{D} \times \partial \mathbb{D}\) be the set of distinct pairs. The \emph{Hopf parametrization}
\[
\mathrm{h} \colon \tilde{\mathcal{G}} \times \mathbb{R} \to \mathrm{T}^1 \mathbb{D}
\]
is defined as follows:
\begin{enumerate}
    \item For each \((x, y) \in \tilde{\mathcal{G}}\), the map \(t \mapsto \mathrm{h}(x, y, t)\) parameterizes the geodesic from \(x\) to \(y\) at unit speed.
    \item \(\mathrm{h}(x, y, 0)\) is the point on the geodesic \([x, y]\) lying on the horocycle centered at \(y\) that passes through \(b_0\).
\end{enumerate}

The space $\tilde{\mathcal{G}}\times \mathbb{R}$ carries a $\Gamma\times \mathbb{R}$ action, where $\Gamma$ acts diagonally on $\tilde{\mathcal{G}}$ and $\mathbb{R}$ acts on each geodesic $[x,y]$ by the geodesic flow.

We define the following measures on \(\tilde{\mathcal{G}} \times \mathbb{R}\):
\[
m^*_{\alpha}(x,y,s) := e^{-\mathcal{J}_2([\xi^{I_2}(x),\xi^{I_2}(y)]_{I_2})} \cdot \mathbb{E}_{\mathcal{ET}}(x,y)\, d\mu_{\alpha}(x) \times d\mu_{\alpha}(y) \times dt(s),
\]
\[
m^*_{\hat{\omega}_{\alpha}}(x,y,s) := e^{-\delta_{\rho'}(\hat{\omega}_1) \hat{\omega}_1([\xi^{I_1}(x),\xi^{I_1}(y)]_{I_1})} \, d\mu_{\hat{\omega}_{\alpha}}(x) \times d\mu_{\hat{\omega}_{\alpha}}(y) \times dt(s),
\]
where \(\mathcal{ET} \subset \Lambda(\Gamma) \times \Lambda(\Gamma)\) is the transverse set defined preceding Corollary~\ref{transverse}, and \(\mathbb{E}_{\mathcal{ET}}(x,y)\) is its indicator function. Clearly \(m^*_{\hat{\omega}_{\alpha}}\) is nonzero, and by Corollary~\ref{transverse} and Lemma~\ref{Gromovproduct}, \(m^*_{\alpha}\) is also nonzero.

Note that $\rho'(\Gamma)\subset \Lambda^n(\mathrm{Sp}(2n,\mathbb{R}))$. Using the defining equation in Definition~\ref{iwasawa}, we observe that for any $x\in \Lambda(\Gamma)$ and $k=1,2$,
\[
B_{I_k}(\rho'(\gamma), \xi^{I_k}(x)) = i|_{\mathfrak{a}_{I_k}} \circ B_{I_k}(\rho'(\gamma), \xi^{I_k}(x)),
\]
where \(i\) denotes the involution on \(\mathfrak{a}_d\) given by
\[
i\left(\mathrm{diag}(\lambda_1, \lambda_2, \dots, \lambda_d)\right) = \mathrm{diag}(-\lambda_d, -\lambda_{d-1}, \dots, -\lambda_1).
\]
Using the cocycle identities of Gromov product and Patterson--Sullivan measures, both \(m^*_{\alpha}\) and \(m^*_{\hat{\omega}_{\alpha}}\) are \(\Gamma\)-invariant and invariant under the geodesic flow. They thus push forward via the Hopf parameterization to $\Gamma$ invariant and flow invariant measures on $\mathrm{T}^1\mathbb{D}$, and thus induce flow-invariant Borel measures \(m_{\alpha}\) and \(m_{\hat{\omega}_{\alpha}}\) on \(\mathrm{T}^1 S\), which we call the \emph{Bowen--Margulis--Sullivan measures}.

The relatively-Anosov property in Corollary~\ref{anosov2} implies that the boundary map \( \xi^{I_1} \) is continuous. Consequently, the function
\[
e^{-\hat{\omega}_1([\xi^{I_1}(x),\xi^{I_1}(y)]_{I_1})}
\]
is continuous, and hence both \( m^*_{\hat{\omega}_{\alpha}} \) and \( m_{\hat{\omega}_{\alpha}} \) are locally finite Borel measures.

Moreover,  we have the identity:
\[
-\mathcal{J}_2([\xi^{I_2}(x),\xi^{I_2}(y)]_{I_2}) = (\omega_2 - 2\omega_1)([\xi^{I_2}(x),\xi^{I_2}(y)]_{I_2}) = \omega_2([\xi^{I_2}(x),\xi^{I_2}(y)]_{I_2}) - 2\omega_1([\xi^{I_1}(x),\xi^{I_1}(y)]_{I_1}).
\]
Since \( \omega_2([\xi^{I_2}(x),\xi^{I_2}(y)]_{I_2}) \le 0 \), it follows that
\[
e^{-\mathcal{J}_2([\xi^{I_2}(x),\xi^{I_2}(y)]_{I_2})}
\]
is locally bounded. Therefore, both \( m^*_{\alpha} \) and \( m_{\alpha} \) are also locally finite Borel measures.

Finally, assuming that \(\rho \colon \Gamma \to \mathrm{Sp}(2n, \mathbb{R})\) is a maximal representation from a lattice, we invoke the following ergodicity result:

\begin{theorem}[{Canary–Zhang–Zimmer~\cite[Theorem~1.1]{CaZhZi-relatively}; \cite[Theorem~11.1 and Corollary~12.1]{CaZhZi3}.}]\label{ergodicities}
\leavevmode
\begin{enumerate}
    \item The diagonal action of $\Gamma$ on $\Lambda(\Gamma)^{(2)}$ is ergodic with respect to $\mu_{\hat{\omega}_{\alpha}}\times \mu_{\hat{\omega}_{\alpha}}$. Consequently, the action of $\Gamma$ on $\Lambda(\Gamma)$ is ergodic with respect to $\mu_{\hat{\omega}_{\alpha}}$.
    \item The geodesic flow on $\mathrm{T}^1 S$ is ergodic with respect to $m_{\hat{\omega}_{\alpha}}$.
\end{enumerate}
\end{theorem}

\subsection{The Shadow Lemma}

Fix a base point \( b_0 \in \mathbb{D} \). For each \( z \in \mathbb{D} \) and \( R > 0 \), define the \emph{shadow} \( \mathcal{O}_R(b_0,z) \subset \Lambda(\Gamma) \) by
\[
\mathcal{O}_R(b_0,z) := \{ x \in \Lambda(\Gamma) \mid [b_0,x) \cap \mathbb{B}(z,R) \ne \emptyset \},
\]
where \( [b_0,x) \) denotes the geodesic ray starting from \( b_0 \) toward \( x \in \partial \mathbb{D} \), and \( \mathbb{B}(z, R) \) is the closed hyperbolic ball of radius \( R \) centered at \( z \).

For a maximal representation $\rho:\Gamma\to \mathrm{Sp}(2n,\mathbb{R})$ from a lattice, the Patterson--Sullivan measures satisfy the following shadow lemmas:

\begin{theorem}\label{shadowlemma}

For any \(r\ge 0 \) and sufficiently large $R>0$, there exists a constant \( C > 1 \) such that for all \( z \in \mathbb{D} \), \( \gamma \in \Gamma \) with \( d_{\mathbb{D}}(z, \gamma(b_0)) \le  r \), the following inequalities hold:
\[
\frac{1}{C} e^{-\alpha(\kappa(\rho(\gamma)))} \le \mu_{\alpha}(\mathcal{O}_R(b_0,z)) \le C e^{-\alpha(\kappa(\rho(\gamma)))},
\]
\[
\frac{1}{C} e^{-\delta_{\rho}(\hat\omega_{\alpha})\hat\omega_{\alpha}(\kappa(\rho(\gamma)))} \le \mu_{\hat\omega_{\alpha}}(\mathcal{O}_R(b_0,z)) \le C e^{-\delta_{\rho}(\hat\omega_{\alpha})\hat\omega_{\alpha}(\kappa(\rho(\gamma)))}.
\]
\end{theorem}

For the inequality of $\mu_{\alpha}$, see Yao~\cite[Corollary~4.6]{yzf} (see also Pozzetti--Sambarino--Wienhard~\cite[Lemma~5.15]{Liplim}). For the inequality of $\mu_{\hat{\omega}_{\alpha}}$, see  Canary--Zhang--Zimmer~\cite[Proposition~7.1]{CaZhZi3}.

To end this section, we restate Theorem \ref{shadowlemma} as follows, which is inspired by Dey--Kim--Oh~\cite{dey2024ahlforsregularitypattersonsullivanmeasures}. We continue to use the following notation: $\Gamma$ is a lattice, $d$ is the dimension of $\bigwedge^n \mathbb{R}^{2n}$, $\kappa_d$ is the Cartan projection in $\mathrm{SL}(\bigwedge^n \mathbb{R}^{2n})$, $\rho'$ is $\bigwedge^n\circ \rho$, and $\xi^{I_1},\xi^1$ are the limit maps from $\Lambda(\Gamma)$ to $\mathcal{F}_{I_1}(\bigwedge^n \mathbb{R}^{2n}),\mathcal{F}_1(\bigwedge^n \mathbb{R}^{2n})$.

\begin{proposition}\label{AhlforsRegularity}
For any \(r\ge 0 \) and sufficiently large $R>0$, there exists a constant \( C > 1 \) such that for all \( z \in \mathbb{D} \) and all \( \gamma \in \Gamma \) with \( d_{\mathbb{D}}(b_0,z) > R \) and \( d_{\mathbb{D}}(z, \gamma(b_0)) \le  r \), the following inequalities hold:
\[
\frac{1}{C} \, e^{\delta_{\rho}(\hat\omega_{\alpha}) \frac{1}{2}\hat\omega_1([\xi^{I_1}(x_z),\xi^{I_1}(y_z)]_{I_1})}
\;\le\;
\mu_{\hat\omega_{\alpha}}(\mathcal{O}_R(b_0,z))
\;\le\;
C \, e^{\delta_{\rho}(\hat\omega_{\alpha}) \frac{1}{2}\hat\omega_1([\xi^{I_1}(x_z),\xi^{I_1}(y_z)]_{I_1})},
\]
where $x_z,y_z$ are the endpoints of $\mathcal{O}_{R}(b_0,z)$. Note that whenever $d_{\mathbb{D}}(b_0,z) > R$, $\mathcal{O}_{R}(b_0,z)$ is a proper arc in $S^1$, so $x_z,y_z$ are well-defined.
\end{proposition}

\begin{proof}
Fix $r\ge 0$ and a sufficiently large $R>0$. From Theorem~\ref{shadowlemma}, it suffices to find $C>1$ such that
\[
\frac{1}{C} \, e^{\delta_{\rho}(\hat\omega_{\alpha})\frac{1}{2}\hat\omega_1([\xi^{I_1}(x_z),\xi^{I_1}(y_z)]_{I_1})}
\;\le\;
e^{-\delta_{\rho}(\hat\omega_{\alpha}) \hat\omega_{\alpha}(\kappa(\rho(\gamma)))}
\;\le\;
C \, e^{\delta_{\rho}(\hat\omega_{\alpha})\frac{1}{2}\hat\omega_1([\xi^{I_1}(x_z),\xi^{I_1}(y_z)]_{I_1})}.
\]
By the definition of $\hat{\omega}_1$, this is equivalent to the existence of $C>1$ such that
\[
\frac{1}{C} \, e^{-\frac{1}{2}\omega_1([\xi^{I_1}(x_z),\xi^{I_1}(y_z)]_{I_1})}
\;\le\;
e^{\omega_{1}(\kappa_d(\rho'(\gamma)))}
\;\le\;
C \, e^{-\frac{1}{2}\omega_1([\xi^{I_1}(x_z),\xi^{I_1}(y_z)]_{I_1})},
\]
where $\rho' := \Lambda^n \circ \rho$.

We first prove the upper bound
\[
e^{\omega_{1}(\kappa_d(\rho'(\gamma)))} \le
C \, e^{-\frac{1}{2}\omega_1([\xi^{I_1}(x_z),\xi^{I_1}(y_z)]_{I_1})}.
\]
The proof of the lower bound is identical.

Suppose that no such $C$ exists. Then there are sequences $\gamma_n \in \Gamma$ and $z_n \in \mathbb{D}$ with
\[
d_{\mathbb{D}}(b_0,z_n) > R, \qquad d_{\mathbb{D}}(z_n,\gamma_n(b_0)) \le r,
\]
and
\[
e^{\omega_{1}(\kappa_d(\rho'(\gamma_n)))} \ge
n \, e^{-\frac{1}{2}\omega_1([\xi^{I_1}(x_{z_n}),\xi^{I_1}(y_{z_n})]_{I_1})}.
\]

Note that $\gamma_n$ is unbounded: otherwise $z_n$ will be a bounded sequence in $\mathbb{D}$, so $(x_{z_n}, y_{z_n})$ will not collapse to each other, so $(\xi^{I_1}(x_{z_n}),\xi^{I_1}(y_{z_n}))$ is uniformly transverse to each other, so  $e^{-\frac{1}{2}\omega_1([\xi^{I_1}(x_{z_n}),\xi^{I_1}(y_{z_n})]_{I_1})}$ is bounded below, so the right-hand side of this inequality is unbounded, contradicting the boundedness of the left-hand side. Let
\[
z_n' = \gamma_n^{-1}(z_n), \quad
x'_n = \gamma_n^{-1}(x_{z_n}), \quad
y'_n = \gamma_n^{-1}(y_{z_n}),
\]
Note that $d_{\mathbb{D}}(z'_n,b_0) = d_{\mathbb{D}}(z_n,\gamma_n(b_0)) \le  r$, so, after passing to a subsequence, we may assume
\[
\gamma_n^{-1}(b_0)\to b'\in \Lambda(\Gamma),\quad  x'_n \to x' \in \Lambda(\Gamma), \qquad
y'_n \to y' \in \Lambda(\Gamma), \qquad
z'_n \to z' \in \mathbb{D}.
\]
Clearly $x',y'$ are the endpoints of $\mathcal{O}_R(b',z')$, so they are distinct points.

By the cocycle property of the Gromov product and the symmetry of the Cartan subalgebra $\mathfrak{a} \subset \mathfrak{sp}(2n,\mathbb{R})$, we have
\begin{align*}
e^{-\frac{1}{2}\omega_1([\xi^{I_1}(x_{z_n}),\xi^{I_1}(y_{z_n})]_{I_1})}
&=
e^{-\frac{1}{2}\omega_1([\rho'(\gamma_n)\xi^{I_1}(x'_n),\rho'(\gamma_n)\xi^{I_1}(y'_n)]_{I_1})} \\
&=
e^{-\frac{1}{2}\omega_1([\xi^{I_1}(x'_n),\xi^{I_1}(y'_n)]_{I_1})
+\frac{1}{2}\omega_1( B_{I_1}(\rho'(\gamma_n), \xi^{I_1}(x'_n)))
+\frac{1}{2}\omega_1( B_{I_1}(\rho'(\gamma_n), \xi^{I_1}(y'_n)))} \\
&=
e^{-\frac{1}{2}\omega_1([\xi^{I_1}(x'_n),\xi^{I_1}(y'_n)]_{I_1})}
\cdot
e^{\frac{1}{2}\omega_1( B_{I_1}(\rho'(\gamma_n), \xi^{I_1}(x'_n)) )}
\cdot
e^{\frac{1}{2}\omega_1( B_{I_1}(\rho'(\gamma_n), \xi^{I_1}(y'_n)) )}.
\end{align*}
Recall from Corollary~\ref{anosov2} that $\rho'$ is relatively $I_1$–Anosov, so $\xi^{I_1}$ is transverse. Since $x'_n,y'_n$ converge to distinct points, the term
\[
e^{-\frac{1}{2}\omega_1([\xi^{I_1}(x'_n),\xi^{I_1}(y'_n)]_{I_1})}
\]
converges to a positive number \(e^{-\frac{1}{2}\omega_1([\xi^{I_1}(x'),\xi^{I_1}(y')]_{I_1})}>0\). From Canary–Zhang–Zimmer~\cite[Lemma~7.3]{CaZhZi3}, there exists a constant $C_0>0$ such that for any $n$ and any $p \in \mathcal{O}_{R}(b_0,z_n)$,
\[
\big|
\omega_1\big(
B_{I_1}(\rho'(\gamma_n)^{-1},\xi^{I_1}(p))
+
\pi_{I_1}\circ\kappa_{d}(\rho'(\gamma_n))
\big)
\big|
\le C_0.
\]
Note that
\[
\omega_1(\pi_{I_1}\circ\kappa_{d}(\rho'(\gamma_n))) = \omega_1(\kappa_{d}(\rho'(\gamma_n))), 
\qquad
B_{I_1}(\rho'(\gamma_n)^{-1},\xi^{I_1}(p))
= -B_{I_1}(\rho'(\gamma_n), \rho'(\gamma_n)^{-1}\xi^{I_1}(p)),
\]
where the latter equality is due to the cocycle identity of $B_{I_1}$. 

Substituting $p=x_{z_n}$ and $p=y_{z_n}$, we obtain
\[
\big|
\omega_1\big(
-B_{I_1}(\rho'(\gamma_n),\xi^{I_1}(x'_n))
+
\kappa_{d}(\rho'(\gamma_n))
\big)
\big|
\le C_0,
\qquad
\big|
\omega_1\big(
-B_{I_1}(\rho'(\gamma_n),\xi^{I_1}(y'_n))
+
\kappa_{d}(\rho'(\gamma_n))
\big)
\big|
\le C_0.
\]
Therefore, there exists a constant $C_1>1$ such that
\[
\frac{1}{C_1} \, e^{\omega_1(\kappa_d(\rho'(\gamma_n)))} 
\le
e^{\, \omega_1(B_{I_1}(\rho'(\gamma_n), \xi^{I_1}(x'_n)))} 
\le
C_1 \, e^{\omega_1(\kappa_d(\rho'(\gamma_n)))},
\]
\[
\frac{1}{C_1} \, e^{\omega_1(\kappa_d(\rho'(\gamma_n)))} 
\le
e^{\, \omega_1(B_{I_1}(\rho'(\gamma_n), \xi^{I_1}(y'_n)))} 
\le
C_1 \, e^{\omega_1(\kappa_d(\rho'(\gamma_n)))}.
\]
Combining these bounds with the factorization above,
\[
e^{-\frac{1}{2}\omega_1([\xi^{I_1}(x_{z_n}),\xi^{I_1}(y_{z_n})]_{I_1})}
=
e^{-\frac{1}{2}\omega_1([\xi^{I_1}(x'_n),\xi^{I_1}(y'_n)]_{I_1})}
\cdot
e^{\frac{1}{2}\omega_1(B_{I_1}(\rho'(\gamma_n),\xi^{I_1}(x'_n)))}
\cdot
e^{\frac{1}{2}\omega_1(B_{I_1}(\rho'(\gamma_n),\xi^{I_1}(y'_n)))},
\]
we obtain a constant $C_2>1$ such that
\[
\frac{1}{C_2} \, e^{\omega_{1}(\kappa_d(\rho'(\gamma_n)))} 
\le
e^{-\frac{1}{2}\omega_1([\xi^{I_1}(x_{z_n}),\xi^{I_1}(y_{z_n})]_{I_1})}
\le
C_2 \, e^{\omega_{1}(\kappa_d(\rho'(\gamma_n)))}.
\]
This contradicts the assumed failure of the desired bound. The proof of the other inequality is identical.
\end{proof}

\section{Proof of Theorem \ref{main}}\label{mainproof}

We now prove Theorem \ref{main}. We have already observed the ``if'' direction in the introduction, so we focus on the ``only if'' direction.

We first follow Canary--Zhang--Zimmer~\cite[Proof of Lemma~13.6]{CaZhZi3} to establish a lemma, which shows that a local comparison of the measures on shadows implies a global comparison of the measures. Similar arguments also appeared in Sullivan~\cite[Section~8]{sullivanacta} and Burger~\cite[Proof of Theorem~1(a)]{burgermanhattan}.

For a fixed $b_0\in \mathbb{D}$, a lattice $\Gamma\subset \mathrm{PSL}(2,\mathbb{R})$, and a Borel measure $\mu$ on $\Lambda(\Gamma)$, we say $R>0$ is \emph{adequate} with respect to $(b_0,\Gamma,\mu)$ if the family of shadows $\{\mathring{\mathcal{O}}_R(b_0,\gamma(b_0))\}_{\gamma\in\Gamma}$ covers $\mu$-almost every $x\in \Lambda(\Gamma)$ with infinite multiplicity.

For a fixed lattice $\Gamma\subset \mathrm{PSL}(2,\mathbb{R})$, we say a Borel measure $\mu$ on $\Lambda(\Gamma)$ has the \emph{$\mathcal{S}$-property} if for each fixed $b_0\in \mathbb{D}$, sufficiently large $R>0$, and $A\ge 1$, there exists a constant $C = C(R,A)>1$ such that for any $\gamma\in \Gamma$,
\begin{equation}\label{Spropertyholds}
\frac{1}{C(R,A)} \mu(\mathcal{O}_{AR}(b_0,\gamma(b_0))) \le \mu(\mathcal{O}_{R}(b_0,\gamma(b_0)))\le  C(R,A) \mu(\mathcal{O}_{AR}(b_0,\gamma(b_0))).    
\end{equation}

\begin{remark}\label{adequateandSproperty}
Note that if \( \rho:\Gamma \to \mathrm{Sp}(2n,\mathbb{R}) \) is a maximal representation of a lattice \( \Gamma \), then both $\mu_{\alpha}$ and $\mu_{\hat{\omega}_{\alpha}}$ are outer regular, atomless Borel measures (see Proposition~\ref{defofjac} and Proposition~\ref{defofps}). They also satisfy the $\mathcal{S}$-property, which is a direct consequence of Theorem~\ref{shadowlemma}.

Moreover, we can find $R$ adequate with respect to both $(b_0, \Gamma,\mu_{\alpha})$ and $(b_0,\Gamma, \mu_{\hat{\omega}_{\alpha}})$ as follows. The \emph{conical limit set} $\Lambda_c(\Gamma)\subset \Lambda(\Gamma)$ is defined as
\[
\Lambda_c(\Gamma) = \{x\in \Lambda(\Gamma)\mid \exists (\gamma_k) \subset \Gamma, \exists A>0 \text{ s.t. } \gamma_k(b_0)\to x \text{ and } d_{\mathbb{D}}(\gamma_k(b_0), \overrightarrow{b_0x})\le A\},
\]
where $\overrightarrow{b_0x}$ denotes the geodesic ray from $b_0$ to $x$. Since $\Gamma$ is a lattice, the set of points not in the conical limit set, i.e.,  $\Lambda(\Gamma) \setminus \Lambda_c(\Gamma)$, consists precisely of the parabolic fixed points, which is countable. Since $\mu_{\alpha}, \mu_{\hat{\omega}_{\alpha}}$ are atomless, $\Lambda(\Gamma) \setminus \Lambda_c(\Gamma)$  has measure zero. Therefore, $\Lambda_c(\Gamma)$ has full measure with respect to both $\mu_{\alpha}$ and $\mu_{\hat{\omega}_{\alpha}}$.

For any fixed $b_0\in \mathbb{D}$, we choose $R>0$ large enough so that $\mathbb{B}_{\mathbb{D}}(b_0,R)$ covers the compact part of a fundamental domain obtained by excluding the cuspidal neighborhoods. Every geodesic ray ending at $\Lambda_c(\Gamma)$ must intersect the $\Gamma$-translates of the compact part infinitely many times, so this choice ensures that the family of shadows $\{\mathring{\mathcal{O}}_R(b_0,\gamma(b_0))\}_{\gamma\in\Gamma}$ covers every point in $\Lambda_c(\Gamma)$ infinitely often. Since $\Lambda_c(\Gamma)$ has full measure with respect to both $\mu_{\alpha}$ and $\mu_{\hat{\omega}_{\alpha}}$, this $R$ is simultaneously adequate with respect to $(b_0,\Gamma, \mu_{\alpha})$ and $(b_0,\Gamma, \mu_{\hat{\omega}_{\alpha}})$.
\end{remark}

\begin{lemma}\label{vitaliargument}
Let $\Gamma\subset \mathrm{PSL}(2,\mathbb{R})$ be a lattice. Let $\mu$ and $\nu$ be Borel measures on $\Lambda(\Gamma)$, and assume that $\mu$ has the $\mathcal{S}$-property and $\nu$ is outer regular. Fix $b_0 \in \mathbb{D}$ and $R>0$ adequate with respect to $(b_0,\Gamma,\mu)$ such that (\ref{Spropertyholds}) holds. Suppose there exists a constant $C>0$ such that for $\mu$-almost every $x \in \Lambda(\Gamma)$, if $(\gamma_k) \subset \Gamma$ is a sequence such that $\gamma_k(b_0) \to x$ and $x\in\mathcal{O}_{R}(b_0,\gamma_k(b_0)) $, then for all sufficiently large $k$,
\[
\mu\bigl(\mathcal{O}_{R}(b_0,\gamma_k(b_0))\bigr)
\;\le\;
C \, \nu\bigl(\mathcal{O}_{R}(b_0,\gamma_k(b_0))\bigr).
\]
Then there exists a constant $C' = C'(C)$ such that
\[
\mu \le C' \, \nu.
\]
Moreover, $\lim_{C \to 0} C'(C) = 0$.
\end{lemma}
\begin{proof}
Let $E\subset \Lambda(\Gamma)$ be a measurable subset. Since $\nu$ is outer regular, for any $\epsilon>0$, we can choose an open set $U \supset E$ such that $\nu(U)<\nu(E)+\epsilon$.

Let $\mathcal{I}$ be the set of $\gamma \in \Gamma$ such that the shadow $\mathcal{O}_{R}(b_0,\gamma(b_0))$ satisfies both conditions:
\[
\mathcal{O}_{R}(b_0,\gamma(b_0))\subset U \quad \text{and} \quad \mu\bigl(\mathcal{O}_{R}(b_0,\gamma(b_0))\bigr)
\;\le\;
C \, \nu\bigl(\mathcal{O}_{R}(b_0,\gamma(b_0))\bigr).
\]
Since $R$ is adequate and by the lemma's hypothesis, for almost every $x\in E$, there exists $\gamma\in \mathcal{I}$ such that $x\in\mathring{\mathcal{O}}_{R}(b_0,\gamma(b_0))$. As a result, the family of open shadows $\{\mathring{\mathcal{O}}_{R}(b_0,\gamma(b_0))\}_{\gamma\in \mathcal{I}}$ forms an open covering of a subset $E' \subset E$ with $\mu(E \setminus E') = 0$.

By a Vitali-type covering argument (see Roblin~\cite[p. 23]{shadowvitali}), there exists a subfamily $\mathcal{J}\subset \mathcal{I}$ such that $\{\mathcal{O}_R(b_0,\gamma(b_0))\mid \gamma \in \mathcal{J}\}$ is a disjoint family and
\[
\bigcup_{\gamma\in \mathcal{I}}\mathcal{O}_{R}(b_0,\gamma(b_0))\subset \bigcup_{\gamma\in \mathcal{J}}\mathcal{O}_{5R}(b_0,\gamma(b_0)).
\]
We can now estimate $\mu(E)$:
\begin{align*}
\mu(E) &= \mu(E') \le \mu\left(\bigcup_{\gamma\in \mathcal{I}}\mathcal{O}_{R}(b_0,\gamma(b_0))\right)
\\
&\le \mu\left(\bigcup_{\gamma\in \mathcal{J}}\mathcal{O}_{5R}(b_0,\gamma(b_0))\right)
\\
&\le \sum_{\gamma\in \mathcal{J}} \mu(\mathcal{O}_{5R}(b_0,\gamma(b_0)))
&& \text{(by subadditivity)}
\\
&\le C(R,5) \sum_{\gamma\in \mathcal{J}} \mu (\mathcal{O}_{R}(b_0,\gamma(b_0)))
&& \text{(by the $\mathcal{S}$-property of $\mu$)}
\\
&\le C(R,5) C \sum_{\gamma\in \mathcal{J}} \nu (\mathcal{O}_{R}(b_0,\gamma(b_0)))
&& \text{(by definition of $\mathcal{I}$)}
\\
&= C(R,5) C \, \nu \left(\bigcup_{\gamma \in \mathcal{J}} \mathcal{O}_{R}(b_0,\gamma(b_0))\right)
&& \text{(by the disjointness property of $\mathcal{J}$)}
\\
&\le C(R,5) C \, \nu (U)
&& \text{(since $\mathcal{O}_R(b_0,\gamma(b_0)) \subset U$ for $\gamma \in \mathcal{J}$)}
\\
&< C(R,5) C (\nu (E)+\epsilon)
&& \text{(by choice of $U$)}.
\end{align*}
Since $\epsilon>0$ was arbitrary, we have $\mu(E) \le C(R,5)C \, \nu(E)$. As $E$ was an arbitrary measurable set, we conclude that $\mu \le C' \nu$ with $C' = C(R,5)C$. It follows directly that $\lim_{C \to 0} C'(C) = 0$.
\end{proof}

As an application, we have the following:
\begin{lemma}\label{domination}
Assume \( \rho:\Gamma \to \mathrm{Sp}(2n,\mathbb{R}) \) is a maximal representation of a lattice \( \Gamma \) and \( \delta_{\rho}(\hat{\omega}_{\alpha}) = 1 \). Then there exists a constant \( C > 0 \) such that
\[
\mu_{\hat{\omega}_{\alpha}} \le C\mu_{\alpha} .
\]
\end{lemma}

\begin{proof}
Fix \( b_0 \in \mathbb{D} \). Pick a sufficiently large $R>0$ that is adequate with respect to $(b_0,\Gamma, \mu_{\hat{\omega}_{\alpha}})$. Theorem~\ref{shadowlemma} implies that there exists a constant \( C_1 > 1 \) such that for all \( \gamma \in \Gamma \)  we have
\[
\frac{1}{C_1} \, e^{-\alpha(\kappa(\rho(\gamma)))}
\;\le\;
\mu_{\alpha}\bigl(\mathcal{O}_R(b_0, \gamma(b_0))\bigr)
\;\le\;
C_1 \, e^{-\alpha(\kappa(\rho(\gamma)))} ,
\]
\[
\frac{1}{C_1} \, e^{-\hat{\omega}_\alpha(\kappa(\rho(\gamma)))}
\;\le\;
\mu_{\hat{\omega}_\alpha}\bigl(\mathcal{O}_R(b_0, \gamma(b_0))\bigr)
\;\le\;
C_1 \, e^{-\hat{\omega}_\alpha(\kappa(\rho(\gamma)))} .
\]
Since
\(
\hat{\omega}_\alpha(\kappa(\rho(\gamma))) \ge \alpha(\kappa(\rho(\gamma)))
\)
for all \( \gamma \in \Gamma \), it follows that
\[
e^{-\alpha(\kappa(\rho(\gamma)))}
\;\ge\;
e^{-\hat{\omega}_\alpha(\kappa(\rho(\gamma)))} .
\]
Combining this with the shadow estimates above, we obtain a constant \( C_2 > 0 \) such that for \( \gamma\in \Gamma \),
\[
\mu_{\hat{\omega}_{\alpha}}(\mathcal{O}_R(b_0,\gamma(b_0))) \le C_2\mu_{\alpha}(\mathcal{O}_R(b_0,\gamma(b_0))).
\]
Applying Lemma~\ref{vitaliargument} completes the proof.
\end{proof}

Next we define the rank function on \( \Lambda(\Gamma) \), which turns out to be closely related to the $\delta_{\rho}(\hat{\omega}_{\alpha}) = 1$ condition.

Let \( \rho:\Gamma\to \mathrm{Sp}(2n,\mathbb{R}) \) be a maximal representation of a lattice, and set \( \rho' = \Lambda^n\circ \rho \), with associated \(1\)-limit map \( \xi^1 \). For any \( x \in \Lambda(\Gamma) \) such that \( \xi^{1}(\Lambda(\Gamma)) \) is differentiable at \( \xi^{1}(x) \), pick \( y,z \in \Lambda(\Gamma) \) such that \( \{x,y,z\} \) are distinct. Choose an arc \( [y,z] \subset \Lambda(\Gamma)\) joining \( y \) to \( z \) with \( x \in (y,z) \), let \( L = \mu_{\alpha}([y,z]) \), and let
\[
p:[0,L]\to [y,z]
\]
be the orientation-preserving arc-length parametrization, which means that
\[
p(0) = y,\quad p(L)=z,\quad \text{and}\quad \mu_{\alpha}(p([0,t])) = t,\quad \forall t\in [0,L].
\]
From Lemma~\ref{lemma:parameterizationofmaps_combined} and Remark~\ref{remark:variationofparameterization}, there exist \( g\in \mathrm{Aut}(\mathfrak{sp}(2n,\mathbb{R})) \) and a continuous map \( f:[0,L)\to \mathrm{Pos}(n)\cup \{0\} \) satisfying
\[
\xi^{1}(y,p(t),z)
=
g\Bigl(
\mathrm{Span}(e_1\wedge\cdots\wedge e_n),\;
U_{f(t)} \mathrm{Span}( e_1\wedge\cdots\wedge e_n),\;
\mathrm{Span}(e_{n+1}\wedge\cdots\wedge e_{2n})
\Bigr),
\]
where
\[
U_{f(t)} = \begin{pmatrix} \mathrm{Id}_n & 0 \\ f(t) & \mathrm{Id}_n \end{pmatrix},
\quad\text{and whenever } t>s,\ \ f(t)-f(s) \in \mathrm{Pos}(n).
\]

Suppose \( p(t_0) = x \). Note that for almost every $x$,  \( f(t) \) is differentiable at \( t_0 \) with \( f'(t_0)\in \overline{\mathrm{Pos}(n)}\setminus\{0\} \).

\begin{definition}\label{defofrank}
In the above setting, if \( x\in \Lambda(\Gamma) \) is a point such that the arc length parameterization of \( \xi^1(\Lambda(\Gamma)) \) is differentiable at \( \xi^1(x) \), then the \emph{rank} of \( x \), denoted \( r(x) \), is defined as the matrix rank of \( f'(t_0) \).
\end{definition}

It is direct to check that \( r(x) \) is well defined: it is independent of the choice of \( y,z \), and independent of the choice of \( g \) in the parametrization. Note that the function \( r(x) \) is defined for \( \mu_{\alpha} \)-almost every \( x \). For every $k\in\{1,\dots,n\}$, let \( R_k \subset \Lambda(\Gamma) \) be the set of points of rank \( k \). Recall that $\mu_{\alpha}$ and $\mu_{\hat{\omega}_{\alpha}}$ are the pullbacks of Patterson--Sullivan measures supported on $\Lambda(\Gamma)$ (see Section~\ref{bmsdef}).

\begin{proposition}\label{lebesgueergodicity}
The following are equivalent:
\begin{enumerate}
    \item $\delta_{\rho}(\hat{\omega}_{\alpha}) = 1$.
    \item $\mu_{\alpha}(R_n)>0$.
    \item $\mu_{\alpha}$ is an ergodic measure and $R_n$ is of full measure.
\end{enumerate}
\end{proposition}

\begin{proof}

Before we begin, we prove a technical lemma.  Recall that \(\|\cdot\|_2\) denotes the Frobenius norm defined by $\|M\|_2 = \sqrt{\mathrm{tr}(M^TM)}$. 

\begin{lemma}\label{leminacute}
For any $x\not= y\in \Lambda(\Gamma)$, let $I\subset \Lambda(\Gamma)$ be a compact neighbourhood of $x$ not containing $y$. Applying Lemma~\ref{lemma:parameterizationofmaps_combined} and Remark~\ref{remark:variationofparameterization}, pick \( g\in \mathrm{Aut}(\mathfrak{sp}(2n,\mathbb{R})) \) and a map
\(
\tilde{f}:I\to \mathrm{Pos}(n)\cup\{0\}\cup \bigl(-\mathrm{Pos}(n)\bigr)
\)
such that for any \( z\in I \) we have
\[
\xi^1(x,z,y)
=
g\bigl(
[e_1\wedge\cdots\wedge e_n],\;
U_{\tilde{f}(z)}[e_1\wedge\cdots\wedge e_n],\;
[e_{n+1}\wedge \cdots\wedge e_{2n}]
\bigr),
\]
where
\(
U_{\tilde{f}(z)} = \begin{pmatrix} \mathrm{Id}_n & 0 \\ \tilde{f}(z) & \mathrm{Id}_n \end{pmatrix},
\) and \( \tilde{f} \) satisfies:
\[
\tilde{f}(x) = 0, \qquad
\tilde{f}(z_1)-\tilde{f}(z_2)\in \mathrm{Pos}(n)
\ \text{whenever } y,z_2,z_1 \text{ are in clockwise order.}
\]
Fix a sufficiently large $R$, then there exists a constant $C>1$ depending on the above choices: For any \( \gamma \) such that \( \mathcal{O}_R(b_0,\gamma(b_0))\subset I \), let \( p_{\gamma} \) and \( q_{\gamma} \) be the endpoints of \( \mathcal{O}_{R}(b_0,\gamma(b_0)) \), ordered such that \( y, p_{\gamma}, q_{\gamma} \) are in clockwise order. We have
\[
\frac{1}{C} \,\|\tilde{f}(q_{\gamma})-\tilde{f}(p_{\gamma}) \|_{2}
\;\le\;
\mu_{\alpha}\bigl(\mathcal{O}_{R}(b_0,\gamma(b_0))\bigr)
\;\le\;
C \,\|\tilde{f}(q_{\gamma})-\tilde{f}(p_{\gamma}) \|_{2}.
\]
\[
\frac{1}{C}\,
\bigl(\det (\tilde{f}(q_{\gamma})-\tilde{f}(p_{\gamma}))\bigr)^{\frac{\delta_{\rho}(\hat{\omega}_{\alpha})}{n}}
\;\le\;
\mu_{\hat{\omega}_{\alpha}}\bigl(\mathcal{O}_{R}(b_0,\gamma(b_0))\bigr)
\;\le\;
C \,
\bigl(\det (\tilde{f}(q_{\gamma})-\tilde{f}(p_{\gamma}))\bigr)^{\frac{\delta_{\rho}(\hat{\omega}_{\alpha})}{n}}.
\]

\end{lemma}

\begin{proof}
     Note that $\tilde{f}(I)$ is parameterized monotonically in $\mathrm{Pos}(n)$, which implies that for any $p,q\in I$, the length of the arc in $\tilde{f}(I)$ connecting $\tilde{f}(p), \tilde{f}(q)$ is coarsely $\| \tilde{f}(p)-\tilde{f}(q)\|_{2}$. Since \( \mu_{\alpha} \) measures the arc length of \( \xi^1(\Lambda(\Gamma)) \), and $\tilde{f}(I)$ is bi-Lipschitz equivalent with $\xi^1(I)$, there exists a constant \( C_1>1 \) such that for any \( \gamma \) satisfying \( \mathcal{O}_R(b_0,\gamma(b_0))\subset I \), we have
\[
\frac{1}{C_1} \,\|\tilde{f}(q_{\gamma})-\tilde{f}(p_{\gamma}) \|_{2}
\;\le\;
\mu_{\alpha}\bigl(\mathcal{O}_{R}(b_0,\gamma(b_0))\bigr)
\;\le\;
C_1 \,\|\tilde{f}(q_{\gamma})-\tilde{f}(p_{\gamma}) \|_{2}.
\]

Recall that we defined the bilinear form $\omega'$ on $\bigwedge^{n}\mathbb{R}^{2n}$ by
\[
\omega'(v,w)=\frac{v\wedge w}{e_1\wedge e_2\wedge \cdots \wedge e_{2n}}.
\]
And as before we denote $I_1 = \{1, \dim \bigwedge^{n}\mathbb{R}^{2n} -1\}$. Using the same type of calculation as in Lemma~\ref{Gromovproduct}, we obtain the following for any \( \gamma \) with \( \mathcal{O}_R(b_0,\gamma(b_0))\subset I \):
\[
\begin{aligned}
e^{\omega_1 ([\xi^{I_1}(q_{\gamma}),\xi ^{I_1}(p_{\gamma})]_{I_1})}
&=
\frac{|\omega'\left(U_{\tilde{f}(p_{\gamma})}(e_1\wedge\cdots\wedge e_n),\,
U_{\tilde{f}(q_{\gamma})}(e_1\wedge\cdots\wedge e_n)\right)|}
{\|gU_{\tilde{f}(p_{\gamma})}(e_1\wedge\cdots\wedge e_n)\|\cdot
\|gU_{\tilde{f}(q_{\gamma})}(e_1\wedge\cdots\wedge e_n)\|}
\\
&=
\frac{\det \bigl(\tilde{f}(q_{\gamma})-\tilde{f}(p_{\gamma})\bigr)}
{\|gU_{\tilde{f}(p_{\gamma})}(e_1\wedge\cdots\wedge e_n)\|\cdot
\|gU_{\tilde{f}(q_{\gamma})}(e_1\wedge\cdots\wedge e_n)\|}.
\end{aligned}
\]

Since \( \tilde{f} \) is continuous and $I$ is compact, there exists \( C_2>1 \) such that for any \( \gamma \) with \( \mathcal{O}_R(b_0,\gamma(b_0))\subset I \), we have
\[
\frac{1}{C_2}\,\det \bigl(\tilde{f}(q_{\gamma})-\tilde{f}(p_{\gamma})\bigr)
\;\le\;
e^{\omega_1 ([\xi^{I_1}(q_{\gamma}),\xi ^{I_1}(p_{\gamma})]_{I_1})}
\;\le\;
C_2 \,\det \bigl(\tilde{f}(q_{\gamma})-\tilde{f}(p_{\gamma})\bigr).
\]
By Proposition~\ref{AhlforsRegularity}, there exists \( C_3>1 \) such that for any \( \gamma \) with \( \mathcal{O}_R(b_0,\gamma(b_0))\subset I \), we have
\[
\frac{1}{C_3}\,
\bigl(\det (\tilde{f}(q_{\gamma})-\tilde{f}(p_{\gamma}))\bigr)^{\frac{\delta_{\rho}(\hat{\omega}_{\alpha})}{n}}
\;\le\;
\mu_{\hat{\omega}_{\alpha}}\bigl(\mathcal{O}_{R}(b_0,\gamma(b_0))\bigr)
\;\le\;
C_3 \,
\bigl(\det (\tilde{f}(q_{\gamma})-\tilde{f}(p_{\gamma}))\bigr)^{\frac{\delta_{\rho}(\hat{\omega}_{\alpha})}{n}}.
\]
Pick $C = \max(C_1, C_3)$, we have finished the proof.
\end{proof}

Now we return to the proof of Proposition~\ref{lebesgueergodicity}.

$(1)\implies (2):$ We prove this by contradiction. Suppose $\delta_{\rho}(\hat{\omega}_{\alpha}) = 1$ but \( r(x) < n \) for \( \mu_{\alpha} \)-almost all \( x \). We saw previously that since \( \Gamma \) is a lattice, there exists \( R>0 \) sufficiently large such that for \( \mu_{\alpha} \)-almost all \( x \), there exist infinitely many \( \gamma \in \Gamma \) with \( x\in \mathcal{O}_R(b_0,\gamma(b_0)) \). Hence for \( \mu_{\alpha} \)-almost all \( x \), both \( r(x) < n \) (which implicitly indicates that $\xi^1(\Lambda(\Gamma))$ is differentiable at $\xi^1(x)$) and \( x\in \mathcal{O}_R(b_0,\gamma(b_0)) \) for infinitely many \( \gamma\in \Gamma \) hold.

Fix such an \( x \), choose $y,I,g,\tilde{f}$ as we did in Lemma~\ref{leminacute}, and from the same lemma, there exists a constant $C>1$ such that for any \( \gamma \) with \( \mathcal{O}_R(b_0,\gamma(b_0))\subset I \), we have
\[
\frac{1}{C} \,\|\tilde{f}(q_{\gamma})-\tilde{f}(p_{\gamma}) \|_{2}
\;\le\;
\mu_{\alpha}\bigl(\mathcal{O}_{R}(b_0,\gamma(b_0))\bigr)
\;\le\;
C \,\|\tilde{f}(q_{\gamma})-\tilde{f}(p_{\gamma}) \|_{2}.
\]
\[
\frac{1}{C}\,
\bigl(\det (\tilde{f}(q_{\gamma})-\tilde{f}(p_{\gamma}))\bigr)^{\frac{1}{n}}
\;\le\;
\mu_{\hat{\omega}_{\alpha}}\bigl(\mathcal{O}_{R}(b_0,\gamma(b_0))\bigr)
\;\le\;
C \,
\bigl(\det (\tilde{f}(q_{\gamma})-\tilde{f}(p_{\gamma}))\bigr)^{\frac{1}{n}}.
\]
where $p_{\gamma}, q_{\gamma}$ are also defined in Lemma~\ref{leminacute}. Note that \( C\) is independent of \( \gamma \), but may depend on \( g \), \( x \), \( y \), and \( I \). Also note that we have substituted $\delta_{\rho}(\hat{\omega}_{\alpha}) = 1$.

Now choose a diverging sequence \( \gamma_k\in \Gamma \) such that \( q_{\gamma_k},p_{\gamma_k}\to x \) from different sides. Since \( r(x)<n \),
\[
\frac{\tilde{f}(q_{\gamma_k})-\tilde{f}(p_{\gamma_k})}{\|\tilde{f}(q_{\gamma_k})-\tilde{f}(p_{\gamma_k}) \|_{2}}
\text{ converges to a degenerate semidefinite matrix},
\]
so
\[
\frac{\det (\tilde{f}(q_{\gamma_k})-\tilde{f}(p_{\gamma_k}))}{\|\tilde{f}(q_{\gamma_k})-\tilde{f}(p_{\gamma_k}) \|_{2}^{\,n}}
\;\longrightarrow\; 0.
\]
After taking the $n$-th root, we get
\[
\frac{\mu_{\hat{\omega}_{\alpha}}(\mathcal{O}_{R}(b_0,\gamma_k(b_0)))}{\mu_{\alpha}(\mathcal{O}_{R}(b_0,\gamma_k(b_0)))}\le C^2 \frac{\big(\det (\tilde{f}(q_{\gamma_k})-\tilde{f}(p_{\gamma_k}))\big)^{\frac{1}{n}}}{\|\tilde{f}(q_{\gamma_k})-\tilde{f}(p_{\gamma_k}) \|_{2}}
\;\longrightarrow\; 0.
\]
Since this holds for \( \mu_{\alpha} \)-almost all \( x \), we see that the Radon--Nikodym derivative \( \frac{d\mu_{\hat{\omega}_{\alpha}}}{d\mu_{\alpha}} \) is zero \( \mu_{\alpha} \)-almost everywhere, so \( \mu_{\hat{\omega}_{\alpha}} = 0 \), which is impossible. Therefore, if \( \delta_{\rho}(\hat{\omega}_{\alpha}) = 1 \), then \(\mu_{\alpha}(R_n) >0\).

$(2)\implies(3)$: If \( \mu_{\alpha}(R_n) > 0 \), then any \( \Gamma \)-invariant subset of positive measure \( A \subset R_n \)  has dense support. Using the notation preceding to Definition~\ref{defofrank},  it is direct that for any fixed $y,z$, and any $0<t_1<t_2<L$, if we let $I_{t_1,t_2} = \{t\in[t_1,t_2]\mid gU_{f(t)}\mathrm{Span}(e_1\wedge \dots e_n)\in \xi^{1}(A)\}$, then $\int_{I_{t_1,t_2}} f'(t)dt\in \mathrm{Pos}(n)$, so $A$ is \( \alpha \)-acute (see Yao~\cite[Definition~4.11]{yzf}), and, as a consequence of \cite[Proposition~4.14]{yzf}, such an \( A \) has full \( \mu_{\alpha} \)-measure. In this case, \( \mu_{\alpha} \) is ergodic and for $\mu_{\alpha}$ almost every $x$, $r(x) = n$.

$(3)\implies(1):$ The argument is similar to the proof of $(1)\implies (2)$. We still prove by contradiction. Assume $\mu_{\alpha}$ is ergodic and $R_n$ is of full measure, and $\delta_{\rho}(\hat{\omega}_{\alpha}) = \delta$ for some $\delta<1$. Fix $b_0\in \mathbb{D}$,  pick \( R>0 \) adequate with respect to $(b_0,\Gamma, \mu_{\alpha})$, so for \( \mu_{\alpha} \)-almost all \( x \), both \( r(x) = n \) and \( x\in \mathcal{O}_R(b_0,\gamma(b_0)) \) for infinitely many \( \gamma\in \Gamma \) hold.

Fix such an \( x \), choose $y,I,g,\tilde{f}$ as we did in Lemma~\ref{leminacute}, and from the same lemma, there exists a constant $C>1$ such that for any \( \gamma \) with \( \mathcal{O}_R(b_0,\gamma(b_0))\subset I \), we have
\[
\frac{1}{C} \,\|\tilde{f}(q_{\gamma})-\tilde{f}(p_{\gamma}) \|_{2}
\;\le\;
\mu_{\alpha}\bigl(\mathcal{O}_{R}(b_0,\gamma(b_0))\bigr)
\;\le\;
C \,\|\tilde{f}(q_{\gamma})-\tilde{f}(p_{\gamma}) \|_{2}.
\]
\[
\frac{1}{C}\,
\bigl(\det (\tilde{f}(q_{\gamma})-\tilde{f}(p_{\gamma}))\bigr)^{\frac{\delta}{n}}
\;\le\;
\mu_{\hat{\omega}_{\alpha}}\bigl(\mathcal{O}_{R}(b_0,\gamma(b_0))\bigr)
\;\le\;
C \,
\bigl(\det (\tilde{f}(q_{\gamma})-\tilde{f}(p_{\gamma}))\bigr)^{\frac{\delta}{n}}.
\]
For any diverging sequence \( \gamma_k\in \Gamma \) such that \( q_{\gamma_k},p_{\gamma_k}\to x \) from different sides, since \( r(x)=n \),
\[
\frac{\tilde{f}(q_{\gamma_k})-\tilde{f}(p_{\gamma_k})}{\|\tilde{f}(q_{\gamma_k})-\tilde{f}(p_{\gamma_k}) \|_{2}}
\text{ converges to a positive definite matrix},
\]
so
\[
\frac{\det (\tilde{f}(q_{\gamma_k})-\tilde{f}(p_{\gamma_k}))}{\|\tilde{f}(q_{\gamma_k})-\tilde{f}(p_{\gamma_k}) \|_{2}^{\,n}}
\; \longrightarrow c
\]
where $c$ is some positive number. So for $k\gg 1$, $\frac{\det (\tilde{f}(q_{\gamma_k})-\tilde{f}(p_{\gamma_k}))}{\|\tilde{f}(q_{\gamma_k})-\tilde{f}(p_{\gamma_k}) \|_{2}^{\,n}}\ge \frac{c}{2}$, so we have
\[
\begin{aligned}
\frac{
    \mu_{\alpha}(\mathcal{O}_{R}(b_0,\gamma_k(b_0)))
}{
    \mu_{\hat{\omega}_{\alpha}}(\mathcal{O}_{R}(b_0,\gamma_k(b_0)))
}
&\le
C^2 \,
\frac{
    \|\tilde{f}(q_{\gamma_k})-\tilde{f}(p_{\gamma_k}) \|_{2}
}{
    \bigl(\det (\tilde{f}(q_{\gamma_k})-\tilde{f}(p_{\gamma_k}))\bigr)^{\frac{\delta}{n}}
}
\\[0.5em]
&=
C^2 \,
\frac{
    \|\tilde{f}(q_{\gamma_k})-\tilde{f}(p_{\gamma_k}) \|_{2}
}{
    \det (\tilde{f}(q_{\gamma_k})-\tilde{f}(p_{\gamma_k}))^{\frac{1}{n}}
}
\cdot
\bigl(\det (\tilde{f}(q_{\gamma_k})-\tilde{f}(p_{\gamma_k}))\bigr)^{\frac{1-\delta}{n}}
\\[0.5em]
&\le
\left(\frac{2}{c}\right)^{\frac{1}{n}}
C^2 \,
\bigl(\det (\tilde{f}(q_{\gamma_k})-\tilde{f}(p_{\gamma_k}))\bigr)^{\frac{1-\delta}{n}}.
\end{aligned}
\]
Since $\delta<1$, the right-hand side converges to $0$ as $k\to \infty$. As the discussion holds for any choice of $x$, $I$, and $\{\gamma_k\}$, applying Lemma~\ref{vitaliargument}, we see that $\mu_{\alpha}=0$, which is a contradiction.
\end{proof}
\begin{corollary}\label{productergodicity}
If \( \delta_{\rho}(\hat{\omega}_{\alpha}) = 1 \), then \( \mu_{\alpha} \times \mu_{\alpha} \) is ergodic with respect to the diagonal \( \Gamma \)-action. Consequently, Corollary \ref{transverse} implies that the set \( \mathcal{ET} \) has full \( \mu_{\alpha} \times \mu_{\alpha} \)-measure.
\end{corollary}

\begin{proof}
By Theorem~\ref{ergodicities} and Proposition~\ref{lebesgueergodicity}~$(3)$,  the quasi-invariant measures \( \mu_{\alpha} \) and \( \mu_{\hat{\omega}_\alpha} \) are ergodic with respect to the \( \Gamma \)-action on \( \Lambda(\Gamma) \). Ergodic measures are either mutually singular or mutually absolutely continuous. Lemma \ref{domination} implies that \(\mu_{\hat{\omega}_{\alpha}}\) is not singular to \(\mu_{\alpha}\), so \( \mu_{\alpha} \) and \( \mu_{\hat{\omega}_\alpha} \) are mutually absolutely continuous.

Hence, the product measures \( \mu_{\alpha} \times \mu_{\alpha} \) and \( \mu_{\hat{\omega}_\alpha} \times \mu_{\hat{\omega}_\alpha} \) are also mutually absolutely continuous. Since the latter is ergodic under the diagonal \( \Gamma \)-action (Theorem \ref{ergodicities}), so is the former.
\end{proof}

\begin{proof}[Proof of Theorem~\ref{main}]

Throughout, we assume $n>2$, $\Gamma\subset \mathrm{PSL}(2,\mathbb{R})$ is a lattice, and $\rho:\Gamma\to\mathrm{Sp}(2n,\mathbb{R})$ is maximal. The case $n = 2$ can be covered by considering $\rho\oplus \rho:\Gamma\to \mathrm{Sp}(8,\mathbb{R})$ which is still maximal.

We set $\rho' := \Lambda^n\circ \rho$. As before, we write $V=\bigwedge^n\mathbb{R}^{2n}$, $d=\dim V$, $I_1=\{1,d-1\}$, and $I_2=\{1,2,d-2,d-1\}$. Let $\xi^1:\Lambda(\Gamma)\to\mathcal{F}_1(V)$ and $\xi^{I_1}:\Lambda(\Gamma)\to\mathcal{F}_{I_1}(V)$ be the $1$-limit and $I_1$-limit maps of $\rho'$, respectively. Let $\xi^{I_2}: \Lambda(\Gamma)\to\mathcal{F}_{I_2}(V)$ be the equivariant measurable map as described in Remark~\ref{ngreaterthan2}.

The main idea of the proof is the following lemma about the quantitative estimate on the Radon--Nikodym derivative $\frac{d\mu_{\hat{\omega}_{\alpha}}}{d\mu_{\alpha}}$. (Note that Lemma~\ref{domination} only implies that it is bounded above.) 

\begin{lemma}\label{radonnikodymupperanddown}
If $\delta_{\rho}(\hat{\omega}_{\alpha}) = 1$, then \(\frac{d\mu_{\hat{\omega}_{\alpha}}}{d\mu_{\alpha}}\) is \red{locally} bounded above and away from $0$ on $\Lambda(\Gamma)$.  
\end{lemma}
\begin{proof}
 Recall that \( m^*_{\alpha} \) and $m^*_{\hat{\omega}_{\alpha}}$ are the $\Gamma\times \mathbb{R}$ invariant measures defined in Section~\ref{bmsdef}. If $\delta_{\rho}(\hat{\omega}_{\alpha}) = 1$, Corollary~\ref{productergodicity} implies that \( m^*_{\alpha} \) is ergodic with respect to the $\Gamma \times \mathbb{R}$ action. Moreover, Theorem~\ref{ergodicities} shows that $m^*_{\hat{\omega}_{\alpha}}$ is also ergodic. Lemma~\ref{domination} shows that \( m^*_{\alpha} \) and \( m^*_{\hat{\omega}_{\alpha}} \) are not mutually singular, thus they are proportional, i.e.,  there exists a constant \( c > 0 \) such that
\[
m^*_{\alpha} = c \, m^*_{\hat{\omega}_{\alpha}}.
\]
Equivalently, for $\mu_{\alpha}\times \mu_{\alpha}$-almost all $(x,y) \in \Lambda(\Gamma)\times \Lambda(\Gamma)$,
\[
e^{-\mathcal{J}_2([\xi^{I_2}(x),\xi^{I_2}(y)]_{I_2})}
=
c \, e^{-\hat{\omega}_1([\xi^{I_1}(x),\xi^{I_1}(y)]_{I_1})}
\cdot
\frac{d\mu_{\hat{\omega}_{\alpha}}}{d\mu_{\alpha}}(x)
\cdot
\frac{d\mu_{\hat{\omega}_{\alpha}}}{d\mu_{\alpha}}(y),
\]
where $\frac{d\mu_{\hat{\omega}_{\alpha}}}{d\mu_{\alpha}}$ is the Radon--Nikodym derivative of these mutually absolutely continuous measures.

Since $\mathcal{J}_2 = 2\omega_1 - \omega_2$ and $\hat{\omega}_1 = \frac{2}{n} \omega_1$, this can be rewritten as
\[
e^{\omega_2([\xi^{I_2}(x),\xi^{I_2}(y)]_{I_2})}
=
c \,
e^{2\left(1 - \frac{1}{n}\right) \omega_1([\xi^{I_1}(x),\xi^{I_1}(y)]_{I_1})}
\cdot
\frac{d\mu_{\hat{\omega}_{\alpha}}}{d\mu_{\alpha}}(x)
\cdot
\frac{d\mu_{\hat{\omega}_{\alpha}}}{d\mu_{\alpha}}(y).
\]

Let \( v_1 = \mathrm{Span}(e_1 \wedge \cdots \wedge e_n) \) and \( v_2 = \mathrm{Span}(e_{n+1} \wedge \cdots \wedge e_{2n}) \). For any \( M \in \overline{\mathrm{Pos}}(n) \), recall the matrices
\[
U^M = \begin{pmatrix}
\mathrm{Id}_n & M \\
0 & \mathrm{Id}_n
\end{pmatrix},
\qquad
U_M = \begin{pmatrix}
\mathrm{Id}_n & 0 \\
M & \mathrm{Id}_n
\end{pmatrix},
\]
and
\[
u^M =
\begin{pmatrix}
0 & M \\
0 & 0
\end{pmatrix},
\qquad
u_M =
\begin{pmatrix}
0 & 0 \\
M & 0
\end{pmatrix}.
\]

For any \((x,y) \in \mathcal{ET}\) (recall that $\mathcal{ET}$ is defined as the set such that $\ell(x)$ and $\mathrm{Ann}(\ell(y))$ exist and are transverse, see the discussion preceding Corollary~\ref{transverse}), choose two disjoint intervals \(J_x\) and \(J_y\) in $\Lambda(\Gamma)$ such that \(J_x\) terminates at \(x\), \(J_y\) terminates at \(y\), and both lie on the clockwise arc from \(x\) to \(y\). Let
\[
L_x = \mu_{\alpha}(J_x), \qquad L_y = \mu_{\alpha}(J_y).
\]
By Lemma~\ref{lemma:parameterizationofmaps_combined}, Remark~\ref{remark:variationofparameterization}, and the arc-length parametrization, there exist:
\begin{itemize}
    \item an element \( g \in \mathrm{Aut}(\mathfrak{sp}(2n,\mathbb{R})) \),
    \item parametrizations \( p_x : [0, L_x] \to J_x \) and \( p_y : [0, L_y] \to J_y \) satisfying
    \[
    \mu_{\alpha}(p_x([0,s])) = s, \qquad
    \mu_{\alpha}(p_y([0,t])) = t,
    \]
    \item monotone and almost everywhere differentiable maps
    \[
    f_x : [0, L_x] \to \mathrm{Pos}(n) \cup \{0\}, \qquad
    f_y : [0, L_y] \to \mathrm{Pos}(n) \cup \{0\},
    \]
\end{itemize}
such that for all \( s \in [0,L_x] \), \( t \in [0,L_y] \),
\[
(\xi^{1}(x),\, \xi^{1}(p_x(s)),\, \xi^{1}(y))
= g\bigl(v_1,\, U_{f_x(s)} v_1,\, v_2 \bigr),\quad 
(\xi^{1}(x),\, \xi^{1}(p_y(t)),\, \xi^{1}(y))
= g\bigl(v_1,\, U^{f_y(t)} v_2,\, v_2 \bigr).
\]

Let \( f_x' : [0,L_x] \to \overline{\mathrm{Pos}}(n) \) and \( f_y':[0,L_y] \to \overline{\mathrm{Pos}}(n) \) denote the derivatives. Since \(p_x,p_y\) are chosen as arc-length parametrizations with respect to \(\mu_\alpha\), and \(g\) is fixed, the matrix norms of \(f'_x\) and \(f'_y\) are uniformly bounded above and away from $0$. By Corollary~\ref{productergodicity}, \(\mathcal{ET}\) has full \(\mu_{\alpha}\times\mu_{\alpha}\)-measure, so we may assume that for almost every \(s,t\), the pair \((\xi^{I_2}(p_x(s)), \xi^{I_2}(p_y(t)))\) is defined and transverse.

By Lemma~\ref{lemma:computations}, it is straightforward to check that for any $0 < s < L_x$, $0 < t < L_y$, and any $\epsilon_1<L_x-s, \epsilon_2<L_y-t$, we have
\begin{align*}
g\bigl(U_{f_x(s+\epsilon_1)} v_1,\; U^{f_y(t+\epsilon_2)} v_2\bigr)
&=
g\bigl(U_{f_x(s+\epsilon_1)} v_1,\; U_{f_y(t+\epsilon_2)^{-1}} v_1\bigr) \\
&=
g\, U_{f_x(s)}\bigl(U_{f_x(s+\epsilon_1)-f_x(s)} v_1,\; U_{\,f_y(t+\epsilon_2)^{-1}-f_x(s)} v_1\bigr) \\
&=
g\, U_{f_x(s)}\bigl(U_{f_x(s+\epsilon_1)-f_x(s)} v_1,\; U^{\,(f_y(t+\epsilon_2)^{-1}-f_x(s))^{-1}} v_2\bigr).
\end{align*}
Here note that as $(x,p_x(s),p_y(t+\epsilon_2),y)$ is positive, we have $(v_1, U_{f_x(s)} v_1,\; U^{f_y(t+\epsilon_2)} v_2, v_2\bigr) = (v_1, U_{f_x(s)} v_1,\; U_{f_y(t+\epsilon_2)^{-1}} v_1, v_2\bigr)$ is positive, so $f_y(t+\epsilon_2)^{-1}-f_x(s)$ is positive definite, hence invertible.

Hence, for almost all $s,t$, let $\epsilon_2\to 0$ to compute the derivative, and by a calculation similar to Proposition~\ref{tangent},
\begin{align*}
&(\ell(p_x(s)),\ell(p_y(t)))\\
=&g\, U_{f_x(s)} \Bigl(
    \mathrm{Span}\bigl(
        v_1,\;
        \sum_{i=1}^{n}
        e_1 \wedge \cdots \wedge e_{i-1}
        \wedge u_{f'_x(s)} e_i
        \wedge e_{i+1} \wedge \cdots \wedge e_n
    \bigr), \\
&
    \red{U^{(f_y(t)^{-1}-f_x(s))^{-1}}}\mathrm{Span}\bigl(
        v_2,\;
        \sum_{i=1}^{n}
        e_{n+1} \wedge \cdots \wedge e_{n+i-1}
        \wedge u^{\tilde{f}_y(s,t)} e_{n+i}
        \wedge e_{n+i+1} \wedge \cdots \wedge e_{2n}
    \bigr)
\Bigr)\\
=& g \, U_{f_x(s)}\red{U^{(f_y(t)^{-1}-f_x(s))^{-1}}} \Bigl(
    \mathrm{Span}\bigl(
        v_1,\;
        \sum_{i=1}^{n}
        e_1 \wedge \cdots \wedge e_{i-1}
        \wedge \red{U^{-(f_y(t)^{-1}-f_x(s))^{-1}}} u_{f'_x(s)} e_i
        \wedge e_{i+1} \wedge \cdots \wedge e_n
    \bigr), \\
&
    \mathrm{Span}\bigl(
        v_2,\;
        \sum_{i=1}^{n}
        e_{n+1} \wedge \cdots \wedge e_{n+i-1}
        \wedge u^{\tilde{f}_y(s,t)} e_{n+i}
        \wedge e_{n+i+1} \wedge \cdots \wedge e_{2n}
    \bigr)
\Bigr)\\
=& g \, U_{f_x(s)}\red{U^{(f_y(t)^{-1}-f_x(s))^{-1}}} \Bigl(
    \mathrm{Span}\bigl(
        v_1,\;
        \sum_{i=1}^{n}
        e_1 \wedge \cdots \wedge e_{i-1}
        \wedge  u_{f'_x(s)} e_i
        \wedge e_{i+1} \wedge \cdots \wedge e_n
    \bigr), \\
&
    \mathrm{Span}\bigl(
        v_2,\;
        \sum_{i=1}^{n}
        e_{n+1} \wedge \cdots \wedge e_{n+i-1}
        \wedge u^{\tilde{f}_y(s,t)} e_{n+i}
        \wedge e_{n+i+1} \wedge \cdots \wedge e_{2n}
    \bigr)
\Bigr)
,
\end{align*}
where
\[
\tilde{f}_y(s,t)
\red{=\frac{\mathrm{d}}{\mathrm{d}t}\left( f_y(t)^{-1}-f_x(s) \right)^{-1}} 
=\bigl(\mathrm{Id}_n-f_y(t)f_x(s)\bigr)^{-1}
\,\, f'_y(t)\, \,
\bigl(\mathrm{Id}_n-f_x(s)f_y(t)\bigr)^{-1}.
\]
Here note again we have shown that $f_y(t+\epsilon_2)^{-1}-f_x(s)$ is invertible. By taking $\epsilon_2 = 0$, we have $f_y(t)^{-1}-f_x(s)$ is invertible, so $\mathrm{Id}_n-f_y(t)f_x(s)$ and its transpose $\mathrm{Id}_n-f_x(s)f_y(t)$ are both invertible.

\begin{redblock}
Since $g \, U_{f_x(s)}\red{U^{(f_y(t)^{-1}-f_x(s))^{-1}}}$ is bounded, the cocycle identity of the Gromov product and Lemma~\ref{Gromovproduct} imply that there exists a constant \( C > 1 \) (depending on \(g,x,y,J_x,J_y\)) such that
\[
\frac{1}{C} \, \mathrm{tr}\bigl(f_x'(s) \tilde{f}_y(s,t)\bigr)
\;\le\;
e^{\omega_2([\xi^{I_2}(p_x(s)), \xi^{I_2}(p_y(t))]_{I_2})}
\;\le\;
C \, \mathrm{tr}\bigl(f_x'(s) \tilde{f}_y(s,t)\bigr).
\]   
\end{redblock}

Since \(\xi^{I_1}\) is continuous and transverse, the function $J_x\times J_y\to \mathbb{R}$ given by
\[
(u,v)\mapsto e^{2\left(1 - \frac{1}{n}\right) \omega_1([\xi^{I_1}(u),\xi^{I_1}(v)]_{I_1})}
\]
is bounded. Therefore, combining the previous estimate with
\[
e^{\omega_2([\xi^{I_2}(x),\xi^{I_2}(y)]_{I_2})}
=
c \,
e^{2\left(1 - \frac{1}{n}\right) \omega_1([\xi^{I_1}(x),\xi^{I_1}(y)]_{I_1})}
\cdot
\frac{d\mu_{\hat{\omega}_{\alpha}}}{d\mu_{\alpha}}(x)
\cdot
\frac{d\mu_{\hat{\omega}_{\alpha}}}{d\mu_{\alpha}}(y),
\]
and by enlarging $C$ if necessary we have that for a.e.\ \((s,t)\in [0,L_x]\times [0,L_y]\),
\begin{redblock}
\[
\frac{1}{C} \, \mathrm{tr}\bigl(f_x'(s) \tilde{f}_y(s,t)\bigr)
\;\le\;
\frac{d\mu_{\hat{\omega}_{\alpha}}}{d\mu_{\alpha}}(p_x(s)) \cdot \frac{d\mu_{\hat{\omega}_{\alpha}}}{d\mu_{\alpha}}(p_y(t))
\;\le\;
C \, \mathrm{tr}\bigl(f_x'(s) \tilde{f}_y(s,t)\bigr).
\]

Since $\mu_{\alpha}(p_y([0,t])) = t$ for all $t \in [0,L_y]$ and $\tilde{f}_y(s,t)=\frac{\mathrm{d}}{\mathrm{d}t}\left( f_y(t)^{-1}-f_x(s) \right)^{-1} $ , integrating in \(t\) over \([0,L_y]\) yields, after further enlarging $C$ if necessary, we have for almost every $s\in [0,L_x]$, 
\[
\frac{1}{C} \, \mathrm{tr}\bigl(f_x'(s) \tilde{M}(s)\bigr)
\;\le\;
\frac{d\mu_{\hat{\omega}_{\alpha}}}{d\mu_{\alpha}}(p_x(s)) \cdot \mu_{\hat{\omega}_{\alpha}}(J_y)
\;\le\;
C \, \mathrm{tr}\bigl(f_x'(s) \tilde{M}(s)\bigr)
\]
where $\tilde{M}(s) = (f_y(L_y)^{-1} -f_x(s))^{-1}$ is positive definite, and bounded from above and away from $0$. Pick $A(s)\in \mathrm{Pos}(n)$ so that $\tilde{M}(s) = A(s)^2$, so $\mathrm{tr}\bigl(f_x'(s) \tilde{M}(s)\bigr) = \mathrm{tr}\bigl(f_x'(s) A(s)^2\bigr) = \mathrm{tr}\bigl(A(s)f_x'(s) A(s)\bigr)$, so after enlarging $C$ if necessary, we have
\[
\frac{1}{\tilde{C}} \, \|f_x'(s)\|_{2}
\;\le\;
\mathrm{tr}\bigl(f_x'(s)\tilde{M}(s)\bigr)
\;\le\;
\tilde{C} \, \|f_x'(s)\|_{2}.
\]
\end{redblock}
Thus, enlarge $C$ again if necessary, we have for a.e.\ \(s\in [0,L_x]\),
\[
\frac{1}{C} \, \|f_x'(s)\|_{2}
\;\le\;
\frac{d\mu_{\hat{\omega}_{\alpha}}}{d\mu_{\alpha}}(p_x(s))
\;\le\;
C \, \|f_x'(s)\|_{2}.
\]

Since $p_x,p_y$ parameterize the limit curve at unit speed, the matrix norms of \( f'_x \) and \( f'_y \) are uniformly bounded above and away from $0$, and this shows that the Radon--Nikodym derivative
\(
\frac{d\mu_{\hat{\omega}_{\alpha}}}{d\mu_{\alpha}}(p_x(s))
\)
is uniformly bounded above and away from $0$ for almost every \( s \in [0,L_x] \). Since $\mu_\alpha$ has full support and $\mathcal{ET}$ has full
$\mu_\alpha\times\mu_\alpha$-measure, for every $z\in\Lambda(\Gamma)$
we can choose $(x,y)\in\mathcal{ET}$ such that $z$ lies in the open
clockwise arc from $x$ to $y$. We may then choose $J_x$ and $J_y$
as above with $z\in\operatorname{int}(J_x)$. The interiors of these
intervals $J_x$ cover $\Lambda(\Gamma)$. Taking a finite subcover
gives uniform upper and positive lower bounds for
$\frac{d\mu_{\hat{\omega}_\alpha}}{d\mu_\alpha}$ almost everywhere.
\end{proof}
 
Now we return to the proof of the main theorem. Fix a sufficiently large $R>0$. For each infinite order $\gamma \in \Gamma$, the ratio of shadows
\[
\frac{\mu_{\hat{\omega}_{\alpha}}\bigl(\mathcal{O}_{R}(b_0,\gamma^k(b_0))\bigr)}{\mu_{\alpha}\bigl(\mathcal{O}_{R}(b_0,\gamma^k(b_0))\bigr)},
\qquad k \in \mathbb{Z}^{+},
\]
is uniformly bounded above and away from $0$, independently of $\gamma$ and $k$. Thus, by Theorem~\ref{shadowlemma}, we have
\[
\frac{e^{-\alpha(\kappa(\rho(\gamma^k)))}}{e^{-\hat{\omega}_\alpha(\kappa(\rho(\gamma^k)))}}
=
e^{-(\alpha - \hat{\omega}_\alpha)(\kappa(\rho(\gamma^k)))},
\]
is also bounded above and away from $0$ independently of $k$. This implies that, for every $\gamma \in \Gamma$, the first $n$ eigenvalues of $\rho(\gamma)$ have the same modulus (so does the last $n$ eigenvalues).

\begin{redblock}

Let $G_{\rho}\subset \mathrm{Sp}(2n,\mathbb{R})$ be the Zariski closure of $\rho(\Gamma)$, which is reductive with compact centralizer (see Burger--Iozzi--Wienhard~\cite[Corollary~4]{burger2008tighthomomorphismshermitiansymmetric}). Let $\mathfrak{g}_{\rho}\subset \mathfrak{sp}(2n,\mathbb{R})$ be its Lie algebra. Write
\[
\mathfrak{g}_{\rho} = \mathfrak{z}_{\rho}\oplus \mathfrak{g}_{\rho,c}\oplus \mathfrak{g}_{\rho,nc},
\]
where $\mathfrak{z}_{\rho}$ is its center, and $\mathfrak{g}_{\rho,c}$ (resp.\ $\mathfrak{g}_{\rho,nc}$) is the compact (resp.\ non-compact) factor of its semisimple part. Since $\mathrm{Ad}_{G_{\rho}}$ permutes the simple factors of the semisimple part, $\mathfrak{g}_{\rho,nc}$ is $\mathrm{Ad}_{G_{\rho}}$-invariant. Let $G_{\rho,nc}\subset \mathrm{Sp}(2n,\mathbb{R})$ be the connected closed subgroup with Lie algebra $\mathfrak{g}_{\rho, nc}$. Burger--Iozzi--Wienhard~\cite[Theorem~5]{burger2008surfacegrouprepresentationsmaximal} show that the symmetric space associated with $G_{\rho,nc}$ is Hermitian of tube type.

Fix a Cartan decomposition of $\mathfrak{g}_{\rho, nc}$ and extend it to $\mathfrak{g}_{\rho}$. Up to $\mathrm{Ad}_{\mathrm{Sp}(2n,\mathbb{R})}$ conjugation, we may assume that the split real Cartan subalgebra of $\mathfrak{g}_{\rho,nc}$ belongs to $\mathfrak{a}$, and that the maximal compact subgroup of $G_{\rho}$ is contained in $\mathrm{U}(n)$.

Let $G$ be a reductive real Lie group with compact center, equipped with a chosen Cartan decomposition, and let $\Gamma'\subset G$ be a subgroup. Recall that the Benoist limit cone $\mathcal{B}(G;\Gamma')$ (see Benoist~\cite{Benoistasym}) is defined as the minimal closed cone in the split real Cartan subalgebra containing all rays along the Jordan projections of $\Gamma'$. Benoist~\cite{Benoistasym} proved that if $\Gamma'$ is Zariski dense in $G$, then $\mathcal{B}(G;\Gamma')$ has a non-empty interior.

The preceding arguments show that $\mathcal{B}(\mathrm{Sp}(2n,\mathbb{R}); \rho(\Gamma))$ is the ray $\mathbb{R}^{\ge0} \mathrm{diag}(\mathrm{Id}_n,-\mathrm{Id}_n)$. On the other hand, $\mathcal{B}(\mathrm{Sp}(2n,\mathbb{R}); \rho(\Gamma))$ equals the intersection between $\overline{\mathfrak{a}^+}$ and the $\mathrm{Sp}(2n,\mathbb R)$-Weyl group translates of $\mathcal{B}(G_{\rho};\rho(\Gamma))$. Since $\rho(\Gamma)\subset G_{\rho}$ is Zariski dense, $G_{\rho}$ has split real rank one. Hence, $\mathfrak{g}_{\rho, nc}$ also has split real rank one. Because $G_{\rho, nc}$ is of Hermitian tube type, we have $\mathfrak{g}_{\rho, nc} \cong \mathfrak{sl}(2,\mathbb{R})$.

Consequently, the line $\mathbb{R} \mathrm{diag}(\mathrm{Id}_n,-\mathrm{Id}_n)$ lies in the Weyl group orbit of the Cartan subalgebra of $\mathfrak{g}_{\rho, nc}$. Up to conjugation, we can thus assume that $\mathfrak{g}_{\rho, nc}$ is generated by the $\mathfrak{sl}_2$ triple
\[
H = \mathrm{diag}(\mathrm{Id}_n,-\mathrm{Id}_n),\quad  E = \begin{pmatrix}0&A\\0&0 \end{pmatrix},\quad F =  \begin{pmatrix}0&0\\A^{-1}&0 \end{pmatrix}
\]
where $A$ is symmetric and invertible. By further taking conjugation and negating $E$ and $F$ if necessary, we can assume $A = A^{-1} =  I_{p,n-p}$ (the diagonal matrix where the first $p$ entries are $1$, and remaining entries are $-1$) where $p\ge n-p$. The only possibility for $\rho(\Gamma)$ to preserve a positive curve in $\mathcal{F}_{\alpha}$ is $p = n$. Indeed, note that the
attracting Lagrangians of the hyperbolic elements of \(G_{\rho,\mathrm{nc}}\)
form the circle
\[
\mathcal C_p
=\{\exp(tF)\mathcal L:t\in\mathbb R\}
 \cup\{\mathcal L^{\mathrm{opp}}\},
\qquad
\exp(tF)=
\begin{pmatrix}
\mathrm{Id}_n & 0\\
tI_{p,n-p} & \mathrm{Id}_n
\end{pmatrix}.
\]
For a hyperbolic \(\gamma\in\rho^{-1}(G_\rho^\circ)\), applying Equation~(\ref{anosovgrowth}) to its powers shows that the noncompact component of \(\rho(\gamma)\) is
hyperbolic. Since the commuting compact component preserves its attracting
Lagrangian, we have \(\xi^\alpha(\gamma^+)\in\mathcal C_{p}\). And note that such $\gamma^+$ forms a dense set in  \(\Lambda(\Gamma) = S^1\), thus the limit
curve equals \(\mathcal C_p\). In particular, the triple
\((\mathcal L,\operatorname{graph}(I_{p,n-p}),\mathcal L^{\mathrm{opp}})\)
is positive, forcing \(I_{p,n-p}\)
to be definite, so $p = n$. As a result, $\mathfrak{g}_{\rho, nc} = \mathrm{d}\rho_D (\mathfrak{sl}(2,\mathbb R))$, so $G_{\rho, nc} = \rho_D (\mathrm{SL}(2,\mathbb R))$. Since $G_{\rho}$ normalizes $G_{\rho, nc}$, $G_{\rho}\subset N_D$.

\end{redblock}
\end{proof}

\bibliographystyle{amsplain}
\bibliography{ref}

@incollection {schapira:hal-01882452,
    AUTHOR = {Schapira, Barbara},
     TITLE = {Dynamics of geodesic and horocyclic Flows},
 BOOKTITLE = {Ergodic Theory and Negative Curvature},
    SERIES = {Lecture Notes in Math.},
    VOLUME = {2164},
     PAGES = {129--155},
 PUBLISHER = {Springer, Cham},
      YEAR = {2017},
      ISBN = {978-3-319-43058-4; 978-3-319-43059-1; 978-2-85629-860-2},
   MRCLASS = {37D40 (53D25)},
  MRNUMBER = {3588134},
MRREVIEWER = {Rafael\ Oswaldo\ Ruggiero},
}

@article {sullivanacta,
    AUTHOR = {Sullivan, Dennis},
     TITLE = {Entropy, {H}ausdorff measures old and new, and limit sets of
              geometrically finite {K}leinian groups},
   JOURNAL = {Acta Math.},
  FJOURNAL = {Acta Mathematica},
    VOLUME = {153},
      YEAR = {1984},
    NUMBER = {3-4},
     PAGES = {259--277},
      ISSN = {0001-5962,1871-2509},
   MRCLASS = {58F11 (28D99 30F40)},
  MRNUMBER = {766265},
MRREVIEWER = {Caroline\ Series},
       DOI = {10.1007/BF02392379},
       URL = {https://doi-org.libproxy1.nus.edu.sg/10.1007/BF02392379},
}

@article {beyrerpozzettisopq,
    AUTHOR = {Beyrer, Jonas and Pozzetti, Maria Beatrice},
     TITLE = {Positive surface group representations in {$\sf{PO}(p,q)$}},
   JOURNAL = {J. Eur. Math. Soc. (JEMS)},
  FJOURNAL = {Journal of the European Mathematical Society (JEMS)},
    VOLUME = {28},
      YEAR = {2026},
    NUMBER = {6},
     PAGES = {2307--2349},
      ISSN = {1435-9855,1435-9863},
   MRCLASS = {22E40 (32G15)},
  MRNUMBER = {5045688},
MRREVIEWER = {Marja\ K.\ Kankaanrinta},
       DOI = {10.4171/jems/1560},
       URL = {https://doi-org.libproxy1.nus.edu.sg/10.4171/jems/1560},
}

@book {imayoshi,
    AUTHOR = {Imayoshi, Y. and Taniguchi, M.},
     TITLE = {An introduction to {T}eichm\"uller spaces},
      NOTE = {Translated and revised from the Japanese by the authors},
 PUBLISHER = {Springer-Verlag, Tokyo},
      YEAR = {1992},
     PAGES = {xiv+279},
      ISBN = {4-431-70088-9},
   MRCLASS = {32G15 (14H15 30F10 30F60)},
  MRNUMBER = {1215481},
MRREVIEWER = {I.\ Kra},
       DOI = {10.1007/978-4-431-68174-8},
       URL = {https://doi-org.libproxy1.nus.edu.sg/10.1007/978-4-431-68174-8},
}

@article {CaZhZi-relatively,
    AUTHOR = {Canary, Richard and Zhang, Tengren and Zimmer, Andrew},
     TITLE = {Patterson-{S}ullivan measures for relatively {A}nosov groups},
   JOURNAL = {Math. Ann.},
  FJOURNAL = {Mathematische Annalen},
    VOLUME = {392},
      YEAR = {2025},
    NUMBER = {2},
     PAGES = {2309--2363},
      ISSN = {0025-5831,1432-1807},
   MRCLASS = {22F10 (22E40)},
  MRNUMBER = {4906322},
MRREVIEWER = {Boris\ Hasselblatt},
       DOI = {10.1007/s00208-025-03137-2},
       URL = {https://doi-org.libproxy1.nus.edu.sg/10.1007/s00208-025-03137-2},
}

@article {maximalearlier,
    AUTHOR = {Burger, Marc and Iozzi, Alessandra and Wienhard, Anna},
     TITLE = {Surface group representations with maximal {T}oledo invariant},
   JOURNAL = {C. R. Math. Acad. Sci. Paris},
  FJOURNAL = {Comptes Rendus Math\'ematique. Acad\'emie des Sciences. Paris},
    VOLUME = {336},
      YEAR = {2003},
    NUMBER = {5},
     PAGES = {387--390},
      ISSN = {1631-073X,1778-3569},
   MRCLASS = {53C35 (32A15 53C24 53C55)},
  MRNUMBER = {1979350},
MRREVIEWER = {David\ Michael\ Fisher},
       DOI = {10.1016/S1631-073X(03)00065-7},
       URL = {https://doi-org.libproxy1.nus.edu.sg/10.1016/S1631-073X(03)00065-7},
}

@misc{yzf,
      title={Critical Exponent Rigidity for $\Theta-$positive Representations}, 
      author={Zhufeng Yao},
      year={2026},
      eprint={2505.17559},
      archivePrefix={arXiv},
      primaryClass={math.DG},
      url={https://arxiv.org/abs/2505.17559}, 
}

@article {burger2008tighthomomorphismshermitiansymmetric,
    AUTHOR = {Burger, Marc and Iozzi, Alessandra and Wienhard, Anna},
     TITLE = {Tight homomorphisms and {H}ermitian symmetric spaces},
   JOURNAL = {Geom. Funct. Anal.},
  FJOURNAL = {Geometric and Functional Analysis},
    VOLUME = {19},
      YEAR = {2009},
    NUMBER = {3},
     PAGES = {678--721},
      ISSN = {1016-443X,1420-8970},
   MRCLASS = {53C24 (32M15 53C35)},
  MRNUMBER = {2563767},
MRREVIEWER = {Michelle\ Bucher-Karlsson},
       DOI = {10.1007/s00039-009-0020-8},
       URL = {https://doi.org/10.1007/s00039-009-0020-8},
}

@article {burgermanhattan,
    AUTHOR = {Burger, Marc},
     TITLE = {Intersection, the {M}anhattan curve, and
              {P}atterson-{S}ullivan theory in rank {$2$}},
   JOURNAL = {Int. Math. Res. Not.},
  FJOURNAL = {International Mathematics Research Notices},
      YEAR = {1993},
    NUMBER = {7},
     PAGES = {217--225},
      ISSN = {1073-7928,1687-0247},
   MRCLASS = {58F17 (53C20 53C35)},
  MRNUMBER = {1230298},
MRREVIEWER = {Boris\ Hasselblatt},
       DOI = {10.1155/S1073792893000236},
       URL = {https://doi.org/10.1155/S1073792893000236},
}

@misc{kim2026measurableboundarymapspattersonsullivan,
      title={Measurable boundary maps and Patterson--Sullivan measures for non-Borel Anosov groups on the Furstenberg boundary}, 
      author={Dongryul M. Kim and Andrew Zimmer},
      year={2026},
      eprint={2603.28694},
      archivePrefix={arXiv},
      primaryClass={math.GT},
      url={https://arxiv.org/abs/2603.28694}, 
}

@book {margulisthesis,
    AUTHOR = {Margulis, Grigoriy A.},
     TITLE = {On some aspects of the theory of {A}nosov systems},
    SERIES = {Springer Monographs in Mathematics},
      NOTE = {With a survey by Richard Sharp: Periodic orbits of hyperbolic
              flows,
              Translated from the Russian by Valentina Vladimirovna
              Szulikowska},
 PUBLISHER = {Springer-Verlag, Berlin},
      YEAR = {2004},
     PAGES = {vi+139},
      ISBN = {3-540-40121-0},
   MRCLASS = {37D20 (37C27 37C30 37C40)},
  MRNUMBER = {2035655},
MRREVIEWER = {Boris\ Hasselblatt},
       DOI = {10.1007/978-3-662-09070-1},
       URL = {https://doi-org.libproxy1.nus.edu.sg/10.1007/978-3-662-09070-1},
}

@article {bishopsteger,
    AUTHOR = {Bishop, Christopher and Steger, Tim},
     TITLE = {Representation-theoretic rigidity in {${\rm PSL}(2,{\bf R})$}},
   JOURNAL = {Acta Math.},
  FJOURNAL = {Acta Mathematica},
    VOLUME = {170},
      YEAR = {1993},
    NUMBER = {1},
     PAGES = {121--149},
      ISSN = {0001-5962,1871-2509},
   MRCLASS = {22E40 (30F35)},
  MRNUMBER = {1208564},
MRREVIEWER = {Gopal\ Prasad},
       DOI = {10.1007/BF02392456},
       URL = {https://doi-org.libproxy1.nus.edu.sg/10.1007/BF02392456},
}

@inproceedings {wienhardtalk,
    AUTHOR = {Wienhard, Anna},
     TITLE = {An invitation to higher {T}eichm\"uller theory},
 BOOKTITLE = {Proceedings of the {I}nternational {C}ongress of
              {M}athematicians---{R}io de {J}aneiro 2018. {V}ol. {II}.
              {I}nvited lectures},
     PAGES = {1013--1039},
 PUBLISHER = {World Sci. Publ., Hackensack, NJ},
      YEAR = {2018},
      ISBN = {978-981-3272-91-0; 978-981-3272-87-3},
   MRCLASS = {22E40 (57M50)},
  MRNUMBER = {3966798},
MRREVIEWER = {Caglar\ Uyanik},
}

@article {sambarino2014orbitalcountingproblemhyperconvex,
    AUTHOR = {Sambarino, Andr\'es},
     TITLE = {The orbital counting problem for hyperconvex representations},
   JOURNAL = {Ann. Inst. Fourier (Grenoble)},
  FJOURNAL = {Universit\'e{} de Grenoble. Annales de l'Institut Fourier},
    VOLUME = {65},
      YEAR = {2015},
    NUMBER = {4},
     PAGES = {1755--1797},
      ISSN = {0373-0956,1777-5310},
   MRCLASS = {37A20 (22E30 22E40)},
  MRNUMBER = {3449196},
MRREVIEWER = {Alexandre\ Jos\'e\ Santana},
       DOI = {10.5802/aif.2973},
       URL = {https://doi.org/10.5802/aif.2973},
}

@article {CaZhZi,
    AUTHOR = {Canary, Richard and Zhang, Tengren and Zimmer, Andrew},
     TITLE = {Entropy rigidity for cusped {H}itchin representations},
   JOURNAL = {J. Topol.},
  FJOURNAL = {Journal of Topology},
    VOLUME = {19},
      YEAR = {2026},
    NUMBER = {1},
     PAGES = {Paper No. e70064, 71},
      ISSN = {1753-8416,1753-8424},
   MRCLASS = {22E40 (30F35 37B40 37F32 53C30)},
  MRNUMBER = {5018584},
MRREVIEWER = {Sadayoshi\ Kojima},
       DOI = {10.1112/topo.70064},
       URL = {https://doi-org.libproxy1.nus.edu.sg/10.1112/topo.70064},
}

@article {Dey_2022,
    AUTHOR = {Dey, Subhadip and Kapovich, Michael},
     TITLE = {Patterson-{S}ullivan theory for {A}nosov subgroups},
   JOURNAL = {Trans. Amer. Math. Soc.},
  FJOURNAL = {Transactions of the American Mathematical Society},
    VOLUME = {375},
      YEAR = {2022},
    NUMBER = {12},
     PAGES = {8687--8737},
      ISSN = {0002-9947,1088-6850},
   MRCLASS = {20F65 (22E40 53C35)},
  MRNUMBER = {4504651},
MRREVIEWER = {Joan\ Porti},
       DOI = {10.1090/tran/8713},
       URL = {https://doi-org.libproxy1.nus.edu.sg/10.1090/tran/8713},
}

@article {dey2024ahlforsregularitypattersonsullivanmeasures,
    AUTHOR = {Dey, Subhadip and Kim, Dongryul M. and Oh, Hee},
     TITLE = {Ahlfors regularity of {P}atterson-{S}ullivan measures of
              {A}nosov groups and applications},
   JOURNAL = {Compos. Math.},
  FJOURNAL = {Compositio Mathematica},
    VOLUME = {162},
      YEAR = {2026},
    NUMBER = {4},
     PAGES = {743--791},
      ISSN = {0010-437X,1570-5846},
   MRCLASS = {22E40 (28A78)},
  MRNUMBER = {5099348},
       DOI = {10.1017/S0010437X26103169},
       URL = {https://doi-org.libproxy1.nus.edu.sg/10.1017/S0010437X26103169},
}

@article {patterson1976limit,
    AUTHOR = {Patterson, S. J.},
     TITLE = {The limit set of a {F}uchsian group},
   JOURNAL = {Acta Math.},
  FJOURNAL = {Acta Mathematica},
    VOLUME = {136},
      YEAR = {1976},
    NUMBER = {3-4},
     PAGES = {241--273},
      ISSN = {0001-5962,1871-2509},
   MRCLASS = {30A58 (10D10 20H10)},
  MRNUMBER = {450547},
MRREVIEWER = {B.\ Maskit},
       DOI = {10.1007/BF02392046},
       URL = {https://doi-org.libproxy1.nus.edu.sg/10.1007/BF02392046},
}

@article {toledo89,
    AUTHOR = {Toledo, Domingo},
     TITLE = {Representations of surface groups in complex hyperbolic space},
   JOURNAL = {J. Differential Geom.},
  FJOURNAL = {Journal of Differential Geometry},
    VOLUME = {29},
      YEAR = {1989},
    NUMBER = {1},
     PAGES = {125--133},
      ISSN = {0022-040X,1945-743X},
   MRCLASS = {57N05 (22E40 32G13 57S25 58E20)},
  MRNUMBER = {978081},
MRREVIEWER = {Robert\ Brooks},
       URL = {http://projecteuclid.org.libproxy1.nus.edu.sg/euclid.jdg/1214442638},
}

@book {Hernandez,
    AUTHOR = {Hernandez Lamoneda, Luis},
     TITLE = {Maximal representations of surface groups in the group of
              isometries of a bounded symmetric domain},
      NOTE = {Thesis (Ph.D.)--The University of Utah},
 PUBLISHER = {ProQuest LLC, Ann Arbor, MI},
      YEAR = {1989},
     PAGES = {35},
   MRCLASS = {99-05},
  MRNUMBER = {2637776},
       URL =
              {http://gateway.proquest.com.libproxy1.nus.edu.sg/openurl?url_ver=Z39.88-2004&rft_val_fmt=info:ofi/fmt:kev:mtx:dissertation&res_dat=xri:pqdiss&rft_dat=xri:pqdiss:9004847},
}

@article {CaZhZi2,
    AUTHOR = {Canary, Richard and Zhang, Tengren and Zimmer, Andrew},
     TITLE = {Cusped {H}itchin representations and {A}nosov representations
              of geometrically finite {F}uchsian groups},
   JOURNAL = {Adv. Math.},
  FJOURNAL = {Advances in Mathematics},
    VOLUME = {404},
      YEAR = {2022},
     PAGES = {Paper No. 108439, 67},
      ISSN = {0001-8708,1090-2082},
   MRCLASS = {22E40 (20H10 30F35 37F32 57M60)},
  MRNUMBER = {4418884},
MRREVIEWER = {Herbert\ Abels},
       DOI = {10.1016/j.aim.2022.108439},
       URL = {https://doi-org.libproxy1.nus.edu.sg/10.1016/j.aim.2022.108439},
}

@article {CaZhZi3,
    AUTHOR = {Canary, Richard and Zhang, Tengren and Zimmer, Andrew},
     TITLE = {Patterson-{S}ullivan measures for transverse subgroups},
   JOURNAL = {J. Mod. Dyn.},
  FJOURNAL = {Journal of Modern Dynamics},
    VOLUME = {20},
      YEAR = {2024},
     PAGES = {319--377},
      ISSN = {1930-5311,1930-532X},
   MRCLASS = {22E40 (20H10 37A05)},
  MRNUMBER = {4799467},
MRREVIEWER = {Cheng\ Zheng},
       DOI = {10.3934/jmd.2024009},
       URL = {https://doi-org.libproxy1.nus.edu.sg/10.3934/jmd.2024009},
}

@article {GW2,
    AUTHOR = {Guichard, Olivier and Wienhard, Anna},
     TITLE = {Generalizing {L}usztig's total positivity},
   JOURNAL = {Invent. Math.},
  FJOURNAL = {Inventiones Mathematicae},
    VOLUME = {239},
      YEAR = {2025},
    NUMBER = {3},
     PAGES = {707--799},
      ISSN = {0020-9910,1432-1297},
   MRCLASS = {22E46 (22E15 22F30 22F50)},
  MRNUMBER = {4861126},
MRREVIEWER = {Alexandre\ Jos\'e\ Santana},
       DOI = {10.1007/s00222-024-01303-y},
       URL = {https://doi-org.libproxy1.nus.edu.sg/10.1007/s00222-024-01303-y},
}

@misc{beyrer2024positivitycrossratioscollarlemma,
 author = {Beyrer, Jonas and Guichard, Olivier and Labourie, Fran{\c{c}}ois and Pozzetti, Beatrice and Wienhard, Anna},
 title = {Positivity, cross-ratios and the collar lemma},
 note = {Preprint, {arXiv}:2409.06294}
 }

@article {fock2006modulispaceslocalsystems,
    AUTHOR = {Fock, Vladimir and Goncharov, Alexander},
     TITLE = {Moduli spaces of local systems and higher {T}eichm\"uller
              theory},
   JOURNAL = {Publ. Math. Inst. Hautes \'Etudes Sci.},
  FJOURNAL = {Publications Math\'ematiques. Institut de Hautes \'Etudes
              Scientifiques},
    NUMBER = {103},
      YEAR = {2006},
     PAGES = {1--211},
      ISSN = {0073-8301,1618-1913},
   MRCLASS = {32G15 (14D22 20F34 57M50)},
  MRNUMBER = {2233852},
MRREVIEWER = {William\ Goldman},
       DOI = {10.1007/s10240-006-0039-4},
       URL = {https://doi-org.libproxy1.nus.edu.sg/10.1007/s10240-006-0039-4},
}

@article {labourie2005anosovflowssurfacegroups,
    AUTHOR = {Labourie, Fran\c cois},
     TITLE = {Anosov flows, surface groups and curves in projective space},
   JOURNAL = {Invent. Math.},
  FJOURNAL = {Inventiones Mathematicae},
    VOLUME = {165},
      YEAR = {2006},
    NUMBER = {1},
     PAGES = {51--114},
      ISSN = {0020-9910,1432-1297},
   MRCLASS = {20F65 (37D20 37F30)},
  MRNUMBER = {2221137},
MRREVIEWER = {Richard\ Kenyon},
       DOI = {10.1007/s00222-005-0487-3},
       URL = {https://doi-org.libproxy1.nus.edu.sg/10.1007/s00222-005-0487-3},
}

@article {quintoriginal,
    AUTHOR = {Quint, J.-F.},
     TITLE = {Mesures de {P}atterson-{S}ullivan en rang sup\'erieur},
   JOURNAL = {Geom. Funct. Anal.},
  FJOURNAL = {Geometric and Functional Analysis},
    VOLUME = {12},
      YEAR = {2002},
    NUMBER = {4},
     PAGES = {776--809},
      ISSN = {1016-443X,1420-8970},
   MRCLASS = {22E40 (37C35)},
  MRNUMBER = {1935549},
MRREVIEWER = {Michel\ Coornaert},
       DOI = {10.1007/s00039-002-8266-4},
       URL = {https://doi-org.libproxy1.nus.edu.sg/10.1007/s00039-002-8266-4},
}

@article {sullivanoriginal,
    AUTHOR = {Sullivan, Dennis},
     TITLE = {The density at infinity of a discrete group of hyperbolic
              motions},
   JOURNAL = {Inst. Hautes \'Etudes Sci. Publ. Math.},
  FJOURNAL = {Institut des Hautes \'Etudes Scientifiques. Publications
              Math\'ematiques},
    NUMBER = {50},
      YEAR = {1979},
     PAGES = {171--202},
      ISSN = {0073-8301,1618-1913},
   MRCLASS = {58F17 (22E40 28C10 30C85)},
  MRNUMBER = {556586},
MRREVIEWER = {Troels\ J\o rgensen},
       URL = {http://www.numdam.org/item?id=PMIHES_1979__50__171_0},
}

@article {bgg1,
    AUTHOR = {Bradlow, Steven B. and Garc\'ia-Prada, Oscar and Gothen, Peter
              B.},
     TITLE = {Surface group representations and {${\rm U}(p,q)$}-{H}iggs
              bundles},
   JOURNAL = {J. Differential Geom.},
  FJOURNAL = {Journal of Differential Geometry},
    VOLUME = {64},
      YEAR = {2003},
    NUMBER = {1},
     PAGES = {111--170},
      ISSN = {0022-040X,1945-743X},
   MRCLASS = {53D30 (57R19)},
  MRNUMBER = {2015045},
MRREVIEWER = {Ignasi\ Mundet-Riera},
       URL = {http://projecteuclid.org.libproxy1.nus.edu.sg/euclid.jdg/1090426889},
}

@article {bgg2,
    AUTHOR = {Bradlow, Steven B. and Garc\'ia-Prada, Oscar and Gothen, Peter
              B.},
     TITLE = {Maximal surface group representations in isometry groups of
              classical {H}ermitian symmetric spaces},
   JOURNAL = {Geom. Dedicata},
  FJOURNAL = {Geometriae Dedicata},
    VOLUME = {122},
      YEAR = {2006},
     PAGES = {185--213},
      ISSN = {0046-5755,1572-9168},
   MRCLASS = {14D20 (53D30 58D29)},
  MRNUMBER = {2295550},
MRREVIEWER = {Anna\ Wienhard},
       DOI = {10.1007/s10711-007-9127-y},
       URL = {https://doi-org.libproxy1.nus.edu.sg/10.1007/s10711-007-9127-y},
}

@article {burger2008surfacegrouprepresentationsmaximal,
    AUTHOR = {Burger, Marc and Iozzi, Alessandra and Wienhard, Anna},
     TITLE = {Surface group representations with maximal {T}oledo invariant},
   JOURNAL = {Ann. of Math. (2)},
  FJOURNAL = {Annals of Mathematics. Second Series},
    VOLUME = {172},
      YEAR = {2010},
    NUMBER = {1},
     PAGES = {517--566},
      ISSN = {0003-486X,1939-8980},
   MRCLASS = {22E41 (20F67 57M07)},
  MRNUMBER = {2680425},
       DOI = {10.4007/annals.2010.172.517},
       URL = {https://doi-org.libproxy1.nus.edu.sg/10.4007/annals.2010.172.517},
}

@article {Goldman1988,
    AUTHOR = {Goldman, William M.},
     TITLE = {Topological components of spaces of representations},
   JOURNAL = {Invent. Math.},
  FJOURNAL = {Inventiones Mathematicae},
    VOLUME = {93},
      YEAR = {1988},
    NUMBER = {3},
     PAGES = {557--607},
      ISSN = {0020-9910,1432-1297},
   MRCLASS = {57M05 (22E40 32G15)},
  MRNUMBER = {952283},
MRREVIEWER = {William\ Harvey},
       DOI = {10.1007/BF01410200},
       URL = {https://doi-org.libproxy1.nus.edu.sg/10.1007/BF01410200},
}

@article {maximal,
    AUTHOR = {Burger, Marc and Iozzi, Alessandra and Labourie, Fran\c cois
              and Wienhard, Anna},
     TITLE = {Maximal representations of surface groups: symplectic {A}nosov
              structures},
   JOURNAL = {Pure Appl. Math. Q.},
  FJOURNAL = {Pure and Applied Mathematics Quarterly},
    VOLUME = {1},
      YEAR = {2005},
    NUMBER = {3},
     PAGES = {543--590},
      ISSN = {1558-8599,1558-8602},
   MRCLASS = {53C24 (53C35 57M99)},
  MRNUMBER = {2201327},
MRREVIEWER = {Jean-Fran\c cois\ Lafont},
       DOI = {10.4310/PAMQ.2005.v1.n3.a5},
       URL = {https://doi-org.libproxy1.nus.edu.sg/10.4310/PAMQ.2005.v1.n3.a5},
}

@article {Liplim,
    AUTHOR = {Pozzetti, Maria Beatrice and Sambarino, Andr\'es and Wienhard,
              Anna},
     TITLE = {Anosov representations with {L}ipschitz limit set},
   JOURNAL = {Geom. Topol.},
  FJOURNAL = {Geometry \& Topology},
    VOLUME = {27},
      YEAR = {2023},
    NUMBER = {8},
     PAGES = {3303--3360},
      ISSN = {1465-3060,1364-0380},
   MRCLASS = {22E40 (37C85 51F30)},
  MRNUMBER = {4668098},
MRREVIEWER = {Boris\ Hasselblatt},
       DOI = {10.2140/gt.2023.27.3303},
       URL = {https://doi-org.libproxy1.nus.edu.sg/10.2140/gt.2023.27.3303},
}

@article {GLW,
    AUTHOR = {Guichard, Olivier and Labourie, Fran\c cois and Wienhard,
              Anna},
     TITLE = {Positivity and representations of surface groups},
   JOURNAL = {Forum Math. Pi},
  FJOURNAL = {Forum of Mathematics. Pi},
    VOLUME = {14},
      YEAR = {2026},
     PAGES = {Paper No. e6, 38},
      ISSN = {2050-5086},
   MRCLASS = {22E40 (22E15 22F30 57S20)},
  MRNUMBER = {5029150},
       DOI = {10.1017/fmp.2025.10022},
       URL = {https://doi-org.libproxy1.nus.edu.sg/10.1017/fmp.2025.10022},
}

@article {HITCHIN1992449,
    AUTHOR = {Hitchin, N. J.},
     TITLE = {Lie groups and {T}eichm\"uller space},
   JOURNAL = {Topology},
  FJOURNAL = {Topology. An International Journal of Mathematics},
    VOLUME = {31},
      YEAR = {1992},
    NUMBER = {3},
     PAGES = {449--473},
      ISSN = {0040-9383},
   MRCLASS = {32G13 (32G15 57M99 58D27 58E15)},
  MRNUMBER = {1174252},
MRREVIEWER = {William\ Goldman},
       DOI = {10.1016/0040-9383(92)90044-I},
       URL = {https://doi-org.libproxy1.nus.edu.sg/10.1016/0040-9383(92)90044-I},
}

@article {Colin,
    AUTHOR = {Davalo, Colin},
     TITLE = {Maximal and {B}orel {A}nosov representations into {${\rm
              Sp}(4,\mathbb{R})$}},
   JOURNAL = {Adv. Math.},
  FJOURNAL = {Advances in Mathematics},
    VOLUME = {442},
      YEAR = {2024},
     PAGES = {Paper No. 109578, 21},
      ISSN = {0001-8708,1090-2082},
   MRCLASS = {22E40 (20F65 57K20)},
  MRNUMBER = {4711845},
MRREVIEWER = {Joan\ Porti},
       DOI = {10.1016/j.aim.2024.109578},
       URL = {https://doi.org/10.1016/j.aim.2024.109578},
}

@article {pozzetti2020conformalityrobustclassnonconformal,
    AUTHOR = {Pozzetti, Maria Beatrice and Sambarino, Andr\'es and Wienhard,
              Anna},
     TITLE = {Conformality for a robust class of non-conformal attractors},
   JOURNAL = {J. Reine Angew. Math.},
  FJOURNAL = {Journal f\"ur die Reine und Angewandte Mathematik. [Crelle's
              Journal]},
    VOLUME = {774},
      YEAR = {2021},
     PAGES = {1--51},
      ISSN = {0075-4102,1435-5345},
   MRCLASS = {22E40 (28A80 37C45 37F40)},
  MRNUMBER = {4250471},
MRREVIEWER = {Thomas\ Ward},
       DOI = {10.1515/crelle-2020-0029},
       URL = {https://doi-org.libproxy1.nus.edu.sg/10.1515/crelle-2020-0029},
}

@article {potrie2017eigenvaluesentropyhitchinrepresentation,
    AUTHOR = {Potrie, Rafael and Sambarino, Andr\'es},
     TITLE = {Eigenvalues and entropy of a {H}itchin representation},
   JOURNAL = {Invent. Math.},
  FJOURNAL = {Inventiones Mathematicae},
    VOLUME = {209},
      YEAR = {2017},
    NUMBER = {3},
     PAGES = {885--925},
      ISSN = {0020-9910,1432-1297},
   MRCLASS = {57M50 (37D40)},
  MRNUMBER = {3681396},
MRREVIEWER = {A.\ N.\ Bezen},
       DOI = {10.1007/s00222-017-0721-9},
       URL = {https://doi-org.libproxy1.nus.edu.sg/10.1007/s00222-017-0721-9},
}

@article {Benoistasym,
    AUTHOR = {Benoist, Y.},
     TITLE = {Propri\'et\'es asymptotiques des groupes lin\'eaires},
   JOURNAL = {Geom. Funct. Anal.},
  FJOURNAL = {Geometric and Functional Analysis},
    VOLUME = {7},
      YEAR = {1997},
    NUMBER = {1},
     PAGES = {1--47},
      ISSN = {1016-443X,1420-8970},
   MRCLASS = {22E15},
  MRNUMBER = {1437472},
MRREVIEWER = {Scot\ Adams},
       DOI = {10.1007/PL00001613},
       URL = {https://doi-org.libproxy1.nus.edu.sg/10.1007/PL00001613},
}

@article {rk2maximal,
    AUTHOR = {Collier, Brian and Tholozan, Nicolas and Toulisse, J\'er\'emy},
     TITLE = {The geometry of maximal representations of surface groups into
              {${\rm SO}_0(2,n)$}},
   JOURNAL = {Duke Math. J.},
  FJOURNAL = {Duke Mathematical Journal},
    VOLUME = {168},
      YEAR = {2019},
    NUMBER = {15},
     PAGES = {2873--2949},
      ISSN = {0012-7094,1547-7398},
   MRCLASS = {20H10 (14H60 22E40 53C50 57M50 58E12)},
  MRNUMBER = {4017517},
MRREVIEWER = {John\ Huerta},
       DOI = {10.1215/00127094-2019-0052},
       URL = {https://doi-org.libproxy1.nus.edu.sg/10.1215/00127094-2019-0052},
}

@article {shadowvitali,
    AUTHOR = {Roblin, Thomas},
     TITLE = {Ergodicit\'e{} et \'equidistribution en courbure n\'egative},
   JOURNAL = {M\'em. Soc. Math. Fr. (N.S.)},
  FJOURNAL = {M\'emoires de la Soci\'et\'e{} Math\'ematique de France.
              Nouvelle S\'erie},
    NUMBER = {95},
      YEAR = {2003},
     PAGES = {vi+96},
      ISSN = {0249-633X,2275-3230},
   MRCLASS = {37D40 (28D10 30F40 37A25 37C40 53C22 53D25 57R30)},
  MRNUMBER = {2057305},
MRREVIEWER = {Marc\ Peign\'e},
       DOI = {10.24033/msmf.408},
       URL = {https://doi-org.libproxy1.nus.edu.sg/10.24033/msmf.408},
}

@article {kapovich2017anosov,
    AUTHOR = {Kapovich, Michael and Leeb, Bernhard and Porti, Joan},
     TITLE = {Anosov subgroups: dynamical and geometric characterizations},
   JOURNAL = {Eur. J. Math.},
  FJOURNAL = {European Journal of Mathematics},
    VOLUME = {3},
      YEAR = {2017},
    NUMBER = {4},
     PAGES = {808--898},
      ISSN = {2199-675X,2199-6768},
   MRCLASS = {22E40 (20F65 53C35)},
  MRNUMBER = {3736790},
MRREVIEWER = {Herbert\ Abels},
       DOI = {10.1007/s40879-017-0192-y},
       URL = {https://doi-org.libproxy1.nus.edu.sg/10.1007/s40879-017-0192-y},
}

@book{steinreal,
 author = {Stein, Elias M. and Shakarchi, Rami},
 title = {{Real Analysis. Measure Theory, Integration, and {Hilbert} Spaces}},
 fseries = {Princeton Lectures in Analysis},
 series = {Princeton Lectures in Analysis},
 volume = {3},
 isbn = {0-691-11386-6},
 year = {2005},
 publisher = {Princeton, NJ: Princeton University Press},
 zbMATH = {2171466},
 Zbl = {1081.28001}
}

\end{document}